\documentclass[12pt, reqno]{amsart}

\usepackage{graphicx} 
\usepackage{amsmath}
\usepackage{amssymb}
\usepackage{amsfonts}
\usepackage{mathtools}
\usepackage{amsthm}  
\usepackage{enumerate}  
\usepackage{dsfont}
\usepackage[colorlinks]{hyperref}
\usepackage{mathrsfs}
\usepackage{xcolor}
\usepackage{anysize}
\usepackage{mathabx}
\usepackage{marginnote}
\usepackage{esint} 
\usepackage{multicol}

\marginsize{3cm}{3cm}{2cm}{2cm}

\usepackage{bm}

\newcommand{\R}{\mathbb{R}}

\newcommand{\Q}{\mathcal{Q}}
\newcommand{\Z}{\mathbb{Z}}
\newcommand{\N}{\mathbb{N}}
\newcommand{\B}{\mathbb{B}}
\newcommand{\T}{\mathbb{T}}

\newcommand{\bbD}{\mathbb{D}}

\newcommand{\tH}{\mathrm{H}}

\newcommand{\tS}{\mathrm{S}}
\newcommand{\tT}{\mathrm{T}}

\newcommand{\fkp}{\mathfrak{p}}

\newcommand{\fkD}{\mathfrak{D}}

\newcommand{\fkP}{\mathfrak{P}}

\newcommand{\fkT}{\mathfrak{T}}

\newcommand{\bdB}{\mathbf{B}}
\newcommand{\bdC}{\mathbf{C}}

\newcommand{\bdE}{\mathbf{E}}

\newcommand{\bdH}{\mathbf{H}}

\newcommand{\bdK}{\mathbf{K}}

\newcommand{\bdQ}{\mathbf{Q}}
\newcommand{\bdR}{\mathbf{R}}

\newcommand{\bdT}{\mathbf{T}}

\newcommand{\bdk}{\mathbf{k}}

\newcommand{\bdx}{\mathbf{x}}

\newcommand{\bdz}{\mathbf{z}}

\newcommand{\bdbeta}{\boldsymbol{\beta}}

\newcommand{\wt}[1]{\widetilde{#1}}

\newcommand{\diam}[1]{\operatorname{diam}\left(#1\right) }
\newcommand{\supp}{\operatorname{supp}}
\newcommand{\dist}{\mathrm{dist}}

\newcommand{\D}{\mathbb{D}}

\newcommand{\ud}{\mathrm{d}}

\def\tnorm#1{\vvvert #1 \vvvert}

\newcommand{\cI}{\mathcal{I}}

\newcommand{\cP}{\mathcal{P}}
\newcommand{\cQ}{\mathcal{Q}}

\newcommand{\cS}{\mathcal{S}}

\newcommand{\cZ}{\mathcal{Z}}

\newcommand{\avg}{\mathrm{avg}}
\newcommand{\hyp}{\mathrm{hyp}}
\newcommand{\mix}{\mathrm{mix}}
\newcommand{\elliptic}{\mathrm{ell}}

\def\inn#1#2{\langle#1,#2\rangle}

\theoremstyle{plain}
\newtheorem{theorem}{Theorem}[section]

\newtheorem*{theorem*}{Theorem}
\newtheorem{proposition}[theorem]{Proposition}
\newtheorem{corollary}[theorem]{Corollary}
\newtheorem{lemma}[theorem]{Lemma}

\theoremstyle{definition} 
\newtheorem{definition}[theorem]{Definition}

\numberwithin{equation}{section}

\allowdisplaybreaks
\title[\texorpdfstring{Oscillatory integral operators and variable propagators}{}]{Oscillatory integral operators and variable Schr\"odinger propagators of arbitrary signature: beyond the universal estimates}

\author[F. Ballesta-Yag\"ue]{Fernando Ballesta-Yag\"ue}
\address[F. Ballesta-Yag\"ue]{Departamento de An\'alisis Matem\'atico y Matem\'atica Aplicada, Facultad de Ciencias Matem\'aticas, Universidad Complutense de Madrid \hfill\break\indent Pl. de las Ciencias 3, 28040 Madrid, Spain}
\email{ferballe@ucm.es}

\author[J. Hickman]{Jonathan Hickman}
\address[J. Hickman]{School of Mathematics and Maxwell Institute for Mathematical Sciences, James
Clerk Maxwell Building, The King’s Buildings, Peter Guthrie Tait Road, Edinburgh,
EH9 3FD, UK}
\email{jonathan.hickman@ed.ac.uk}

\author[H. Roy]{Hrit Roy}
\address[H. Roy]{Department of Mathematics and Statistics, University of Jyväskylä \hfill\break\indent PO Box 35
FI-40014 University of Jyväskylä,
Finland}
\email{hrit.h.roy@jyu.fi}

\newcommand{\eps}{\varepsilon}

\begin{document}

\begin{abstract}
    We consider H\"ormander oscillatory integral operators on $\R^n$. In particular, we are interested in proving bounds with dependence on the signature of the phase. The sharp universal range of $L^p\to L^p$ estimates was obtained for any given signature by the second author and Iliopoulou.  Under the mild non-compression hypotheses on the phase introduced in Chen \textit{et al.}, we prove $L^p\to L^p$ estimates that go beyond this universal range in odd dimensions. We also obtain the sharp universal range of $L^p\to L^q$ estimates for odd $n$. The techniques are based on the breakthrough work of Wang--Wu on the Fourier restriction conjecture. 
\end{abstract}

\thanks{
F. B-Y. is supported by grants PID2024-155917NB-I00 funded by MCIN/AEI/10.13039/501100011033,  CEX2019-000904-S funded by MCIN/AEI/ 10.13039/501100011033 and Grupo UCM-970966 (Spain), and benefited from an FPU Grant FPU21/06111 and FPU Travelling Grant EST25/00120 from Ministerio de Universidades (Spain). J. H. is supported by New Investigator Award UKRI097. H. R. is supported by the European Research Council (ERC) under the European Union’s Horizon Europe research and innovation programme (grant agreement No 101087499) and by the Research Council of Finland via the project \emph{Approximate incidence geometry}, grant no. 355453.
}

\title[\texorpdfstring{Oscillatory integral operators}{}]{Oscillatory integral operators of arbitrary signature: beyond the universal estimates}

\maketitle




\section{Introduction}




\subsection{Oscillatory integral operators}

This article is a follow up to the recent work \cite{CGGHIW} and also to \cite{HI2022}. Our primary goal is to study the $L^p$ mapping properties of certain variable coefficient oscillatory integral operators. We begin by recalling the relevant definitions. 

For $n \geq 2$ and $\rho > 0$ write $\D^n_{\rho}\coloneq \B_{\rho}^{n-1} \times \B^1_{\rho} \times \B_{\rho}^{n-1}$, where  $\B^d_{\rho}$ denotes the closed ball in $\R^d$ of radius $\rho$ centred at the origin. When $\rho = 1$, we abbreviate $\D^d_{\rho}$ and $\B^d_{\rho}$ to $\D^d$ and $\B^d$, respectively. 

\begin{definition}
For $n \geq 2$, let $\phi \in C^{\infty}(\D^n)$ be a real-valued function. We say $\phi$ is a \textit{H\"ormander-type phase} if it satisfies the following conditions: 
\begin{enumerate}
    \item [H1)] $\det \partial_{y x}^2 \phi(x, t ; y)\neq 0$ for all $(x, t ; y) \in \D^n$. 
    \item [H2)] Defining the map $G: \D^n \to \mathbb{S}^{n-1}$ by $G(x,t; y)\coloneqq \frac{G_0(x,t ; y)}{\left|G_0(x,t ; y)\right|}$ where
    \[
    G_0(x,t ; y)\coloneqq\bigwedge_{j=1}^{n-1} \partial_{y_j} \partial_{x,t} \phi(x,t ; y),
    \]
    the curvature condition
    \begin{equation}\label{eq:curvatureCondition}
    \left.\operatorname{det} \partial_{y y}^2\langle\partial_{x,t} \phi(x, t ; y), G(x,t ; y_0)\rangle\right|_{y=y_0} \neq 0    
    \end{equation}
    holds for all $(x,t ; y_0) \in \D^n$.    
\end{enumerate}
\end{definition}

For the purposes of this article, it is important to refine H2) by specifying the \textit{signature} of the matrix appearing in \eqref{eq:curvatureCondition}. Given a symmetric matrix $A \in \mathrm{GL}(n-1, \R)$ we define its signature to be $|\sigma_+ - \sigma_-|$  where $\sigma_+$ and $\sigma_-$ are the numbers of positive and negative eigenvalues, respectively. 

\begin{definition}
Let $\phi$ be a H\"ormander-type phase. The \textit{signature} of $\phi$ is defined to be the common signature of the matrix
\begin{equation*}
\left.\partial_{y y}^2\langle\partial_{x,t} \phi(x, t ; y), G(x, t ; y_0)\rangle\right|_{y=y_0} \qquad \textrm{for any $(x,t; y_0) \in \D^n$.}
\end{equation*}
\end{definition}

Note that, by continuity, the signature is constant over the domain $\D^n$ and so the signature of $\phi$ is well-defined as a single integer value $0 \leq \sigma \leq n-1$.

Often we will pair the phase with a smooth cutoff function $a \in C_c^{\infty}(\D^n)$, in which case we refer to the pairing $[\phi; a]$ as a \textit{phase-amplitude pair}. 

It is convenient to work with a rescaled setup. For any $\lambda \geq 1$ and phase-amplitude pair $[\phi; a]$, let 
\begin{equation}\label{eq: rescaled}
  \phi^{\lambda}(x,t;y) \coloneq\lambda\phi(x/\lambda, t/\lambda;y) \qquad \textrm{and} \qquad   a^{\lambda}(x, t; y) \coloneq a(x/\lambda, t/\lambda; y).
\end{equation}
\begin{definition} Let $[\phi; a]$ be a phase-amplitude pair and $\lambda \geq 1$.
The \textbf{H\"ormander-type oscillatory integral operator} $H^{\lambda}[\phi; a]$ associated to $[\phi; a]$ is given by
\begin{equation*}
H^{\lambda}[\phi; a]f(x,t) \coloneq  \int_{\B^{n-1}} e^{i \phi^{\lambda}(x,t; y)}a^{\lambda}(x,t;y) f(y)\,\ud y,
\end{equation*}
defined initially as acting on $f \in C^{\infty}_c(\B^{n-1})$.
\end{definition}

H\"ormander-type operators $H^{\lambda}[\phi; a]$ were introduced by H\"ormander~\cite{hormander} as a means to unify the Carleson--Sj\"olin \cite{CarlesonSjolin} approach to the Bochner--Riesz problem and the nascent Fourier restriction theory of Fefferman and Stein \cite{Fefferman1970}. These operators have found a wide range of applications to problems in harmonic analysis and beyond. Here we are interested in determining the $L^p$ mapping properties for the operators. This question has a long and rich history and has seen an intense burst of activity in recent years. For additional context, we refer the reader to the introductory sections of \cite{GHI, CGGHIW, GWZ2024} and also \cite{Nadjimzadah}.

In order to state and compare results, we introduce the following notation.

\begin{definition}
Let $\phi: \D^n \to \R$ be a H\"ormander-type phase, $1 \leq p, q \leq \infty$ and $\alpha \geq 0$. We let $\bdH_{p \to q}(\phi ; \alpha)$ denote the statement: there exists $\rho=\rho(\phi)>0$ such that for all choices of $a \in C_c^{\infty}\left(\D_\rho^n\right)$ and all $\varepsilon > 0$, the inequality\footnote{See \S\ref{sec: notation} below for the definition of the asymptotic notation.} 
    \begin{equation}\label{eq: H pq}
    \|H^\lambda[\phi ; a] f\|_{L^q(B_R)} \lesssim_{\phi, a, \varepsilon}
    R^{\alpha + \varepsilon}
    \|f\|_{L^p(\B^{n-1})}    
    \end{equation}
    holds uniformly over $1 \leq R \leq \lambda$, all $R$-balls $B_R \subseteq \R^n$ and all $f \in C_c^{\infty}(\B^{n-1})$. When $\alpha=0$, we write $\bdH_{p \to q}(\phi)$ in place of $\bdH_{p \to q}(\phi ; 0)$.
\end{definition}




\subsection{Universal estimates} For $n \geq 2$ and $0 \leq \sigma \leq n-1$ define the \textit{universal exponent}
\begin{equation*}
q_{n, \sigma} \coloneq
\begin{cases}
2 \cdot \frac{\sigma+2(n+1)}{\sigma+2(n-1)} &\text { if } n \text { is odd,} \\[4pt]
2 \cdot \frac{\sigma+2 n+3}{\sigma+2 n-1} &\text { if } n \text { is even.}       
\end{cases}
\end{equation*}  
With this definition, the following theorem describes the \textit{universal} range of estimates, satisfied for any choice of phases with a prescribed signature.

\begin{theorem}[\cite{HI2022}]\label{thm: HI2022} For $n \geq 2$ and $0 \leq \sigma \leq n-1$, let $\phi \colon \D^n \to \R$ be a H\"ormander-type phase with signature $\sigma$. Then $\bdH_{q\to q}(\phi)$ holds for all $q\geq q_{n,\sigma}$.
\end{theorem}

For $n = 2$, Theorem~\ref{thm: HI2022} coincides with H\"ormander's theorem \cite{hormander}. Moreover, prior to \cite{HI2022}, the minimal signature case where $n$ is odd and $\sigma = 0$ was treated by Stein~\cite{stein_beijing} and where $n$ is even and $\sigma = 1$ by Bourgain--Guth \cite{bourgain--guth}. On the other extreme, the positive-definite case (where $n \geq 2$ is arbitrary and $\sigma = n-1$) was treated in \cite{GHI}. 

Theorem~\ref{thm: HI2022} is known to be sharp in the following sense: for any permissible choice of dimension $n \geq 2$ and signature $0 \leq \sigma \leq n-1$, there exists a phase $\phi_{\tH}$ with signature $\sigma$ such that $\bdH_{p\to q}(\phi_{\tH})$ fails for all $1 \leq p \leq \infty$ whenever $q < q_{n, \sigma}$. This was established in \cite{HI2022} by adapting a well-known example due to Bourgain~\cite{BourgainGAFA1991_no4}.

The sharp examples for $\bdH_{q\to q}(\phi)$ are based on the \textit{Kakeya compression phenomenon}. In particular, the sharp examples all have the property that the wavepackets associated to the operators can be arranged so as to concentrate in a small neighbourhood of a low dimensional, low degree variety.




\subsection{Geometric maximal functions}\label{subsec: geometric max}  Let $\phi \colon \D^n \to \R$ be a H\"ormander-type phase and, fixing a ball $B= B(\bdx_B, R) \subseteq \R^n$ for $R \geq 1$, define the recentred phase
\begin{equation*}
    \phi^\lambda_B(\bdx;y) \coloneq \phi^\lambda(\bdx_B + \bdx;y) - \phi^\lambda(\bdx_B;y).
\end{equation*}
Given $y \in R^{-1/2}\Z^{n-1} \cap \B^{n-1}$ and $\omega \in R^{1/2}\Z^{n-1}$, define the associated \textit{$R$-tube} by
\begin{equation}\label{eq: original tube}
    T_{y,\omega} \coloneq \big\{ \bdx = (x,t) \in \R^{n-1} \times (-2R, 2R) : |\partial_y \phi_B^{\lambda}(\bdx;y) - \omega| < R^{1/2} \big\}.
\end{equation}
We let $\bdT^R$ denote the set of all such $R$-tubes. Note that $\bdT^R$ depends on a choice of $R$-ball $B$, but we suppress this in the notation. 

\begin{definition}Given a set of $R$-tubes $\T \subseteq \bdT^R$, we say $\T$ is \textit{direction-separated} if for every $y \in R^{-1/2}\Z^{n-1} \cap \B^{n-1}$ there exists at most one $\omega \in R^{1/2}\Z^{n-1}$ such that $T_{y, \omega} \in \T$.
\end{definition}

The key geometric/combinatorial inputs are estimates of the form
\begin{equation}\label{eq: multiplicity bound}
      \big\|\sum_{T \in \T} \chi_T\big\|_{L^s(\R^n)} \lesssim_{\phi, \varepsilon} R^{(n+1)/(2s)+\beta/2 + \varepsilon} \big[\#\T \big]^{1/s},
    \end{equation}
where $\T$ ranges over certain families of tubes associated to some phase $\phi$ and $\beta$, $\varepsilon \geq 0$. To make this more precise, we introduce the following definitions.

\begin{definition} Let $\phi \colon \D^n \to \R$ be a H\"ormander-type phase, $1 \leq s \leq \infty$ and $\beta \geq 0$. We let $\bdK_{s \to s}^*(\phi; \beta)$ denote the statement:\medskip
    
    \noindent For all $\varepsilon > 0$, the inequality
    \eqref{eq: multiplicity bound} holds uniformly over all $\lambda\geq 1$, all $1\leq R \leq \lambda$, and all direction-separated families $\T$ of $R$-tubes associated to $\phi$. 
\end{definition}

We have the following \textit{universal} geometric estimates, which hold for all H\"ormander-type phases $\phi$.

\begin{proposition}[Universal geometric estimate]\label{prop: univ Kak/Nik} Let $n \geq 3$ be odd and $\phi \colon \D^n \to \R$ be a H\"ormander-type phase. Then $\bdK_{s \to s}^*\big(\phi;  \beta(n, s)\big)$ holds for all $1 \leq s \leq \frac{n+1}{n-1}$ and $\beta(n, s) \coloneq \frac{n-1}{2s'}$.
\end{proposition}

Proposition~\ref{prop: univ Kak/Nik} can be proved using a variant of Bourgain's bush argument: see \cite{wisewell} (and also Lemma~\ref{lem: univ two-ends} below). For $n \geq 2$ even, Proposition~\ref{prop: univ Kak/Nik} is still valid, but stronger universal estimates in fact hold: see \cite{bourgain--guth}.



 
\subsection{Main results} Although Theorem~\ref{thm: HI2022} is sharp over the class of H\"ormander-type phases of fixed signature, there are specific examples where one can do better. Our main result is a conditional statement which shows that geometric inequalities of the form $\bdK_{s \to s}^*(\phi; \beta)$ for certain $s$ and $\beta$ imply corresponding oscillatory operator bounds $\bdH_{p \to q}(\phi)$. This conditional result can be combined with the geometric maximal estimates from \cite{CGGHIW} to obtain improved estimates for a wide class of oscillatory operators. 

Unfortunately, the numerology in the main result is somewhat involved. For this reason, we postpone the definitions of the exponents involved to \S\ref{subsec: numerology} below. This will allow us to proceed directly to the corollaries of our theorem, which are of primary interest.

\begin{theorem}\label{thm: main}
Let $n \geq 3$ be odd and $\phi \colon \D^n \to \R$ be a H\"ormander-type phase of signature $\sigma$. Further let $1 \leq s \leq \frac{n+1}{n-1}$ and $\beta \geq 0$ and define $q \coloneq q_{n, \sigma}(s, \beta)$ and $p\coloneq p_{n,\sigma}(s,\beta)$ as in \S\ref{subsec: numerology} below. Then
$\bdK_{s \to s}^*(\phi; \beta) \implies \bdH_{p \to q}(\phi).$
\end{theorem}

As a first application, we obtain sharp universal off-diagonal bounds in odd dimensions. These bounds appear to be new at this level of generality, and recover the universal oscillatory estimates of Theorem~\ref{thm: HI2022} in odd dimensions. 

\begin{corollary}\label{cor: univ off diag}
    Let $n \geq 3$ be odd and $\phi \colon \D^n \to \R$ be a H\"ormander-type phase of signature $\sigma$. Then $\bdH_{p \to q}(\phi)$ holds for all $1 \leq p$, $q \leq \infty$ satisfying $q\geq q_{n,\sigma}$ and
\begin{equation}\label{eq: universal range off diagonal}
        \frac{n+1}{q}\leq\frac{n-\sigma-1}{2}+\frac{\sigma}{p'},\qquad 
        p'\leq \frac{n-1}{n+1}q.
    \end{equation}
\end{corollary}

Corollary~\ref{cor: univ off diag} is a very direct consequence of Theorem~\ref{thm: main} and we shall give the proof later in the introduction, in \S\ref{subsec: numerology} below. One interesting feature of these bounds is that they are sharp. In particular, there exists a phase $\phi$ satisfying the hypotheses of Corollary~\ref{cor: univ off diag} such that $\bdH_{p \to q}(\phi)$ fails for $p$, $q$ outside the range of \eqref{eq: universal range off diagonal}. See Lemma~\ref{lem: off diag sharp}.

Our second application is to obtain estimates beyond the universal range in Theorem \ref{thm: HI2022} for a certain class of phases.

In \cite{BourgainGAFA1991_no4}, Bourgain sought weak additional hypotheses on $\phi$ which guaranteed estimates beyond the universal range of Theorem~\ref{thm: HI2022} in the $n =3$ and $\sigma = 0$ case. This line of enquiry was further pursued in \cite{CGGHIW} where, \textit{inter alia}, the authors extended Bourgain's results to general odd dimension $n \geq 3$, again in the $\sigma = 0$ case. The key is to obtain conditions on the phase which prevent (or at least partially prevent) the Kakeya compression phenomenon observed in the sharp examples. In \cite{CGGHIW}, such conditions were found and were referred to as
\begin{enumerate}
    \item[]\textit{Hypothesis I): Kakeya non-compression} (\cite[Definition 2.6]{CGGHIW}).
\end{enumerate}
Since the statement of this hypothesis is somewhat involved, we refer the reader to \cite{CGGHIW} for the details.

Both \cite{BourgainGAFA1991_no4} and \cite{CGGHIW} work under the additional assumption that the phase is \textit{translation invariant}. That is, we assume there exists some $\psi \in C^{\infty}(\B^1 \times \B^{n-1})$ such that
    \begin{equation*}
    \phi(x,t;y)=\langle x,y\rangle+\psi(t;y)\qquad \textrm{where} \qquad \psi(0;y)=0.
\end{equation*}

By combining Theorem~\ref{thm: main} with the geometric estimates from \cite[Theorem 2.9]{CGGHIW}, we obtain the following corollary. 

\begin{corollary}\label{cor: beyond univ} Let $n \geq 3$ be odd and $\phi \colon \D^n \to \R$ be a real analytic H\"ormander-type phase of signature $\sigma$. Suppose $\phi$ is translation-invariant and satisfies Hypothesis I). Then there exists $q_\phi<q_{n,\sigma}$ such that $\bdH_{q \to q}(\phi)$ holds for all $q \geq q_\phi$.
\end{corollary}

We show how Corollary~\ref{cor: beyond univ} follows from Theorem~\ref{thm: main} and \cite[Theorem 2.9]{CGGHIW} in \S\ref{subsec: cor} below. We remark that this argument involves some significant additional technicalities owing to the required uniformity in $R$ in the estimate $\bdK_{s \to s}^*(\phi; \beta)$. This uniformity does not feature in the geometric estimates from \cite{CGGHIW} and forces the use of additional arguments. 

Corollary~\ref{cor: beyond univ} improves over \cite[Theorem~2.8]{CGGHIW} in the following way:
Theorem~2.8 in \cite{CGGHIW} only leads to an improved range of $q$ in the minimal signature $\sigma = 0$ case. On the other hand, Corollary~\ref{cor: beyond univ} is effective for every choice of $\sigma$. 

We remark that the case $n = 3$ and $\sigma = 2$ is also considered in \cite[Theorem~2.5]{CGGHIW}, although the precise details of the proof are not given. It is unclear whether the methods used to prove this specific case in \cite{CGGHIW} can be generalised. All remaining cases of Corollary~\ref{cor: beyond univ} (that is, the cases $\sigma>0$, other than $n=3$, $\sigma=2$) are new.

The value of $q_{\phi}$ in Corollary~\ref{cor: beyond univ} is not explicit, since the exponents featured in the geometric estimates from \cite[Theorem 2.9]{CGGHIW} are not explicit. Following Bourgain~\cite{BourgainGAFA1991_no4}, the focus here is on finding weak conditions which guarantee \textit{some} improvement over the universal range. Hypothesis I) is indeed a weak condition, which is satisfied by many example operators. The power and interest of Corollary~\ref{cor: beyond univ} lies in this range and flexibility of application. Furthermore, there are important problems for which it is crucial only to know \textit{some} improvement holds over the universal estimate, without needing to quantify the precise value of $q_{\phi}$. One striking instance of this is the recent resolution in \cite{BGH_PY} of a question posed by Pierce--Yung \cite{PY2019} concerning maximally modulated singular Radon transforms.




\subsection{Proof strategy}

To prove Theorem~\ref{thm: main}, we rely on a method introduced by Wang--Wu~\cite{Wang_Wu} in their breakthrough work on the Fourier restriction problem (see also \cite{MKN, GHLW}). The approach of \cite{Wang_Wu}, adpated to our purpose, cleanly combines refined decoupling with certain key geometric maximal bounds from \cite{CGGHIW} to produce the desired result.

The Wang--Wu~\cite{Wang_Wu} method involves two main ingredients. Here we describe how these ingredients manifest and are treated in the proof of Theorem~\ref{thm: main}. 

\begin{enumerate}
    \item \textit{Refined decoupling}. The oscillatory parts of the problem are handled using \textit{refined decoupling inequalities}. Such inequalities were introduced in \cite{GIOW2018} (and observed independently by Du--Zhang) for constant coefficient, positive-definite phases, as part of work on the Falconer distance problem. It was later extended to variable coefficient positive-definite phases in \cite{ILX}. To prove Theorem~\ref{thm: main}, we formulate a refined decoupling inequality in the mixed signature setting. The desired refined decoupling inequality then follows from a corresponding decoupling statement due to Schippa~\cite{Schippa}.
    
    \item \textit{Two-ends Kakeya maximal estimates}. The second ingredient is a geometric Kakeya-type bound. For the Wang--Wu \cite{Wang_Wu} method, $\bdK_{s \to s}^*(\phi; \beta)$ alone is insufficient for our purposes. Instead, we interpolate these bounds with a refined bush estimate, which is valid under a two-ends hypothesis in the sense of Wolff~\cite{Wolff1995} (see Lemma~\ref{lem: univ two-ends} below). 
\end{enumerate}

Rather than use the Wang--Wu~\cite{Wang_Wu} method, the earlier work \cite{CGGHIW} relies on a classical argument due to Bourgain \cite{BourgainGAFA1991_no4}. This approach does not apply in the mixed signature case. Indeed, Bourgain's method \cite{BourgainGAFA1991_no4} is heavily tied to the Stein--Tomas exponent $2 \cdot \frac{n+1}{n-1}$ and for signature $\sigma>0$ we need to work away from this. 




\subsection{Numerology and proof of Corollary~\ref{cor: univ off diag}}\label{subsec: numerology} To conclude the introduction, we define the exponents featured in Theorem~\ref{thm: main}. 

Let $n \geq 3$ be odd and $0 \leq \sigma \leq n-1$. Define
\begin{equation}\label{eq: dec exp 0}
r_{n,\sigma}\coloneq\frac{2(n+\sigma+3)}{n+\sigma-1}  \quad \textrm{and} \quad t_{n,\sigma}\coloneq\frac{n+\sigma+3}{2n+\sigma+2} \quad \textrm{so} \quad \frac{1}{q_{n,\sigma}}=\frac{1-t_{n,\sigma}}{2}+\frac{t_{n,\sigma}}{r_{n,\sigma}}.
\end{equation}
The exponent $r_{n, \sigma}$ is the decoupling exponent for Hörmander-type phases of signature $\sigma$ (see Theorem \ref{thm: dec schippa} below). Given $1 \leq s \leq \frac{n+1}{n-1}$, let $\beta(n, s) \coloneq \frac{n-1}{2s'}$ as in Proposition~\ref{prop: univ Kak/Nik}. To simplify certain formul\ae, for $\beta \geq 0$, we define
\begin{equation}\label{eq: A exp}
    A(s, \beta) \coloneq 
    (s-1)\Big(\frac{n-1}{4}-\frac{\beta s'}{2}\Big)
    =
    \frac{s}{2}\big(\beta(n,s)-\beta\big).
\end{equation}

The $L^p\to L^q$-bound will be obtained by interpolating an $L^2\to L^2$-estimate with an $L^{a}\to L^{r_{n,\sigma}}$-estimate. Hence, the exponents have the form 
\begin{align}
\label{eq: target exp arrival}
\frac{1}{q_{n,\sigma}(s,\beta)}
&=
\frac{1-t_{n,\sigma}(s,\beta)}{2}+\frac{t_{n,\sigma}(s,\beta)}{r_{n,\sigma}},
\quad \\
\label{eq: target exp domain}
\frac{1}{p_{n,\sigma}(s,\beta)}
&=
\frac{1-t_{n,\sigma}(s,\beta)}{2}+\frac{t_{n,\sigma}(s,\beta)}{a_{n,\sigma}(s,\beta)},
\end{align}
where the interpolation parameter is chosen to be
\begin{equation}\label{eq: u exp}
    t_{n,\sigma}(s,\beta) \coloneq t_{n,\sigma}\cdot u_{n,\sigma}(s,\beta),\quad \text{with}\quad u_{n, \sigma}(s,\beta) \coloneq  \frac{1-A(s,\beta)}{1-A(s,\beta)t_{n,\sigma}},
\end{equation}
and $a = a_{n,\sigma}(s,\beta)$ is given by 
\begin{equation}\label{eq: a n sigma p beta}
    a_{n,\sigma}(s,\beta)\coloneq\frac{2(n+\sigma+3)}{n+\sigma+1-2B(s,\beta)} \quad
\textrm{where} \quad
    B(s,\beta)\coloneq\frac{s-1}{4}\frac{n-1}{1-A(s,\beta)}.
\end{equation}

Having defined the relevant exponents, we may now derive Corollary~\ref{cor: univ off diag} from Theorem~\ref{thm: main}. First, we recall Stein's estimate \cite{stein_beijing}, which is $\bdH_{2\to q_{\mathrm{ST}}}(\phi),$ where $q_{\mathrm{ST}}\coloneq\tfrac{2(n+1)}{n-1}=q_{n,0}$. Define the exponent $p_{n,\sigma}\coloneq\tfrac{2(2n+\sigma+2)}{2n+\sigma}$.

\begin{proof}[Proof of Corollary \ref{cor: univ off diag}]
    Take $s=1$ and $\beta=\beta(n,1)=0$. For these choices, $q_{n,\sigma}(s,\beta)=q_{n,\sigma}$ and $p_{n,\sigma}(s,\beta)=p_{n,\sigma}$. 
    
    By Theorem \ref{thm: main} and the trivial $L^1$ estimate $\mathbf{K}^*_{1\to 1}(\phi;0)$, we obtain $\bdH_{p_{n,\sigma}\to q_{n,\sigma}}(\phi)$. When $\sigma=0$, Stein's estimate gives $\bdH_{2\to q_{\tS\tT}}(\phi)$ (which is stronger than $\bdH_{p_{n,0}\to q_{n,0}}(\phi)$). For  $\sigma>0$, we interpolate $\bdH_{p_{n,\sigma}\to q_{n,\sigma}}(\phi)$ with Stein's estimate $\bdH_{2\to q_{\mathrm{ST}}}(\phi)$, and pass to larger $p$ using the nesting of $L^p(\B^{n-1})$. This gives 
     $\bdH_{p\to q}(\phi)$ for all $$q_{n,\sigma}\leq q\leq q_{\mathrm{ST}},\qquad \frac{n+1}{q}\leq\frac{n-\sigma-1}{2}+\frac{\sigma}{p'}.$$
    On the other hand, interpolating Stein's estimate with the trivial estimate $\bdH_{1\to \infty}(\phi)$, and passing to larger $p$ gives $\bdH_{p\to q}(\phi)$ for all $$q\geq q_{\mathrm{ST}},\qquad p'\leq\frac{n-1}{n+1}q,$$
    concluding the proof.
\end{proof}

\subsection*{AI Declaration} Very minor use of AI was made for final proofreading of the article.   



\section{Notation and general preliminaries}\label{sec: notation}

\subsubsection*{Fourier transform} Given $f \in L^1(\R^{n-1})$, we define
\begin{equation}\label{eq: FT}
    \widehat{f}(\xi)  \coloneq \int_{\R^{n-1}} f(x) e^{-i \xi \cdot x}\,\ud x \quad \textrm{and} \quad \widecheck{f}(x)  \coloneq \frac{1}{(2\pi)^{n-1}}\int_{\R^{n-1}} f(\xi) e^{i x \cdot \xi}\,\ud \xi.
\end{equation}

\subsubsection*{Cubes} For $\rho > 0$, a \textit{lattice $\rho$-cube} is a set $Q \subseteq \R^n$ of the form $\rho\bdk +[-\rho/2, \rho/2)^n$ for some $\bdk \in \Z^n$. We let $\bdQ^{\rho}$ denote the set of all lattice $\rho$-cubes. Given $\cQ \subseteq \bdQ^{\rho}$, we let $\bdE(\cQ) \coloneq \bigcup_{Q \in \cQ} Q$ denote their union. 

\subsubsection*{Caps} Let $c_n \coloneq 1/\sqrt{n}$.  For $0 < r \leq 1$ we let $\Theta(r)$ denote the boundedly overlapping covering of $\R^{n-1}$ formed by balls of radius $r$ centred on lattice points in $c_n r\Z^{n-1}$. We refer to each $\theta \in \Theta(r)$ as a \textit{cap}. For $R \geq 1$, we fix a $C^{\infty}$ partition of unity $(\psi_{\theta})_{\theta \in \Theta(R^{-1/2})}$ subordinate to the covering $\Theta(R^{-1/2})$ of $\B^{n-1}$. We also fix $(\tilde{\psi}_{\theta})_{\theta \in \Theta(R^{-1/2})}$, a family of $C^{\infty}(\R^{n-1})$ functions satisfying $\tilde{\psi}_{\theta}(y) = 1$ if $y \in \supp \psi_{\theta}$ and $\supp \tilde{\psi}_{\theta} \subseteq 2 \cdot \theta$, where $2 \cdot \theta$ denotes the cap concentric to $\theta$ but with twice the radius. We may choose these functions so that 
\begin{equation*}
   |\partial_y^{\bdbeta} \psi_{\theta}(y)| + |\partial_y^{\bdbeta} \tilde{\psi}_{\theta}(y)| 
   \lesssim_{\bdbeta}
   R^{|\bdbeta|/2} \qquad \textrm{for all $\bdbeta \in \N_0^{n-1}$.}
\end{equation*}
\subsubsection*{Tubes} For technical reasons, we shall use a slightly different definition for our tubes than that in \S\ref{subsec: geometric max}. Let $R \geq 1$ and $\delta \geq 0$. For $(\theta, \omega) \in \Theta(R^{-1/2}) \times R^{1/2}\Z^{n-1}$, we let
\begin{equation}\label{eq: tube I}
    T_{\theta,\omega} \coloneq \big\{ \bdx \in \R^{n-1} \times (-2R, 2R) : |\partial_y \phi_B^{\lambda}(\bdx;y_{\theta}) - \omega| < R^{1/2 + \delta} \big\}, 
\end{equation}
where $y_{\theta} \in R^{-1/2}\Z^{n-1}$ denotes the centre of the cap $\theta$.

Here the width $R^{1/2 + \delta}$ of the tube is slightly larger than $R^{1/2}$, but other than that the definition \eqref{eq: tube I} agrees with \eqref{eq: original tube}. We let $\bdT^R$ denote the set of all $R$-tubes as defined in \eqref{eq: tube I}. Note that $\bdT^R$ depends on a choice of $\delta  \geq 0$, $\lambda \geq 1$ and ball $B$, but we suppress this in the notation. 

Given $T = T_{\theta, \omega} \in \bdT^R$, we let $2 \cdot T$ denote the tube with the same direction $\theta$ and centre $\omega$ as $T$, but with width $2R^{1/2 + \delta}$ rather than $R^{1/2 + \delta}$. 

Suppose \eqref{eq: multiplicity bound} holds for some set of $R$-tubes as defined in \eqref{eq: original tube} (that is, \eqref{eq: tube I} with $\delta = 0$). Then \eqref{eq: multiplicity bound} still holds if we replace each $\chi_T$ with $\chi_{2 \cdot T}$. Similarly, \eqref{eq: multiplicity bound} also holds with the corresponding sets of $R$-tubes as defined in \eqref{eq: tube I}, provided we also replace the $\varepsilon$ exponent with $\varepsilon + O(\delta)$.

\subsubsection*{Weight function} Given $R \geq 1$ and $\bdx_0 \in \R^n$, we define the weight $w_{B(\bdx_0, R)} \colon \R^n \to \R$ adapted to the ball $B(\bdx_0, R)$ by
\begin{equation}\label{eq: weight}
    w_{B(\bdx_0, R)}(\bdx) \coloneq \big(1 + R^{-1}|\bdx - \bdx_0|\big)^{-100n}.
\end{equation}

\subsubsection*{Asymptotic notation} Given a list of objects $L$ and real numbers $A$, $B \geq 0$, here and throughout we write $A \lesssim_L B$ or $B \gtrsim_L A$ to indicate $A \leq C_L B$ for some constant $C_L$ which depends only on items in the list $L$ and possibly a choice of dimension $n$, H\"ormander-type phase $\phi$ and amplitude $a$. We write $A \sim_L B$ to indicate $A \lesssim_L B$ and $B \lesssim_L A$. Throughout the following, we shall work with some choice of large parameter $R \geq 1$. We also write $A \lessapprox_L B$ or $B \gtrapprox_L A$ if for all $\eta > 0$ we have $A \lesssim_{L, \eta} R^{\eta} B$.




\section{Analytic preliminaries}\label{sec: anal prelim}




\subsection{Wave packet decomposition}\label{sec: wavepackets}

Let $\phi \colon \D^n \to \R$ be a H\"ormander-type phase, $1 \leq R \leq \lambda$ and fix a ball $B:=B(\bdx_B, R) \subseteq \R^n$. Let $\tilde{f}\coloneq e^{i \phi^{\lambda}(\bdx_B;\,\cdot\,)}f$. Observe that
\begin{equation*}
    H^{\lambda}[\phi, a]f(\bdx_B + \bdx) = \int_{\R^{n-1}} e^{i\phi_B^{\lambda}(\bdx; y)} a_B^{\lambda}(\bdx; y) e^{i \phi^{\lambda}(\bdx_B; y)} f(y)\,\ud y = H^{\lambda}[\phi_B, a_B]\tilde{f}(\bdx),
\end{equation*}
where we define
$$\phi_B(\bdx;y):=\phi(\bdx_B/\lambda+\bdx;y)-\phi(\bdx_B/\lambda;y),$$ and we define 
$$a_B(\bdx;y):=a(\bdx_B/\lambda+\bdx;y).$$
Since $|\tilde{f}|=|f|$, it follows that \eqref{eq: H pq} holds if and only if 
$$\|H^\lambda[\phi_B;a_B]\tilde{f}\|_{L^q(B(0,R))}\lesssim_{\phi,a,\eps}R^{\alpha+\eps}\|\tilde{f}\|_{L^p(\B^{n-1})}.$$
With the ball $B(\bdx_B,R)$ fixed, for the remainder of the section we let $H^\lambda$ denote the recentred operator $H^\lambda[\phi_B;a_B]$, and $B_R$ denote $B(0,R)$, and write $f$ for $\tilde{f}$.

Given $f \in C_c^{\infty}(\B^{n-1})$, using the definitions of the previous section, we write $f_{\theta} \coloneq f \cdot \psi_{\theta}$ for all caps $\theta \in \Theta(R^{-1/2})$. By performing a Fourier series decomposition, we may express $f_{\theta} = f_{\theta} \cdot \tilde{\psi}_{\theta}$ as
\begin{equation*}
    f_{\theta}(y) = \sum_{\omega \in R^{1/2}\Z^{n-1}} f_{\theta, \omega}(y) \quad \textrm{where} \quad f_{\theta,\omega}(y) \coloneq R^{(n-1)/2} (f_{\theta})\;\widecheck{}\;(\omega)e^{-i \omega \cdot y}\tilde{\psi}_{\theta}(y)
\end{equation*}
for Fourier coefficients $(f_{\theta})\;\widecheck{}\;(\omega)$ as defined in \eqref{eq: FT}. Given a choice of $\delta > 0$, we let $\bdT^R$ be as in \S\ref{sec: notation} and define 
\begin{equation*}
   f_T \coloneq f_{\theta,\omega}  \qquad \textrm{for} \quad T = T_{\theta, \omega} \in \bdT^R.
\end{equation*}
For $\theta \in \Theta(R^{-1/2})$, let $\bdT^R[\theta] \coloneq \{T_{\theta, \omega} \in \bdT^R : \omega \in R^{1/2}\Z^{n-1}\}$, so that
\begin{equation*}
    f=\sum_{\theta\in\Theta(R^{-1/2})}\sum_{\omega\in R^{1/2}\mathbb{Z}^{n-1}} f_{\theta,\omega}=\sum_{T\in\bdT^R}f_T. 
\end{equation*}



\subsection{Properties of the wavepackets} We now discuss various (standard) properties of the functions $f_T$ introduced above. 

\subsubsection*{Orthogonality} For $T_{\theta, \omega} \in \bdT^R$ we have
\begin{equation}\label{eq: wp orthog 1}
    |(f_{\theta})\;\widecheck{}\;(\omega)|^2 \sim R^{-(n-1)/2}\|f_{T_{\theta, \omega}}\|_{L^2(\R^{n-1})}^2.
\end{equation}
Given $\T \subseteq \bdT^R$ and $\theta \in \Theta(R^{-1/2})$, if we define
\begin{equation*}
   \T[\theta] \coloneq \T \cap \bdT^R[\theta], 
\end{equation*}
 then it follows from \eqref{eq: wp orthog 1} and Plancherel's theorem for Fourier series that
\begin{equation}\label{eq: wp orthog 2}
   \Big\|\sum_{T \in \T[\theta]} f_T\Big\|_2^2 \lesssim \sum_{T \in \T[\theta]} \|f_T\|_2^2.
\end{equation}
Since the $(\tilde{\psi}_{\theta})_{\theta \in \Theta(R^{-1/2})}$ have boundedly overlapping support, it further follows that
\begin{equation}\label{eq: wp orthog 3}
     \Big\|\sum_{T \in \T} f_T\Big\|_{L^2(\R^{n-1})}^2 \lesssim  \sum_{T \in \T} \|f_T\|_{L^2(\R^{n-1})}^2 \lesssim \|f\|_{L^2(\R^{n-1})}^2. 
\end{equation}

\subsubsection*{Spatial concentration} We now investigate the action of our operator $H^{\lambda}$ on the wavepackets.

\begin{lemma}  For any $T \in \bdT^R$, we have
\begin{equation*}
    \sup_{\bdx \in B(0,R)\setminus T} |H^{\lambda} f_T(\bdx)|
    \lesssim_{\delta, N}
   ( 1+ R^{-1}\dist(T, B(0,R)))^{-N} R^{-N}\|f\|_{L^2(\R^{n-1})}
\end{equation*}
for all $N \in \N$.
\end{lemma}

The lemma follows from a straightforward non-stationary phase (integration-by-parts) argument that uses the fact that $\phi_B^\lambda(0;y)=0$ for all $y$. We omit the proof.

\subsubsection*{$L^p$ estimates} Finally, we describe $L^p$ estimates for individual wavepackets. 

\begin{lemma}\label{lem: Bernstein} Let $2\leq p\leq\infty$. For all $T\in\bdT^R$, we have
    \begin{equation*}
 \|H^\lambda f_T\|_{L^p(w_{B_R})}\lesssim R^{-(n+1)(1/4-1/(2p))}R^{1/2}\|f_T\|_{L^2(\R^{n-1})}.  
\end{equation*}
\end{lemma}
Here $w_{B_R}$ is the weight function defined in \eqref{eq: weight}. The case $p = 2$ of Lemma~\ref{lem: Bernstein} follows directly from H\"ormander's $L^2$-bound \cite{hormander}. For $p > 2$, morally the estimate follows from Bernstein's inequality using the microlocalisation of $f_T$ to lower the $L^p$ to an $L^2$ norm and then applying the $L^2$ result. 

\begin{proof}[Proof (of Lemma~\ref{lem: Bernstein})]  We have already noted the $p = 2$ case follows from H\"ormander's $L^2$-bound \cite{hormander}. The $L^{\infty}$ estimate follows from the Cauchy--Schwarz inequality:
\begin{equation*}
  \|H^{\lambda}f_T\|_{L^{\infty}(B_R)} \lesssim \|f_T\|_{L^1(\R^{n-1})} \lesssim R^{-(n-1)/4}\|f_T\|_{L^2(\R^{n-1})}.  
\end{equation*}
The result now follows by interpolation (via logarithmic convexity of $L^p$).
\end{proof}




\subsection{Refined decoupling in the mixed signature case} As mentioned in the introduction, the proof of Theorem \ref{thm: main} is based on a refined decoupling estimate for H\"ormander-type phases of general signatures. To state the relevant results, let
\begin{equation}\label{eq: dec exp}
    r_{n,\sigma}
\coloneqq 
\frac{2(n+3+\sigma)}{n-1+\sigma}, \qquad
d(n,\sigma)
\coloneqq 
\frac{n-1-\sigma}{2(n+3+\sigma)},
\end{equation}
so that $r_{n, \sigma}$ is as in \eqref{eq: dec exp 0}. Given a set of $R$-tubes $\T \subseteq \bdT^R$, we define the \textit{multiplicity function}
\begin{equation}\label{eq: mult fn}
   M_{\T} \colon \bdQ^{R^{1/2}} \to \N_0, \qquad M_{\T}(Q) \coloneq \#\big\{T \in \T : T \cap Q \neq \varnothing \big\} \quad \textrm{for all $Q \in \bdQ^{R^{1/2}}$.}
\end{equation}
Recall also that for a set $\cQ \subseteq \bdQ^{R^{1/2}}$ we let $\bdE(\cQ)$ denote their union.

Refined decoupling inequalities are obtained by iterative applications of a corresponding `classical' decoupling inequality. In the case of \cite[Theorem 4.2]{GIOW2018}, the classical inequality is the Bourgain--Demeter theorem \cite{bourgain--demeter}. In our setting, the classical inequality is the following.

\begin{theorem}[Schippa \cite{Schippa}]\label{thm: dec schippa} Let $1\leq R\leq \lambda$ and $\phi:\D^n\to \R$ be a H\"ormander-type phase of signature $\sigma$, and $a\in C^\infty_c(\D^n)$. Let $B(\bdz,R)\subset\R^n$ be any ball of radius $R$. Then
$$\|H^\lambda[\phi;a] f\|_{L^{r_{n,\sigma}}(B(\bdz,R))}\lessapprox R^{d(n,\sigma)}\big(\sum_{\theta\in\Theta(R^{-1/2})}\|H^\lambda[\phi;a] f_\theta\|_{L^{r_{n,\sigma}}(w_{B(\bdz,R)})}^2\big)^{1/2},$$ 
where $w_{B(\bdz,R)}$ is a rapidly decaying weight as in \eqref{eq: weight}.
\end{theorem}
By applying the above theorem at many scales, we obtain a refinement which accounts for the multiplicity of intersections between the tubes associated to the wavepackets. We state the result for the recentred operator $H^\lambda$ and the origin-centred ball $B_R$. Recall the wavepacket decomposition of the function $f$ from \S\ref{sec: wavepackets}.

\begin{theorem}\label{thm: refined dec} Let $\Q$ be a set of lattice $R^{1/2}$-cubes contained in a fixed dilate of $B_R$ and let $\T \subseteq \bdT^R$ be such that 
   \begin{enumerate}[i)]
       \item $\|H^{\lambda} f_{T_1}\|_{L^{r_{n, \sigma}}(w_{B_R})} \sim \|H^{\lambda} f_{T_2}\|_{L^{r_{n, \sigma}}(w_{B_R})}$ \text{for all $T_1$, $T_2\in\T$;}\\[2pt]
       \item $M_{\T}(Q) \leq M$ for all $Q \in \cQ$.
   \end{enumerate}
   If we set $E \coloneq \bdE(\cQ)$, then it follows that 
   \begin{equation*}
       \Big\|\sum_{T \in \T} H^{\lambda} f_T \Big\|_{L^{r_{n, \sigma}}(E)}
   \lessapprox
   \bigg(\frac{M}{\#\T}\bigg)^{1/2-1/{r_{n,\sigma}}} R^{d(n,\sigma)}\Big(\sum_{T\in\T }\|H^{\lambda} f_T\|_{L^{r_{n, \sigma}}(w_{B_R})}^2\Big)^{1/2}.
   \end{equation*}
\end{theorem}

Here $w_{B_R}$ is the rapidly decaying weight function, as defined in \eqref{eq: weight}. We omit the proof of Theorem~\ref{thm: refined dec} as it follows very similar arguments to those in the literature: in particular, we refer the reader to the proofs of \cite[Theorem 4.2]{GIOW2018} and \cite[Theorem 5.1]{ILX}.




\section{Geometric preliminaries}\label{sec: geo prelim}




\subsection{Shadings} We begin by recalling some basic definitions used commonly in the Kakeya literature.

\begin{definition}[Shading] Let $R \geq 1$ and $\T \subseteq \bdT^R$. A \textit{shading} is a mapping $Y\colon \T\to \cP(\R^n)$ such that for all $T \in \T$ we have 
\begin{equation*}
    Y(T) = \bigcup_{Q \in \cQ_Y[T]} Q \quad \textrm{for some $\cQ_Y[T] \subseteq \bdQ^{R^{1/2}}$ with $T \cap Q \neq \varnothing$ for all $Q \in \cQ_Y[T]$. }
\end{equation*}
In this case, the pair $(\T,Y)$ is called a \textit{shaded family of $R$-tubes} and we write
\begin{equation*}
    \cQ_Y \coloneq \bigcup_{T \in \T} \cQ_Y[T].
\end{equation*}
For each $Q\in\cQ_Y$, define the \textit{restricted multiplicity function} 
\begin{equation}
    \label{eq: res mult fn}
    M_Y(Q)\coloneq\#\big\{T\in\T: Q\in \cQ_Y[T]\big\}.
\end{equation}   
\end{definition}
\begin{definition}[Density]\label{dfn: dense}  Let $R \geq 1$, $\rho >0$ and $(\T,Y)$ be a shaded family of $R$-tubes. We say $(\T,Y)$ is \textit{$\rho$-dense} if $\rho|T|\leq|Y(T)|<2\rho|T|$ for all $T\in\T$.
\end{definition}

\begin{definition}[Two-ends \cite{Wolff1995}]\label{dfn: two ends}
     Let $R \geq 1$, $0 < \varepsilon_2 \leq \varepsilon_1 \leq 1$ and $(\T,Y)$ be a shaded family of $R$-tubes. We say that $(\T,Y)$ is \textit{$(\varepsilon_1,\varepsilon_2)$-two-ends} if there exists an absolute constant $C$ such that
\begin{equation*}
     \sup_{x \in \R^n} |Y(T)\cap B(x,R^{1-\varepsilon_1})|
    \leq C R^{-\varepsilon_2}|Y(T)|\qquad\text{for all $T\in\T$.} 
\end{equation*}
\end{definition}




\subsection{Universal two-ends multiplicity bound} Given any H\"ormander-type phase $\phi$ and any associated set of $R$-tubes $\T \subseteq \bdT^R$, we have the following universal estimate. 

\begin{lemma}[Universal two-ends multiplicity bound]\label{lem: univ two-ends} Let $\T \subseteq \bdT^R$ and $(\T, Y)$ be an $(\varepsilon_1, \varepsilon_2)$-two-ends shaded family of $R$-tubes. Then
\begin{equation*}
    \#\{Q \in \cQ_Y : M_Y(Q) > \mu \} \lesssim_{\varepsilon} R^{(\varepsilon_1+\delta)(n-1)} \mu^{-1} (\#\T)^{1/2} \#\cQ_Y \qquad \textrm{for all $\mu > 0$.}  
\end{equation*}
\end{lemma}

Here $\delta > 0$ is the parameter appearing in the definition of the $R$-tubes in $\bdT^R$. Lemma~\ref{lem: univ two-ends} corresponds to Bourgain's $\frac{n+1}{2}$-dimension bound for Kakeya sets \cite{Bourgain1991a} and can be proved directly using a dual version of Bourgain's bush argument \cite{Bourgain1991a}. 

\begin{proof}[Proof (of Lemma~\ref{lem: univ two-ends})] By Markov's inequality,
\begin{equation}\label{eq: univ two-ends 1}
  \#\{Q \in \cQ_Y : M_Y(Q) > \mu \} \leq \mu^{-1}\sum_{Q \in \cQ_Y} M_Y(Q).    
\end{equation}

To estimate the right-hand side of \eqref{eq: univ two-ends 1}, we double-count the set of triples
\begin{equation*}
    \fkT \coloneq \big\{(T, Q_1, Q_2) \in \T \times \cQ_Y^2 : Q_j \in \cQ_Y[T] \textrm{ for $j = 1$, $2$, } d(T \cap Q_1, T \cap Q_2) \geq R^{1- \varepsilon_1}/2\big\},
\end{equation*}
where here $d(F_1, F_2) \coloneq \inf\{|x_1-x_2| : x_1 \in F_1, x_2 \in F_2\}$ for nonempty $F_1$, $F_2 \subseteq \R^n$. On the one hand,
\begin{equation}\label{eq: univ two-ends 2}
    \#\fkT \leq \sum_{\substack{Q_1, Q_2 \in \cQ_Y \\ d(Q_1, Q_2) \geq R^{1- \varepsilon_1}/4}} \#\big\{T \in \T : Q_j \in \cQ_Y[T] \textrm{ for $j = 1$, $2$} \big\}.
\end{equation}
If $T_j = T_{\theta_j, \omega_j} \in \bdT^R$, $j = 1$, $2$, satisfy $\mathrm{diam}(2\cdot T_1 \cap 2\cdot T_2) \gtrsim R^{1 - \varepsilon_1}$, then 
\begin{equation*}
  d(\theta_1, \theta_2)\lesssim R^{-1/2+\varepsilon_1+\delta},\quad\text{and}\quad|\omega_1 - \omega_2| \lesssim R^{1/2+\varepsilon_1 + \delta}.   
\end{equation*} 
This follows for a class of quadratic curves from \cite[p.49, Lemma 20]{WisewellThesis}; for general curves, the authors are unaware of an explicit reference and so a proof is presented in Appendix~\ref{sec: tube intersect}.\footnote{Strictly speaking, at this point we should assume the setup is suitably localised (see the hypotheses in Lemma~\ref{lem: tube intersect}). This is possible by breaking the domain and the tubes into a bounded number of pieces.} In light of the above, by the separation of the centres and directions,
\begin{equation*}
  \max_{\substack{Q_1, Q_2 \in \cQ_Y \\ d(Q_1, Q_2) \geq R^{1- \varepsilon_1}/4}}  \#\big\{T \in \T : Q_j \in \cQ_Y[T] \textrm{ for $j = 1$, $2$} \big\} \lesssim R^{2(n-1)(\delta+\varepsilon_1)}.
\end{equation*}
 Thus, recalling \eqref{eq: univ two-ends 2}, we have
\begin{equation}\label{eq: univ two-ends 3}
    \#\fkT \lesssim R^{2(n-1)(\delta+\varepsilon_1)} [\#\cQ_Y]^2.
\end{equation}

On the other hand, by the two-ends hypothesis,
\begin{align*}
    \#\fkT &= \sum_{T \in \T} \#\big\{(Q_1, Q_2) \in \cQ_Y^2 : Q_j \in \cQ_Y[T], \textrm{ $j = 1$, $2$, } d(T \cap Q_1, T \cap Q_2) \geq R^{1- \varepsilon_1}/2\big\} \\
    &\gtrsim \sum_{T \in \T} \big[\#\big\{Q \in \cQ_Y : Q \in \cQ_Y[T] \big\}\big]^2.
\end{align*}
Thus, by the Cauchy--Schwarz inequality and a change of order of summation, 
\begin{equation}\label{eq: univ two-ends 4}
    \#\fkT \gtrsim [\#\T]^{-1} \Big(\sum_{Q \in \cQ_Y} M_Y(Q)\Big)^2
\end{equation}
Combining \eqref{eq: univ two-ends 3} and \eqref{eq: univ two-ends 4}, we have
\begin{equation}\label{eq: univ two-ends 5}
    \sum_{Q \in \cQ_Y} M_Y(Q) \lesssim R^{(\varepsilon_1+\delta)(n-1)} [\#\T]^{1/2} \#\cQ_Y.
\end{equation}
Applying \eqref{eq: univ two-ends 5} to the right-hand side of \eqref{eq: univ two-ends 1} gives the desired bound. \end{proof}




\subsection{Conditional multiplicity bound} Any Kakeya estimate as in \eqref{eq: multiplicity bound} can be converted into a corresponding bound on the level sets of the restricted multiplicity function $M_Y$.

\begin{lemma}[Conditional multiplicity bound]\label{lem: cond mult} Let $R \geq 1$, $0 < \rho \leq 1$, $\beta \geq 0$ and $1 < s < \infty$. Suppose $(\T, Y)$ is a $\rho$-dense shaded family of $R$-tubes which satisfies
    \begin{equation}\label{eq: improv mult 1}
      \big\|\sum_{T \in \T} \chi_{2 \cdot T}\big\|_{L^s(\R^n)} \leq A R^{(n+1)/(2s)+\beta/2} \big[\#\T \big]^{1/s}
    \end{equation}
for some $A \geq 1$. Then
\begin{equation*}
    \#\{Q \in \cQ_Y : M_Y(Q) > \mu \} \lesssim A^{s'} \mu^{-1} \rho^{-(s'-1)} R^{\beta s'/2} \#\cQ_Y \qquad \textrm{for all $\mu > 0$.} 
\end{equation*}
\end{lemma}

\begin{proof}[Proof (of Lemma~\ref{lem: cond mult})] By Markov's inequality and a change of order of summation,
\begin{align}
\nonumber
  \#\{Q \in \cQ_Y : M_Y(Q) > \mu \} &\leq \mu^{-1}\sum_{T \in \T} \#\{Q \in \cQ_Y : Q \in \cQ_Y[T]\} \\
  \label{eq: improv mult 3}
  &= \mu^{-1} R^{-n/2} \sum_{T \in \T} |Y(T)|.
\end{align}
To estimate the right-hand side of \eqref{eq: improv mult 3}, we apply the hypothesis \eqref{eq: improv mult 1}. Let $E \coloneq \bdE(\cQ_Y)$ and note that
\begin{equation*}
  \sum_{T \in \T} |Y(T)| \leq  \int_E \sum_{T \in \T} \chi_{2 \cdot T} \leq |E|^{1/s'} \Big\|\sum_{T \in \T} \chi_{2 \cdot T} \Big\|_{L^s(\R^n)}. 
\end{equation*}
Assuming \eqref{eq: improv mult 1} holds, we have
\begin{equation}\label{eq: improv mult 4}
    \sum_{T \in \T} |Y(T)| \leq A R^{\beta/2} \big[R^{(n+1)/2}\#\T \big]^{1/s} |E|^{1/s'}.
\end{equation}
Since $(\T, Y)$ is $\rho$-dense, it follows that 
\begin{equation*}
    R^{(n+1)/2} \#\T \lesssim \rho^{-1}\sum_{T \in \T} |Y(T)|.
\end{equation*}
 Applying this to \eqref{eq: improv mult 4} and rearranging the resulting expression yields
\begin{equation}\label{eq: improv mult 5}
   \sum_{T \in \T} |Y(T)| \lesssim A^{s'} \rho^{-(s'-1)} R^{\beta s'/2} |E| \lesssim A^{s'}\rho^{-(s'-1)} R^{n/2 + \beta s'/2} \#\cQ_Y.
\end{equation}
Applying \eqref{eq: improv mult 5} to the right-hand side of \eqref{eq: improv mult 3} gives the desired bound. 
\end{proof}




\section{Proof of Theorem~\ref{thm: main}: Setting up the argument}




\subsection{Induction-on-scale} Given $E \subseteq \R^{n-1}$ a measurable set with $0 < |E| < \infty$, define
\begin{equation*}
    \|f\|_{L_{\avg}^2(E)} \coloneqq \Big(|E|^{-1}\int_E |f|^2 \Big)^{1/2}.
\end{equation*}
Rather than work directly with $L^q \to L^q$ inequalities, it is useful to replace the right-hand norm in \eqref{eq: H pq} with the following functional.

\begin{definition}
Let $2 \leq p < \infty$ and $R \geq 1$. For $f \in L^2(\B^{n-1})$, let 
    \begin{equation*}
        \tnorm{f}_{L_R^p(\B^{n-1})} \coloneq \|f\|_{L^2(\B^{n-1})}^{2/p}\sup_{\theta \in \Theta(R^{-1/2})}\|f\|_{L_{\avg}^2(\theta)}^{1-2/p}.
    \end{equation*}
\end{definition}

If $E \subseteq \B^{n-1}$ is a measurable set, then 
\begin{equation}\label{eq: induction 1}
   \tnorm{\chi_{E}}_{L_R^p(\B^{n-1})} \leq \|\chi_{E}\|_{L^p(\B^{n-1})},
\end{equation}
which partially motivates the choice of notation.

Henceforth we let $n \geq 3$ be odd, fix the H\"ormander-type phase $\phi \colon \D^n \to \R$ with signature $\sigma$ and an amplitude $a\in C^\infty_c(\mathbb{D}^n)$. We fix some $\varepsilon > 0$ and let $\bdR_{\varepsilon} \geq 1$ and $\bdC_{\varepsilon} \geq 1$ be constants chosen sufficiently large to satisfy the forthcoming requirements of the proof. Without loss of generality we may assume $0 < \varepsilon < \varepsilon_{\circ}$, where $0 < \varepsilon_{\circ} < 1/100n$ is a small dimensional constant. Our goal is to establish the following conditional bounds. 

\begin{proposition}\label{prop: induction} Let $1 \leq s \leq \frac{n+1}{n-1}$, $\beta \geq 0$ and define $q \coloneq q_{n, \sigma}(s, \beta)$ and $p \coloneq p_{n, \sigma}(s, \beta)$ as in \eqref{eq: target exp arrival} and \eqref{eq: target exp domain}. Let $1 \leq R \leq \lambda$ and $R$-ball $B\subset \R^n$. If $\bdK_{s \to s}^*(\phi; \beta)$ holds, then 
    \begin{equation}\label{eq: induction I}
        \|H^{\lambda}[\phi;a] f\|_{L^q(B)}\leq \bdC_{\varepsilon} R^{\varepsilon} \tnorm{f}_{L_R^p(\B^{n-1})} \qquad \textrm{for all $f\in C_c^{\infty}(\B^{n-1})$.}
    \end{equation}
\end{proposition}

Once \eqref{eq: induction I} is established, we can apply it to characteristic functions of measurable sets, using \eqref{eq: induction 1}, to deduce that $H^{\lambda}[\phi;a]$ satisfies a restricted weak-type $L^q$ inequality. Using elementary real interpolation techniques, this leads to Theorem~\ref{thm: main}.  

The proof of Proposition~\ref{prop: induction} follows an induction-on-scale argument, similar in structure to that of Wolff~\cite{Wolff2001}. \medskip

\noindent\textbf{Base case.} Suppose $1 \leq R \leq \bdR_{\varepsilon}$. Provided $\bdC_{\varepsilon}$ is chosen sufficiently large, depending on $\bdR_{\varepsilon}$, the base case follows from the trivial estimate $\bdH_{1\to\infty}(\phi)$.\medskip

We now fix some $\bdR_{\varepsilon} < R \leq \lambda$ and assume the following induction hypothesis. 

\noindent\textbf{Induction hypothesis.} For all $1 \leq R' \leq R/2$ and all $R'$-ball $B'\subset \R^n$, we have
\begin{equation}\label{eq: induction 3}
        \|H^{\lambda}[\phi;a] f\|_{L^q(B')}\leq \bdC_{\varepsilon} (R')^{\varepsilon} \tnorm{f}_{L_{R'}^p(\mathbb{B}^{n-1})}\qquad \textrm{for all $f\in C_c^{\infty}(\B^{n-1})$.}
    \end{equation}

The goal is now to show \eqref{eq: induction I} holds for our fixed choice of $R$ and an arbitrary $R$-ball $B(\bdx_B,R)$. As in \S\ref{sec: wavepackets}, we let $H^\lambda$ denote the recentred operator $H^\lambda[\phi_B;a_B]$ and let $B_R=B(0,R)$. For the induction step, we will establish the following 
$$ \|H^{\lambda}f\|_{L^q(B_R)}\leq \bdC_{\varepsilon} R^{\varepsilon} \tnorm{f}_{L_R^p(\B^{n-1})} \qquad \textrm{for all $f\in C_c^{\infty}(\B^{n-1})$.}$$

For technical reasons, we introduce a further constant $1 \leq \bar{C} \leq 2\bdC_{\varepsilon}$, again chosen sufficiently large to satisfy the forthcoming requirements of the proof. Without loss of generality, we may assume the following:
\begin{flalign}\label{eq: non deg}
&\textrm{\textbf{Non-degeneracy hypothesis:} $\|H^{\lambda} f\|_{L^q(B_R)}\geq 2\bar{C} \|f\|_{L^2(\mathbb{B}^{n-1})}$.} &
\end{flalign}
Indeed, if \eqref{eq: non deg} were to fail, then it is easy to see that, provided $\bdC_{\varepsilon}$ is chosen large relative to $\bar{C}$, we automatically have a favourable estimate.  




\subsection{Pigeonholing}\label{subsec: pigeon} As groundwork for the main argument, we perform a wavepacket decomposition of the operator and then dyadic pigeonhole in various statistics associated to this decomposition. 

Let $\bdT^R$ be as in \S\ref{sec: notation}, where we take $\delta \coloneq \varepsilon^{10} < \varepsilon/100n$. Recall the wavepacket decompositions defined in \S\ref{sec: wavepackets}. Define $\T_1 \coloneqq 
\{T\in\bdT^R: T \cap B(0,R) \neq \varnothing\}$.
Let $f_T$ denote the wavepackets of $f \in C_c^{\infty}(\B^{n-1})$, for all $T \in \T_1$.\smallskip

By the spatial concentration property of the wavepackets, provided $\bar{C} \geq 1$ is chosen sufficiently large,
\begin{equation}\label{eq: pigeon -1}
   \|H^{\lambda}f\|_{L^q(B_R)} \leq  \Big\|
\sum_{T \in \T_1} H^{\lambda}f_T \Big\|_{L^q(B_R)} + \bar{C}\|f\|_{L^2(\mathbb{B}^{n-1})}.
\end{equation}
By the non-degeneracy hypothesis \eqref{eq: non deg}, the second term on the right-hand side of the above display may be absorbed into the left-hand side to give
\begin{equation}\label{eq: pigeon 0}
   \|H^{\lambda}f\|_{L^q(B_R)} \lesssim  \Big\|
\sum_{T \in \T_1} H^{\lambda}f_T \Big\|_{L^q(B_R)}. 
\end{equation}

\subsubsection*{1. Individual wavepacket mass} We begin by pigeonholing the wavepackets to ensure the first hypothesis of the refined decoupling inequality (Theorem~\ref{thm: refined dec}) holds. 

Let $r_{n, \sigma}$ be the exponent defined in \eqref{eq: dec exp}. For any $T \in \T_1$, we trivially have
\begin{equation*}
    \|H^{\lambda}f_T\|_{L^{r_{n,\sigma}}(w_{B_R})} \lesssim R^{n/r_{n, \sigma}} \|H^{\lambda}f_T\|_{L^{\infty}(\R^n)} \leq R^{10n}\|f\|_{L^2(\mathbb{B}^{n-1})}.
\end{equation*}
On the other hand, the $T \in \T_1$ satisfying
\begin{equation*}
    \|H^{\lambda}f_T\|_{L^{r_{n,\sigma}}(w_{B_R})} \lesssim R^{-10n}\|f\|_{L^2(\mathbb{B}^{n-1})}
\end{equation*}
contribute minimally to the right-hand side of \eqref{eq: pigeon 0} and can be dealt with using the non-degeneracy hypothesis as in \eqref{eq: pigeon -1} and \eqref{eq: pigeon 0}. We may therefore apply dyadic pigeonholing to find some dyadic number $R^{-10n} \lesssim \beta \lesssim R^{10n}$ and a subset $\T_2 \subseteq \T_1$ such that the following hold:
\begin{enumerate}
    \item[P1 i)] $\|H^{\lambda} f_T\|_{L^{r_{n, \sigma}}(w_{B_R})} \sim \beta\|f\|_{L^2(\mathbb{B}^{n-1})}$ for all $T \in \T_2$;
    \item[P1 ii)]
    \begin{equation}\label{eq: pigeon 1}
   \|H^{\lambda}f\|_{L^q(B_R)}
\lessapprox
\Big\|
\sum_{T\in\T_2}H^{\lambda}f_T
\Big\|_{L^q(B_R)}.
\end{equation}
\end{enumerate}

\subsubsection*{2. Number of parallel wavepackets} Recall the definition of $\bdT^R[\theta]$ from \S\ref{sec: anal prelim}.  By a further dyadic pigeonholing, we can find some dyadic $1 \leq m \lesssim R^{(n-1)/2}$ and a subset $\T_{3} \subseteq \T_2$ with the following properties:
\begin{enumerate}
    \item[P2 i)]  For all $\theta \in \Theta(R^{-1/2})$, the set $\T_3[\theta] \coloneq \T_3 \cap \bdT^R[\theta]$ satisfies
    \begin{equation*}
        \textrm{either} \quad \# \T_3[\theta] = 0 \quad \textrm{or} \quad \# \T_3[\theta] \sim m.
    \end{equation*}
    
    \item[P2 ii)]
    \begin{equation}\label{eq: pigeon 2}
   \Big\|
\sum_{T\in\T_2}H^{\lambda}f_T
\Big\|_{L^q(B_R)}
\lessapprox
\Big\|
\sum_{T\in\T_3}H^{\lambda}f_T
\Big\|_{L^q(B_R)}.
\end{equation}
\end{enumerate}

\subsubsection*{3. Total mass over a cube} We introduce a new parameter and scale
\begin{equation}\label{eq: epsilon 1}
    \varepsilon_1 \coloneq \varepsilon^2 < \varepsilon_{\circ}\varepsilon \qquad \textrm{and} \qquad R' \coloneq R^{1- \varepsilon_1}.
\end{equation}
By pigeonholing, we can find a set $\mathcal{D}$ of lattice $R'$-cubes such that $D \cap B_R \neq \varnothing$ for all $D \in \mathcal{D}$, if $D_1$, $D_2 \in \mathcal{D}$ satisfy $d(D_1, D_2) < 10R'$, then $D_1 = D_2$ and 
\begin{equation*}
    \Big\|
\sum_{T\in\T_3}H^{\lambda}f_T
\Big\|_{L^q(B_R)}^q \lesssim \sum_{D \in \mathcal{D}} \Big\|
\sum_{T\in\T_3}H^{\lambda}f_T
\Big\|_{L^q(D)}^q.
\end{equation*}
By dyadic pigeonholing, there exists a set of lattice $R^{1/2}$-cubes $\cQ_1$ such that for each $Q \in \cQ_1$ there is a unique $D \in \mathcal{D}$ satisfying $Q \cap D \neq \varnothing$ and the following hold:
\begin{enumerate}
    \item[P3 i)] For all $Q_1, Q_2 \in\cQ_1$, we have
    \begin{equation*}
        \Big\|\sum_{T \in \T_3} H^{\lambda} f_T \Big\|_{L^q(Q_1)} \sim  \Big\|\sum_{T \in \T_3} H^{\lambda} f_T \Big\|_{L^q(Q_2)}.
    \end{equation*}
    \item[P3 ii)] Defining $E_1\coloneqq \bdE(\cQ_1)$,
    \begin{equation}\label{eq: pigeon 3}
       \Big\|\sum_{T \in \T_3} H^{\lambda} f_T\Big\|_{L^q(B_R)} 
    \lessapprox
    \Big\|\sum_{T \in \T_3} H^{\lambda} f_T\Big\|_{L^q(E_1)} .
    \end{equation}
\end{enumerate}

\subsubsection*{4. Segment density} For $T \in \T_3$ and $D \in \mathcal{D}$, let
\begin{align}
    \label{eq: Q_1 def}
    \cQ_1[T] &\coloneq\{Q\in\cQ_1:Q\cap T \neq\varnothing\}, \\
    \nonumber
    \cQ_{1,D}[T]&\coloneq\{Q\in\cQ_1[T]:Q \cap D \neq \varnothing\}.
\end{align} 
Note that the properties of $\mathcal{D}$ and $\cQ_1$ ensure that the family of sets $\{\cQ_{1,D}[T] : D \in \mathcal{D}\}$ forms a partition of $\cQ_1[T]$.  For each $D \in \mathcal{D}$, we define a shading
\begin{equation}\label{eq: pigeon 4}
   \bar Y_D(T)\coloneq\bigcup_{Q\in \cQ_{1,D}[T]}Q \qquad \textrm{for all $T \in \T_3$.} 
\end{equation}
For each $T\in\T_3$ and a dyadic number $\kappa$, let 
\begin{align}
    \label{eq: pigeon 5}
    \mathcal{D}(T)&\coloneq\{D\in\mathcal{D}:D\cap T\neq\varnothing\}, \\ 
    \nonumber
    \mathcal{D}(T;\kappa)&\coloneq\{D\in\mathcal{D}(T):\kappa|T| \leq |\bar Y_D(T)| < 2\kappa|T|\}.
\end{align}
By dyadic pigeonholing, for each $Q\in\cQ_1$, there exists $R^{-1/2 - (n-1)\delta}\lesssim\kappa_Q\lesssim 1$ dyadic such that
\begin{enumerate}
      \item[P4 i)]
      \begin{equation*}
\Big\|\sum_{T\in\T_3}H^\lambda f_T\Big\|_{L^q(Q)}\lessapprox\Big\|\sum_{T\in\T_3}\sum_{D\in\mathcal{D}(T;\kappa_Q)}(H^\lambda f_T)\cdot\chi_D\Big\|_{L^q(Q)}.
\end{equation*}
\end{enumerate}
We now perform a dyadic pigeonholing of the cubes $Q \in \cQ_1$, both the value of $\kappa_Q$ and the size of the right-hand side above. In particular, there exists a subset $\cQ_2 \subseteq \cQ_1$ and a dyadic number $\bar{\kappa} >0$ such that,
\begin{enumerate}
    \item[P4 ii)] For all $Q$, $Q_1$, $Q_2 \in \cQ_2$, we have
    \begin{equation*}
        \kappa_Q = \bar{\kappa} \quad \textrm{and} \quad \Big\|\sum_{T \in \T_3} \sum_{D\in\mathcal{D}(T;\bar\kappa)} \big(H^{\lambda} f_T\big) \cdot \chi_D\Big\|_{L^q(Q_1)} \sim \Big\|\sum_{T \in \T_3} \sum_{D\in\mathcal{D}(T;\bar\kappa)} \big(H^{\lambda} f_T\big) \cdot \chi_D\Big\|_{L^q(Q_2)}.
    \end{equation*}
    \item[P4 iii)] $\# \cQ_2 \gtrapprox \#\cQ_1$. 
\end{enumerate}
Consequently, in view of P3 i) and P4 iii) above,
\begin{equation}\label{eq: pigeon 6}
     \Big\|\sum_{T\in\T_3}H^{\lambda} f_T\Big\|_{L^q(E_1)}^q =\sum_{Q\in\cQ_1} \Big\|\sum_{T\in\T _{3}}H^{\lambda} f_T\Big\|_{L^q(Q)}^q \lessapprox \sum_{Q\in\cQ_2}\Big\|\sum_{T\in\T_3}H^{\lambda} f_T\Big\|_{L^q(Q)}^q.
\end{equation}
Here we have used the non-degeneracy hypothesis \eqref{eq: non deg} to remove the error terms. The right-hand side of \eqref{eq: pigeon 6} can be bounded using P4 i) and P4 ii) to give
\begin{equation*}
    \Big\|\sum_{T\in\T_3}H^{\lambda} f_T\Big\|_{L^q(E_1)}^q \lessapprox\sum_{Q\in\cQ_2}\Big\|\sum_{T\in\T_3}\sum_{D\in\mathcal{D}(T;\bar\kappa)} \big(H^{\lambda} f_T\big) \cdot \chi_D\Big\|_{L^q(Q)}^q.
\end{equation*}
Thus, if we define $E_2 \coloneq \bdE(\cQ_2)$, then
\begin{equation}\label{eq: pigeon 8}
    \Big\|\sum_{T\in\T_3}H^{\lambda} f_T\Big\|_{L^q(E_1)}\lessapprox \Big\|\sum_{T\in\T_3}\sum_{D\in\mathcal{D}(T;\bar\kappa)} \big(H^{\lambda} f_T\big) \cdot \chi_D\Big\|_{L^q(E_2)}.
\end{equation}

\subsubsection*{5. Number of segments}  We partition $\T_3$ into sets $\T_3(\alpha)$ where $\alpha$ varies over dyadic numbers satisfying $1 \lesssim \alpha \lesssim R^{\varepsilon_1}$ and
 \begin{equation}\label{eq: pigeon 9}
      \#\mathcal{D}(T;\bar\kappa) \sim \alpha  \qquad \textrm{for all $T \in \T_3(\alpha)$.}
 \end{equation}
By yet another dyadic pigeonholing, for all $Q \in \cQ_2$ there exists a dyadic number $\alpha_Q \in \N$ such that:
\begin{enumerate}
    \item[P5 i)] 
\begin{equation}\label{eq: pigeon 10}
    \Big\| \sum_{T \in \T_3} \sum_{D\in\mathcal{D}(T;\bar\kappa)}
(H^{\lambda} f_T)\cdot \chi_D
\Big\|_{L^q(Q)}
\lessapprox
\Big\| \sum_{T\in\T_3(\alpha_Q)}
\sum_{D\in\mathcal{D}(T;\bar\kappa)} \big(H^{\lambda} 
f_T\big) \cdot \chi_D \Big\|_{L^q(Q)}.
\end{equation}
\end{enumerate}
We now perform a dyadic pigeonholing of the cubes $Q \in \cQ_2$, both the value of $\alpha_Q$ and the size of the right-hand side \eqref{eq: pigeon 10}. In particular, there exists a subset $\cQ_3 \subseteq \cQ_2$ and a dyadic number $\bar{\alpha} \geq 1$ such that, if we define
\begin{equation}\label{eq: pigeon 11}
    \T \coloneq \T_3(\bar{\alpha}),
\end{equation}
then the following hold:
\begin{enumerate}
    \item[P5 ii)] For all $Q$, $Q_1$, $Q_2 \in \cQ_3$, we have
    \begin{equation*}
        \alpha_Q = \bar{\alpha} \quad \textrm{and} \quad \Big\|\sum_{T \in \T} \sum_{D\in\mathcal{D}(T;\bar\kappa)} \big(H^{\lambda} f_T\big) \cdot \chi_D\Big\|_{L^q(Q_1)} \sim \Big\|\sum_{T \in \T} \sum_{D\in\mathcal{D}(T;\bar\kappa)} \big(H^{\lambda} f_T\big) \cdot \chi_D\Big\|_{L^q(Q_2)}.
    \end{equation*}
    \item[P5 iii)] $\# \cQ_3 \gtrapprox \#\cQ_2$. 
\end{enumerate}
Consequently, in view of the second part of P4 ii) and P5 iii) above,
\begin{equation*}
     \Big\|\sum_{T \in \T_3}\sum_{D\in\mathcal{D}(T;\bar\kappa)}\big(H^{\lambda} f_T\big)\cdot\chi_D\Big\|_{L^q(E_2)}^q \lessapprox \sum_{Q\in\cQ_3} \Big\|\sum_{T \in \T_3}\sum_{D\in\mathcal{D}(T;\bar\kappa)}\big(H^{\lambda} f_T\big)\cdot\chi_D\Big\|_{L^q(Q)}^q.
\end{equation*}
The right-hand side can be bounded using P5 i). Combining it with the first part of P5 ii) and the definition \eqref{eq: pigeon 11}, we obtain
\begin{equation}\label{eq: pigeon 12}
    \Big\|\sum_{T \in \T_3}\sum_{D\in\mathcal{D}(T;\bar\kappa)}\big(H^{\lambda} f_T\big)\cdot\chi_D\Big\|_{L^q(E_2)}^q \lessapprox \sum_{Q\in\cQ_3}\Big\|\sum_{T \in \T}\sum_{D\in\mathcal{D}(T;\bar\kappa)} \big(H^{\lambda} f_T\big) \cdot \chi_D\Big\|_{L^q(Q)}^q.
\end{equation}
Thus, if we define $E_3 \coloneq \bdE(\cQ_3)$, we conclude that
\begin{equation}\label{eq: pigeon 13}
   \|H^{\lambda} f\|_{L^q(B_R)}
   \lessapprox
   \Big\|\sum_{T\in\T}
   \sum_{D\in\mathcal{D}(T;\bar\kappa)}\big(H^{\lambda} f_T \big) \cdot \chi_D \Big\|_{L^q(E_3)},
\end{equation}
which follows from \eqref{eq: pigeon 1}, \eqref{eq: pigeon 2}, \eqref{eq: pigeon 3}, \eqref{eq: pigeon 8} and \eqref{eq: pigeon 12}. 




\subsection{Two-ends dichotomy} Recall that the statistic $\bar{\alpha}$ from P5 in \S\ref{subsec: pigeon} essentially corresponds to $\#\mathcal{D}(T, \bar{\kappa})$ for each tube $T \in \T$. Since the cubes $D \in \mathcal{D}(T, \bar{\kappa})$ have side-length $R^{1-\varepsilon_1}$, it follows that $1 \leq \bar{\alpha} \lesssim R^{\varepsilon_1}$. The proof splits into two cases depending on where $\bar{\alpha}$ lies in this range of values. 

To describe the dichotomy, we introduce a further parameter
\begin{equation}\label{eq: epsilon 2}
    \varepsilon_2 \coloneq \varepsilon^4 < \varepsilon_{\circ} \varepsilon^3.
\end{equation} 
Define a family of cubes
\begin{equation}\label{eq: shading 1}
 \bar{\cQ}
 \coloneq
\{Q \in \cQ_1: \exists\, T \in \T \textrm{ and } \exists\, D\in\mathcal{D}(T;\bar\kappa)\text{ s.t. } Q\cap T\cap D\neq\varnothing\}
\end{equation}
and, recalling \eqref{eq: Q_1 def}, define a shading $\bar{Y}$ on $\T$ by setting
\begin{equation}\label{eq: shading 2}
  \bar{\cQ}[T] \coloneqq 
\bigcup_{D\in\mathcal{D}(T;\bar\kappa)}\cQ_{1,D}[T] \quad \textrm{and} \quad    \bar{Y}(T)
\coloneqq\bigcup_{Q \in \bar{\cQ}[T]} Q,
\qquad \textrm{for all $T\in \T$.}   
\end{equation}
Thus, $\cQ_{\bar{Y}}[T] = \bar{\cQ}[T]$ for all $T \in \T$ and $\cQ_{\bar{Y}} = \bar{\cQ}$. Recall $\bar{\kappa}$ is the statistic introduced from P4 in \S\ref{subsec: pigeon}. Using the language of Definition~\ref{dfn: dense} and Definition~\ref{dfn: two ends}, the shaded family $(\T, \bar{Y})$ satisfies the following properties. 

\begin{lemma}\label{lem: shading} Let $(\T, \bar{Y})$ be as defined above.
\begin{enumerate}[a)]
    \item $(\T, \bar{Y})$ is $\rho$-dense for some $\rho \sim \bar{\kappa} \cdot \bar{\alpha}$.
    \item If $\bar{\alpha}\geq R^{\varepsilon_2}$, then $(\T, \bar{Y})$ is $(\varepsilon_1,\varepsilon_2)$-two-ends.
\end{enumerate}
\end{lemma}

In light of the above, we say
\begin{itemize}
    \item The \textit{concentrated case} holds if $1 \leq \bar{\alpha} < R^{\varepsilon_2}$;
    \item The \textit{two-ends case} holds if $ \bar{\alpha} \geq R^{\varepsilon_2}$.
\end{itemize}
In the concentrated case the operator is essentially local at scale $R^{1-\varepsilon_1}$. This allows a rather direct appeal to the induction hypothesis \eqref{eq: induction 3}. The two-ends case contains the real difficulties: it is here that the refined decoupling and Kakeya estimates play a role.

\begin{proof}[Proof (of Lemma~\ref{lem: shading})] Fix $T \in \T$. \medskip

\noindent a) It is clear from the definitions that the sets $\{\cQ_{1, D}[T] : D \in \mathcal{D}(T; \bar{\kappa})\}$ form a partition of $\bar{\cQ}[T]$. By \eqref{eq: pigeon 5}, we have
\begin{equation*}
    |\bar{Y}(T)|=\sum_{D\in\mathcal{D}(T;\bar\kappa)}|\bar Y_D(T)|\sim\#\mathcal{D}(T;\bar\kappa)\cdot|\bar Y_D(T)|\sim \bar{\alpha}\cdot\bar\kappa|T|.
\end{equation*} 
This establishes the first claim.\medskip

\noindent b) By the separation of the balls in $\mathcal{D}(T; \bar{\kappa})$, if $D_1$, $D_2 \in \mathcal{D}$ are distinct, then it follows that $d(\bar{Y}_{D_1}(T), \bar{Y}_{D_2}(T)) > 2R^{1-\varepsilon_1}$. Consequently, 
\begin{equation*}
   \max_{\substack{D \subseteq \R^n \\ \textrm{$R^{1-\varepsilon_1}$-ball}}} |\bar{Y}(T) \cap D| \leq \max_{D \in \mathcal{D}(T; \bar{\kappa})} |\bar{Y}_D(T)| \sim \bar{\kappa}|T|,
\end{equation*}
where in the last step we used the result of part a). On the other hand, 
\begin{equation*}
   |\bar{Y}(T)|\sim\bar{\alpha}\cdot\bar\kappa |T|,\quad\text{so that}\quad |\bar{Y}(T)\cap D|\lesssim \bar{\alpha}^{-1}|\bar{Y}(T)|. 
\end{equation*}
Assuming $\bar{\alpha} \geq R^{\varepsilon_2}$, the desired conclusion is now immediate. 
\end{proof}




\subsection{Passing to low multiplicity cubes} In the two-ends case, we must further refine the set of cubes $\cQ_3$ from \S\ref{subsec: pigeon} in order to control the multiplicity $\mu$ of tube-cube intersections. Here we require more quantitative information about the multiplicity $\mu$, so simple pigeonholing no longer suffices. 

Recall the set of cubes $\bar{\cQ}$ with shading $\bar{Y}$ defined in \eqref{eq: shading 1} and \eqref{eq: shading 2}, respectively. We now fix $\mu \geq 1$ and assume that it is chosen large enough so that
\begin{equation}\label{eq: Q mult 1}
    \#\big\{Q \in \bar{\cQ} : M_{\bar{Y}}(Q) >\mu \big\} \lesssim_{\varepsilon} R^{-\varepsilon_1} \#\bar{\cQ}.
\end{equation}
Later in the proof, in \S\ref{subsec: choose mu}, we shall use the multiplicity form of the hypothesised Kakeya estimates from \S\ref{sec: geo prelim} to determine explicit values of $\mu$ for which \eqref{eq: Q mult 1} holds. 

We further define 
\begin{equation}\label{eq: Q mult 2}
\cQ \coloneqq \big\{Q \in \bar{\cQ} : M_{\bar{Y}}(Q) \leq \mu \big\} \cap \cQ_3 \qquad \textrm{and} \qquad E \coloneq \bdE(\cQ).
\end{equation}
The condition \eqref{eq: Q mult 1} on $\mu$ ensures that $\cQ$ contains a large proportion of the cubes in $\cQ_3$. More precisely,  from P5 iii) and P4 iii) and the definition \eqref{eq: Q mult 1}, we have 
\begin{equation*}
    \#\cQ_3 \gtrapprox \#\cQ_2 \gtrapprox \#\cQ_1 \geq \#\bar{\cQ} \qquad \textrm{and so} \qquad \#\bar{\cQ} \lesssim_{\varepsilon} R^{\varepsilon_1/2} \#\cQ_3. 
\end{equation*}
Combining this with our hypothesis \eqref{eq: Q mult 1}, we deduce that there exists some constant $C_{\varepsilon} \geq 1$ such that 
\begin{equation*}
  \#  \big\{Q \in \bar{\cQ} : M_{\bar{Y}}(Q) >\mu \big\}  \leq C_{\varepsilon} R^{-\varepsilon_1/2} \#\cQ_3
\end{equation*}
Thus, provided $\bdR_{\varepsilon} \geq  1$ is chosen sufficiently large so as to ensure $C_{\varepsilon} R^{-\varepsilon_1/2} \leq 1/2$ for $R \geq \bdR_{\varepsilon}$, we have
\begin{equation}\label{eq: Q mult 3}
    \#\cQ \geq \# \cQ_3 - \#\big\{Q \in \bar{\cQ} : M_{\bar{Y}}(Q) >\mu \big\} \gtrsim \#\cQ_3. 
\end{equation}

Combining \eqref{eq: Q mult 3} with \eqref{eq: pigeon 13} and the second part of P5 ii) from \S\ref{subsec: pigeon}, we obtain 
\begin{equation}\label{eq: Q mult 4}
    \|H^{\lambda} f\|_{L^q(B_R)}\lessapprox \Big\|\sum_{T\in\T }\sum_{D\in\mathcal{D}(T;\bar\kappa)}\big(H^{\lambda} f_T\big) \cdot \chi_D \Big\|_{L^q(E)}.
\end{equation}
This (finally!) concludes the initial setup for the argument.




\section{Proof of Theorem~\ref{thm: main}: key estimates}




\subsection{Concentrated case via induction-on-scale} We begin with the concentrated case
\begin{equation}\label{eq: non 2 ends 1}
    1 \leq \bar{\alpha}<R^{\varepsilon_2}.
\end{equation}
This is settled by a direct application of the induction hypothesis. 
   
For each $D \in \mathcal{D}$, we define 
 \begin{equation}\label{eq: non 2 ends 2}
   f_D \coloneqq \sum_{T\in\T[D]}f_T \qquad \textrm{where} \qquad    \T[D]
       \coloneqq
       \{T\in\T:D\in\mathcal{D}(T;\bar\kappa)\}.
   \end{equation}
In view of \eqref{eq: pigeon 13}, by changing the order of summation we may bound 
\begin{equation*}
    \|H^{\lambda}f\|_{L^q(B_R)}^q
   \lessapprox\int_{E_3}\Big|\sum_{D\in\mathcal{D}} \big(H^{\lambda} f_D\big)\cdot \chi_D\Big|^q
  .
\end{equation*}
Since the cubes belonging to $\mathcal{D}$ are pairwise disjoint, it follows that
\begin{equation}\label{eq: non 2 ends 3}
\|H^{\lambda}f\|_{L^q(B_R)}^q \lessapprox
   \sum_{D\in\mathcal{D}}\|H^{\lambda} f_D\|_{L^q(D)}^q.
\end{equation}
    We wish to estimate the right-hand expressions by invoking the induction hypothesis. In order to do that, we recall that 
    $$\|H^\lambda f_D\|_{L^q(D)}=\|H^\lambda[\phi;a]\big(e^{-i\phi^\lambda(\bdx_B;\,\cdot\,)}f_D\big)\|_{L^q(\bdx_B+D)}.$$
    Thus, applying the induction hypothesis \eqref{eq: induction 3} at scale $R' \coloneq R^{1-\varepsilon_1}$ to each summand in \eqref{eq: non 2 ends 3} we obtain
    \begin{equation}
        \label{eq: non 2 ends 4}
\|H^{\lambda}f\|_{L^q(B_R)}^q
        \lessapprox
        \bdC_\varepsilon^q R^{(1-\varepsilon_1)\varepsilon q} 
\sum_{D\in\mathcal{D}}\tnorm{f_D}_{L_{R'}^p(\mathbb{B}^{n-1})}^q,
    \end{equation}
    noting that $\tnorm{e^{-i\phi^\lambda(\bdx_B;\,\cdot\,)}f_D}_{L^p_{R'}}=\tnorm{f_D}_{L^p_{R'}}$.
    We claim that
\begin{equation}\label{eq: non 2 ends 5}
     \sum_{D\in\mathcal{D}}\tnorm{f_D}_{L_{R'}^p(\mathbb{B}^{n-1})}^q \lesssim \bar{\alpha}^{q/p} \tnorm{f}_{L_{R}^p(\mathbb{B}^{n-1})}^q.
\end{equation}
Once \eqref{eq: non 2 ends 5} is established, we can combine it with \eqref{eq: non 2 ends 4} and the concentration hypothesis \eqref{eq: non 2 ends 1} to conclude that
\begin{equation*}
        \|H^{\lambda}f\|_{L^q(B_R)}^q
        \leq C_{\varepsilon}^q R^{(1 + q/p)\varepsilon_2-\varepsilon \varepsilon_1q}
        \bdC_\varepsilon^q R^{ \varepsilon q} 
         \tnorm{f}_{L_{ R}^p(\mathbb{B}^{n-1})}^q.
    \end{equation*}
for some constant $C_{\varepsilon} \geq 1$, which is independent of $R$. By our choice of $\varepsilon_1$ and $\varepsilon_2$ from \eqref{eq: epsilon 1} and \eqref{eq: epsilon 2}, we know that $\varepsilon \varepsilon_1q - (1+ q/p)\varepsilon_2 \geq \varepsilon^3(1-2/100n)$. Hence, provided $\bdR_{\varepsilon} \geq 1$ is chosen sufficiently large, depending only on $\varepsilon$, since $R > \bdR_{\varepsilon}$ we may ensure $C_{\varepsilon}^qR^{(1+ q/p)\varepsilon_2-\varepsilon \varepsilon_1q} \leq 1$. This closes the induction in this case.

It remains to prove the claim \eqref{eq: non 2 ends 5}. Recalling the definition of $f_D$ from \eqref{eq: non 2 ends 2}, by applying $L^2$ orthogonality from \eqref{eq: wp orthog 3} it follows that 
$$\sum_{D\in\mathcal{D}}\|f_D\|_{L^2(\R^{n-1})}^2\lesssim \sum_{D\in\mathcal{D}}\sum_{T \in \T[D]}\|f_T\|_{L^2(\R^{n-1})}^2=\sum_{T\in\T}\sum_{D\in\mathcal{D}(T;\bar\kappa)}\|f_T\|_{L^2(\R^{n-1})}^2.
$$
But by \eqref{eq: pigeon 9} and \eqref{eq: pigeon 11}, we know that $\#\mathcal{D}(T;\bar\kappa)\sim\bar{\alpha}$ for all $T\in\T$. Thus, by applying \eqref{eq: wp orthog 3} again, we have 
\begin{equation}\label{eq: non 2 ends 6}
\sum_{D\in\mathcal{D}}\|f_D\|_{L^2(\R^{n-1})}^2\lesssim \bar{\alpha}\,\|f\|_{L^2(\mathbb{B}^{n-1})}^2.
\end{equation}
 Recall from \S\ref{sec: wavepackets} that each $T \in \T$ has the form $T=T_{\theta,\omega}$, and that $f_{T_{\theta,\omega}}$ is supported in $2 \cdot \theta$. For each $\tau\in\Theta((R')^{-1/2})$, define 
\begin{equation*}
  \T[D;\tau]\coloneq\{T_{\theta,\omega}\in\T[D]:2\cdot\theta\cap\tau\neq\varnothing\}.   
\end{equation*}
Then, by the $L^2$ orthogonality as in \eqref{eq: wp orthog 3}, we have
\begin{equation*}
\|f_D\|_{L^2(\tau)}^2=\Big\|\sum_{T\in\T[D;\tau]}f_T\Big\|_{L^2(\tau)}^2 \lesssim\sum_{T\in\T[D;\tau]}\|f_T\|_{L^2(\R^{n-1})}^2\leq\sum_{\substack{\theta \in \Theta(R^{-1/2}) \\ 2\cdot\theta\cap\tau\neq\varnothing}}\sum_{T\in\bdT^R[\theta]}\|f_T\|_{L^2(\R^{n-1})}^2.  
\end{equation*}
On the other hand, for each fixed $\theta \in \Theta(R^{-1/2})$, we can apply $L^2$ orthogonality in the form \eqref{eq: wp orthog 2} to deduce 
\begin{equation*}
    \sum_{T\in\bdT^R[\theta]}\|f_T\|_{L^2(\R^{n-1})}^2 \sim \|f_\theta\|_{L^2(\theta)}^2.
\end{equation*}
It follows immediately that
\begin{equation}\label{eq: non 2 ends 7}
    \|f_D\|_{L_{\avg}^2(\tau)}^2 \lesssim\sup_{2\cdot\theta \cap\tau\neq\varnothing} \|f_{\theta}\|_{L_{\avg}^2(\theta)}^2
    \leq \sup_{\theta\in\Theta(R^{-1/2})} \|f\|_{L_{\avg}^2(\theta)}^2.
\end{equation}
Using the nesting of $\ell^p$ and applying \eqref{eq: non 2 ends 6} and \eqref{eq: non 2 ends 7} to definition of $\tnorm{\cdot}_{L_R^{p}}$, we obtain
\begin{equation*}
\sum_{D\in\mathcal{D}}\tnorm{f_D}_{L_{R'}^p(\B^{n-1})}^q \leq \Big(\sum_{D\in\mathcal{D}}\tnorm{f_D}_{L_{R'}^p(\B^{n-1})}^p \Big)^{q/p} \lesssim \bar{\alpha}^{q/p} \tnorm{f}_{L_R^p(\B^{n-1})}^q,
\end{equation*}
which is precisely the claimed bound \eqref{eq: non 2 ends 5}.



\subsection{Two-ends case: decoupling}
We now turn to the two-ends case. That is, we assume 
\begin{equation}\label{eq: 2 ends hyp}
    R^{\varepsilon_2} \leq \bar{\alpha} \lesssim R^{\varepsilon_1}. 
\end{equation}
Here we estimate the operator directly, without use of the induction hypothesis, by a mixture of refined decoupling, geometric estimates and $L^2$ theory.

 Recall that $E = \bdE(\cQ) \subseteq \R^n$ is the union of lattice $R^{1/2}$-cubes introduced in \eqref{eq: Q mult 2}. Observe the following:
\begin{enumerate}[i)]
    \item From \S\ref{subsec: pigeon}, we see that $\T \subseteq \T_3 \subseteq \T_2$. Thus, by P1 i), it follows that
    \begin{equation*}
        \|H^{\lambda} f_{T_1}\|_{L^{r_{n, \sigma}}(w_{B_R})} \sim \|H^{\lambda} f_{T_2}\|_{L^{r_{n, \sigma}}(w_{B_R})} \qquad \textrm{for all $T_1$, $T_2 \in \T$;}
    \end{equation*}
    \item By the definitions \eqref{eq: Q mult 2} and \eqref{eq: mult fn}, we have $M_{\bar{Y}}(Q) \lesssim \mu$ for all $Q \in \cQ$.
\end{enumerate}
This can be used to show that $\T$ and $\cQ$ satisfy the hypotheses of the refined decoupling inequality from Theorem~\ref{thm: refined dec}.\footnote{There is a slight discrepancy, since $M_{\bar{Y}}(Q)$ is the \textit{restricted} multiplicity, rather than the multiplicity $M_{\T}(Q)$; however, we shall deal with this minor point in the proof below.} For $r_{n, \sigma}$ the decoupling exponent from \eqref{eq: dec exp}, we thereby arrive at the following estimate. 

\begin{proposition}\label{prop: 2 ends dec} With the above setup,
\begin{multline*}
\Big\|\sum_{T\in\T }\sum_{D\in\mathcal{D}(T;\bar\kappa)} \big(H^{\lambda} f_T\big)\cdot \chi_D\Big\|_{L^{r_{n, \sigma}}(E)}\\
\lessapprox
R^{n\varepsilon_1}
\mu^{1/2 - 1/r_{n, \sigma}} R^{d(n,\sigma)} \max_{D \in \mathcal{D}}\Big(\sum_{T \in \T[D]}\|H^{\lambda} f_T\|_{L^{r_{n, \sigma}}(w_{B_R})}^{r_{n, \sigma}}\Big)^{1/r_{n, \sigma}}.
\end{multline*}  
\end{proposition}

\begin{proof} Recalling the definitions in \eqref{eq: non 2 ends 2}, by a change of order of summation,
\begin{equation*}
\Big\|\sum_{T\in\T }\sum_{D\in\mathcal{D}(T;\bar\kappa)} \big(H^{\lambda} f_T\big)\cdot \chi_D\Big\|_{L^{r_{n, \sigma}}(E)} = \Big\|\sum_{D\in\mathcal{D}} \big(H^{\lambda} f_D\big)\cdot \chi_D\Big\|_{L^{r_{n, \sigma}}(E)}.
\end{equation*}
Since $\mathcal{D}$ is a collection of $O(R^{n \varepsilon_1})$ pairwise disjoint cubes, it follows that
\begin{equation*}
\Big\|\sum_{T\in\T}\sum_{D\in\mathcal{D}(T;\bar\kappa)} \big(H^{\lambda} f_T\big)\cdot \chi_D\Big\|_{L^{r_{n, \sigma}}(E)} \lesssim R^{n \varepsilon_1} \max_{D\in \mathcal{D}} \big\|H^{\lambda} f_D\big\|_{L^{r_{n, \sigma}}(E_D)}
\end{equation*}
where $E_D \coloneq \bdE(\cQ_D)$ for $\cQ_D \coloneq \{Q \in \cQ : Q \cap D \neq \varnothing\}$. 

Fixing $D \in \mathcal{D}$, we wish to apply refined decoupling over $E_D$. First, observe that 
\begin{equation}
    \label{eq: res to unres}
    M_{\T[D]}(Q)\leq M_{\bar{Y}}(Q)\qquad\textrm{for all $Q \in \cQ_D$.}
\end{equation}
Here $M_{\T[D]}(Q)$ is the multiplicity as defined in \eqref{eq: mult fn} and $M_{\bar{Y}}(Q)$ is the restricted multiplicity as defined in \eqref{eq: res mult fn}. To see why \eqref{eq: res to unres} holds, suppose $Q \in \cQ_D$ and $T \in \T[D]$ satisfy $Q \cap T \neq \varnothing$. From the definitions \eqref{eq: Q mult 2} and \eqref{eq: shading 1}, we see that $\cQ \subseteq \cQ_1$. Furthermore, it follows from the definition \eqref{eq: Q_1 def} that $Q\in\cQ_{1,D}[T]$. On the other hand, $T\in\T[D]$ implies that $D\in\mathcal{D}(T;\bar\kappa)$, and so by \eqref{eq: shading 2}, it follows that $Q \in \bar{\cQ}[T] = \cQ_{\bar{Y}}[T]$.

In light of \eqref{eq: res to unres}, and the discussion preceding the proposition, we may apply Theorem~\ref{thm: refined dec} on each set $E_D$ and thereby bound
\begin{equation*}
\Big\|\sum_{T\in\T[D]}H^{\lambda} f_T\Big\|_{L^{r_{n, \sigma}}(E_D)} \lessapprox \mu^{1/2 - 1/r_{n, \sigma}} R^{d(n,\sigma)} \Big(\sum_{T \in \T[D]}\|H^{\lambda} f_T\|_{L^{r_{n, \sigma}}(w_{B_R})}^{r_{n, \sigma}}\Big)^{1/r_{n, \sigma}}.
\end{equation*}
This concludes the proof.
\end{proof}
To estimate the $\ell^{r_{n,\sigma}}$-sum, we define
\begin{equation*}
    \fkD
\coloneqq
\{\theta\in\Theta(R^{-1/2}):\T_3[\theta]\neq\varnothing\}.
\end{equation*}
\begin{proposition}\label{prop: ell r n estimate}
    With the above setup, for any $\vartheta\in [0,1]$, we have
    \begin{multline}\label{eq: ell r n estimate}
    \Big(\sum_{T \in \T[D]}\|H^{\lambda} f_T\|_{L^{r_{n,\sigma}}}^{r_{n, \sigma}}\Big)^{1/r_{n, \sigma}}
    \lesssim    
    (\#\T_3)^{\frac{1}{r_{n,\sigma}}-\frac12}
    R^{\frac{n+1}{2}(\frac{1}{r_{n,\sigma}}-\frac12)}
    (\#\fkD)^{\frac{1}{2}(1-\vartheta)}R^{\frac{n-1}{2}(\frac{1}{2/\vartheta}-\frac12)}R^{\frac{1}{2}}\tnorm{f}_{L_{R}^{2/\vartheta}}.
    \end{multline}
    Here, the $L^{r_n,\sigma}$ norm is taken against the weight $w_{B_R}$.
\end{proposition}

\begin{proof}
    We relax the summation in the expression of the left-hand side of \eqref{eq: ell r n estimate} to a sum over all $T\in\T_3$. Using the comparability property P1 i), we have
    \begin{equation*}
       \Big(\sum_{T \in \T_3}\|H^{\lambda} f_T\|_{L^{r_{n, \sigma}}(w_{B_R})}^{r_{n, \sigma}}\Big)^{1/r_{n, \sigma}}
    \sim 
    (\#\T_3)^{-(\frac{1}{2}-\frac{1}{r_{n,\sigma}})}\Big(\sum_{T\in\T_3}\|H^{\lambda}f_T\|_{L^{r_{n,\sigma}}(w_{B_R})}^{2}\Big)^{\frac{1}{2}}. 
    \end{equation*}
    By Lemma \ref{lem: Bernstein}, we can bound
    \begin{equation*}
\Big(\sum_{T\in\T_3}\|H^{\lambda}f_T\|_{L^{r_{n,\sigma}}(w_{B_R})}^{2}\Big)^{\frac{1}{2}}
    \lesssim 
    R^{\frac{1}{2}-\frac{n+1}{2}(\frac{1}{2}-\frac{1}{r_{n,\sigma}})}\Big(\sum_{T \in \T_3}\|f_T\|_{L^2(\R^{n-1})}^2\Big)^{1/2}.
    \end{equation*}
    Recall we have basic $L^2$ orthogonality as in \eqref{eq: wp orthog 3}, giving
    \begin{equation}\label{eq: off diag many directions}
        \Big(\sum_{T\in \T_3}\|f_T\|_{L^2(\R^{n-1})}^2\Big)^{1/2}
        \lesssim
        \|f\|_{L^2(\R^{n-1})}.
        \end{equation}
We use \eqref{eq: wp orthog 2} and the definition of $f_{\theta}$ to estimate 
\begin{equation}\label{eq: off diag few directions case I}
        \Big(\sum_{T\in \T_{3}}\|f_T\|_2^2\Big)^{1/2}
        \lesssim 
        \Big(\sum_{\theta\in\fkD}\|f_{\theta}\|_2^2\Big)^{1/2}
        \lesssim
        (\#\fkD)^{\frac{1}{2}}R^{-\frac{n-1}{4}}\sup_{\theta\in\Theta(R^{-1/2})}\|f\|_{L_{\avg}^2(\theta)}.
        \end{equation}
For $\vartheta\in [0,1]$, we take the $\vartheta$-geometric mean of \eqref{eq: off diag many directions} and \eqref{eq: off diag few directions case I}.
Noting that 
    $$
    R^{-\frac{n-1}{4}(1-\vartheta)}
    =
    R^{-\frac{n-1}{2}(\frac{1}{2}-\frac{1}{2/\vartheta})},$$
    we deduce the desired result.\qedhere
\end{proof}

\begin{corollary}\label{cor: two ends L r n}
    With the above setup, for every $\vartheta\in [0,1]$, if $a\coloneqq 2/\vartheta$, we have
    \begin{multline*}
    \Big\|\sum_{T\in\T }\sum_{D\in\mathcal{D}(T;\bar\kappa)} \big(H^{\lambda} f_T\big)\cdot \chi_D\Big\|_{L^{r_{n, \sigma}}(E)}\\
\lessapprox
R^{n\varepsilon_1}
    (\mu^{-1}\#\T_3)^{-(\frac{1}{2}-\frac{1}{r_{n,\sigma}})}
    (\#\fkD)^{\frac{1}{2}(1-\vartheta)}R^{-\frac{n-1}{2}(\frac{1}{2}-\frac{1}{a})}\tnorm{f}_{L_{R}^{a}(\mathbb{B}^{n-1})}.
    \end{multline*}    
    
\end{corollary}

\begin{proof} Since 
\begin{equation*}
  d(n,\sigma)-\frac{n+1}{2}\Big(\frac{1}{2}-\frac{1}{r_{n,\sigma}}\Big)+\frac{1}{2}=0,  
\end{equation*}
this follows by combining Proposition \ref{prop: 2 ends dec} with Proposition \ref{prop: ell r n estimate}.
\end{proof}




\subsection{\texorpdfstring{$L^2$}{} theory} Recall that $\rho \sim \bar{\kappa} \bar{\alpha}$ denotes the density of the shaded family $(\T, \bar{Y})$, as determined in Lemma~\ref{lem: shading}. As in the previous subsection, $E = \bdE(\cQ) \subseteq \R^n$ is the union of $R^{1/2}$-cubes introduced in \eqref{eq: Q mult 2}.

\begin{lemma}[{\textit{cf.} \cite[Lemma 4.5]{Wang_Wu}, \cite[Lemma 2.11]{MKN}}]\label{lem: 2 ends L2} With the above setup, 
     \begin{equation*}
    \Big\|\sum_{T\in\T }\sum_{D\in\mathcal{D}(T;\bar\kappa)} \big(H^{\lambda} f_T\big)\cdot \chi_D \Big\|_{L^2(E)}
     \lesssim
     (\rho R)^{1/2} R^{O(\delta)}\tnorm{f}_{L_{R}^2(\mathbb{B}^{n-1})}.
     \end{equation*}
\end{lemma}

\begin{proof} From the definitions \eqref{eq: Q mult 2} and \eqref{eq: shading 1}, we see that $\cQ \subseteq \cQ_1$. For $Q \in \cQ_1$, let $D_Q$ denote the unique $D \in \mathcal{D}$ satisfying $Q \cap D \neq \varnothing$ and set
\begin{equation*}
    f_Q \coloneqq  \sum_{T \in \T[Q]} f_T \quad \textrm{where} \quad \T[Q]
\coloneqq 
\{T \in \T: T \cap Q \neq \varnothing \textrm{ and } D_Q\in\mathcal{D}(T;\bar\kappa)\}.
\end{equation*}
By the spatial concentration property of the wavepackets, 
\begin{equation}\label{eq: key 2 0} 
\int_Q\Big|\sum_{T\in\T}\sum_{D\in\mathcal{D}(T;\bar\kappa)}\big(H^{\lambda} f_T\big)\cdot\chi_D\Big|^2 \lesssim_{\varepsilon} \int_Q\Big|\sum_{\substack{T\in\T \\T \cap Q \neq \varnothing}}\sum_{D\in\mathcal{D}(T;\bar\kappa)}\big(H^{\lambda} f_T\big)\cdot\chi_D\Big|^2 + R^{-100n} \|f\|_{L^2(\mathbb{B}^{n-1})}^2.
\end{equation}

Focusing on the first term on the right-hand side of \eqref{eq: key 2 0}, by interchanging the order of summation,
\begin{equation}\label{eq: key 2 1} \int_Q\Big|\sum_{\substack{T\in\T \\T \cap Q \neq \varnothing}}\sum_{D\in\mathcal{D}(T;\bar\kappa)}\big(H^{\lambda} f_T\big)\cdot\chi_D\Big|^2 \leq \|H^{\lambda} f_Q\|_{L^2(Q)}^2.
\end{equation}
By H\"ormander's $L^2$-bound \cite{hormander} at scale $R^{1/2}$ and \eqref{eq: wp orthog 3}, we have
\begin{equation}\label{eq: key 2 2}
    \|H^{\lambda} f_Q\|_{L^2(Q)}^2
\lesssim R^{1/2}\|f_Q\|_{L^2(\R^{n-1})}^2
\lesssim
R^{1/2}\sum_{T\in \T[Q]}\|f_T\|_{L^2(\R^{n-1})}^2.
\end{equation}
We apply \eqref{eq: key 2 2} to the right-hand side of \eqref{eq: key 2 1} and sum the resulting bound over all $Q\in \cQ$. A simple argument using \eqref{eq: Q mult 4} and the non-degeneracy hypothesis \eqref{eq: non deg} allows us to absorb the rapidly decreasing terms from \eqref{eq: key 2 0} into the left-hand side of the resulting inequality, giving 
\begin{equation}\label{eq: key 2 3}
  \Big\|\sum_{T\in\T}\sum_{D\in\mathcal{D}(T;\bar\kappa)}\big(H^{\lambda} f_T\big)\cdot \chi_D\Big\|_{L^2(E)}^{2}
\lesssim 
R^{1/2}\sum_{Q\in\cQ}\sum_{T\in\T[Q]} \|f_T\|_{L^2(\R^{n-1})}^2. 
\end{equation}

Recall from \eqref{eq: shading 2} that
\begin{equation*}
    \bar{\cQ}[T] \coloneq \big\{Q \in \cQ_1 : \exists\, D \in \mathcal{D}(T;\bar\kappa) \textrm{ s.t. } D\cap Q \cap T \neq\varnothing  \big\} = \big\{Q \in \cQ_1 : T \in \T[Q]\big\}
\end{equation*}
for all $T \in \T$. We relax the summation in \eqref{eq: key 2 3} to the larger family $\cQ_1$ and interchange the order of summation to obtain
\begin{equation}
    \label{eq: key 2 4}
    \Big\|\sum_{T\in\T }\sum_{D\in\mathcal{D}(T;\bar\kappa)}\big(H^{\lambda} f_T\big)\cdot \chi_D\Big\|_{L^2(E)}^{2}\lesssim R^{1/2}\sum_{T\in\T}\#\bar{\cQ}[T]\|f_T\|_{L^2(\R^{n-1})}^2.
\end{equation}
Recall from \eqref{eq: shading 2} that the shading $\bar{Y}(T)$ is the union of the cubes from $\bar{\cQ}[T]$. Thus, since we know from Lemma~\ref{lem: shading} that $(\T, \bar{Y})$ is a $\rho$-dense shaded family,
\begin{equation}\label{eq: key 2 5}
    \#\bar{\cQ}[T] = R^{-n/2} |\bar{Y}(T)| \lesssim \rho R^{-n/2}|T| \sim \rho R^{1/2 + O(\delta)}.
\end{equation}
Combining \eqref{eq: key 2 4} and \eqref{eq: key 2 5} and once again applying $L^2$ orthogonality as in \eqref{eq: wp orthog 3}, we have
\begin{equation*}
  \Big\|\sum_{T\in\T}\sum_{D\in\mathcal{D}(T;\bar\kappa)}\big(H^{\lambda} f_T\big)\cdot \chi_D\Big\|_{L^2(E)}^{2}
\lesssim 
\rho R^{1+ O(\delta)}\sum_{T\in\T } \|f_T\|_{L^2(\R^{n-1})}^2 \lesssim 
\rho R^{1  + O(\delta)}\|f\|_{L^2(\mathbb{B}^{n-1})}^2. 
\end{equation*}
This concludes the proof. 
\end{proof}




\subsection{Choosing the multiplicity: Kakeya estimates}\label{subsec: choose mu}

Throughout the analysis of the two-ends case, we have worked over the set $E$, which is defined in \eqref{eq: Q mult 2} with respect to some choice of multiplicity parameter $\mu$. This parameter must be taken to be sufficiently large so as to ensure \eqref{eq: Q mult 1} holds; that is,
\begin{equation}\label{eq: choose mu 1}
    \#\big\{ Q \in \cQ_{\bar{Y}} : M_{\bar{Y}}(Q) > \mu \big\} \lesssim_{\varepsilon} R^{-\varepsilon_1} \#\cQ_{\bar{Y}},
\end{equation}
where we have used the fact that $\bar{\cQ} = \cQ_{\bar{Y}}$. As a final ingredient in the proof, we apply the Kakeya theory from \S\ref{sec: geo prelim} to determine a suitable choice of $\mu$. 

The hypothesis \eqref{eq: 2 ends hyp}, together with Lemma~\ref{lem: shading}, implies that $(\T, \bar{Y})$ is $(\varepsilon_1, \varepsilon_2)$-two ends. By the definition \eqref{eq: pigeon 11} and P2 i) from \S\ref{subsec: pigeon}, we have 
\begin{equation*}
    \#\T \leq \sum_{\theta \in \Theta(R^{-1/2})} \#\T_3[\theta] \lesssim m \cdot \#\Theta(R^{-1/2}) \sim m R^{(n-1)/2},
\end{equation*}
Thus, Lemma~\ref{lem: univ two-ends} implies that
\begin{equation}\label{eq: choose mu 2}
     \#\{Q \in \cQ_{\bar{Y}} : M_{\bar{Y}}(Q) > \mu \} \lesssim_{\varepsilon} \mu^{-1} (\#\T)^{\frac{1}{2}} R^{O(\varepsilon_1)}  \#\cQ_{\bar{Y}} \qquad \textrm{for all $\mu > 0$.}
\end{equation}
On the other hand, we also know from Lemma~\ref{lem: shading} that $(\T, \bar{Y})$ is a $\rho$-dense shaded family. Furthermore, recalling $\T \subseteq \T_3$ from \eqref{eq: pigeon 11}, by P2 i) from \S\ref{subsec: pigeon} it follows that:\smallskip

For each fixed $\theta \in \Theta(R^{-1/2})$, there are $O(m)$ choices of $\omega \in R^{1/2}\Z^{n-1}$ such that $T_{\theta, \omega} \in \T$. We can therefore partition $\T$ into $O(m)$ subfamilies to which the hypothesised inequality $\bdK_{s \to s}^*(\phi; \beta)$ applies.
In light of the above, it follows from H\"older's inequality that \eqref{eq: improv mult 1} holds for $A = O(m^{1/s'}R^{O(\delta)})$. For $1 < s < \infty$, Lemma~\ref{lem: cond mult} then implies that
\begin{equation}\label{eq: choose mu 3}
     \#\{Q \in \cQ_{\bar{Y}} : M_{\bar{Y}}(Q) > \mu \} \lesssim \mu^{-1} m \rho^{-(s'-1)}R^{\frac{\beta s'}{2}+O(\varepsilon_1)} \#\cQ_{\bar{Y}}
\end{equation}
holds for all $\mu > 0$.

For $1 < s < \infty$, we may take the geometric mean of the estimates from \eqref{eq: choose mu 3} and \eqref{eq: choose mu 2} to give 
\begin{equation*}
    \#\{Q \in \cQ_{\bar{Y}} : M_{\bar{Y}}(Q) > \mu \} \lesssim_{\varepsilon} R^{C\varepsilon_1} \mu^{-1} (\#\T)^{\frac{1}{2}(1-\varphi)} m^{\varphi} R^{\frac{\beta s'}{2}\varphi} \rho^{-(s'-1)\varphi} \#\cQ_{\bar{Y}}
\end{equation*}
for any $0 \leq \varphi \leq 1$, for some dimensional constant $C \geq 1$. Therefore, if we define 
\begin{equation}\label{eq: choose mu 4 general case}
    \bar{\mu}(\varphi)
    \coloneq
    R^{(C+1)\varepsilon_1} (\#\T)^{\frac{1}{2}(1-\varphi)} m^{\varphi} R^{\frac{\beta s'}{2}\varphi}\rho^{-(s'-1)\varphi} \qquad \textrm{for $0 \leq \varphi \leq 1$}, 
\end{equation}
then \eqref{eq: choose mu 1} holds with $\mu = \bar{\mu}(\varphi)$. For $s = 1$, we always take $\varphi \coloneq 0$ and consider \eqref{eq: choose mu 2} only. That is, we let $\bar{\mu}(0)
    \coloneq R^{(C+1)\varepsilon_1} (\#\T)^{\frac{1}{2}}$ in this case.




\subsection{Interpolation and concluding the proof}
 
For a given $t\in [0,1]$, let $1 \leq p$, $q \leq \infty$ be such that
\begin{equation*}
    \frac{1}{q}=\frac{1-t}{2}+\frac{t}{r_{n,\sigma}},\qquad \frac{1}{p}=\frac{1-t}{2}+\frac{t}{a}.
\end{equation*}
Notice that 
\begin{equation*}
    \|f\|_{L^2(\mathbb{B}^{n-1})}^{1-t}\tnorm{f}_{L_{R}^{a}(\mathbb{B}^{n-1})}^t  = \tnorm{f}_{L_{R}^p(\mathbb{B}^{n-1})}.
\end{equation*}
Hence, using logarithmic convexity, we interpolate between the $L^2\to L^2$ estimate from Lemma~\ref{lem: 2 ends L2} and the $L^{a}\to L^{r_{n,\sigma}}$ estimate from Corollary~\ref{cor: two ends L r n} to obtain
\begin{multline*}
    \Big\|\sum_{T\in\T }\sum_{D\in\mathcal{D}(T;\bar\kappa)} \big(H^{\lambda} f_T\big)\cdot \chi_D\Big\|_{L^{q}(E)}\\
    \lesssim_{\varepsilon}
    \rho^{\frac{(1-t)}{2}}
    (\mu^{-1}\#\T_3)^{-(\frac{1}{2}-\frac{1}{r_{n,\sigma}})t}
    (\#\fkD)^{\frac{1}{2}(1-\vartheta)t}
    R^{-\frac{n-1}{2}(\frac{1}{2}-\frac{1}{a})t+\frac{1}{2}(1-t)}\tnorm{f}_{L_{R}^{p}(\mathbb{B}^{n-1})}.
\end{multline*}    
Recall from \S\ref{subsec: choose mu} that $\mu = \bar{\mu}(\varphi)$, as defined in \eqref{eq: choose mu 4 general case}, is an admissible multiplicity for any $\varphi \in [0,1]$, in the sense that this choice guarantees \eqref{eq: Q mult 1} holds. Moreover, taking $\mu = \bar{\mu}(\varphi)$, the right-hand side of the above estimate is bounded by
\begin{multline}\label{eq: complicated RHS}
\rho^{\frac{1-t}{2}-(s'-1)\varphi t(\frac{1}{2}-\frac{1}{r_{n,\sigma}})}
[(\#\T_3)^{-\frac{1}{2}(1+\varphi)(\frac{1}{2}-\frac{1}{r_{n,\sigma}})t}
m^{\varphi t(\frac{1}{2}-\frac{1}{r_{n,\sigma}})}
(\#\fkD)^{\frac{1}{2}(1-\vartheta)t}]
\\R^{O(\varepsilon_1)+t[\frac{\beta s'}{2}(\frac{1}{2}-\frac{1}{r_{n,\sigma}})\varphi-\frac{n-1}{2}(\frac{1}{2}-\frac{1}{a})]+\frac{1}{2}(1-t)}\tnorm{f}_{L_{R}^p(\mathbb{B}^{n-1})}, 
\end{multline}
having used that $R^{(C+1)\varepsilon_1 t(\frac{1}{2}-\frac{1}{r_{n,\sigma}})}=R^{O(\varepsilon_1)}$.

Since $\rho\in (0,1)$, we need its exponent to be non-negative. This happens if and only if
\begin{equation}\label{eq: condition rho}
\varphi \leq \frac{s-1}{2}\Big(\frac{1}{t}-1\Big)\Big(\frac{1}{2}-\frac{1}{r_{n,\sigma}}\Big)^{-1}.    
\end{equation}
Using that $\#\T_3 \sim m\cdot (\#\fkD)$, the expression in square brackets in \eqref{eq: complicated RHS} becomes
\begin{equation*}
 m^{(\varphi-\frac{1}{2}(1+\varphi))(\frac{1}{2}-\frac{1}{r_{n,\sigma}})t}
(\#\fkD)^{(\frac{1}{2}(1-\vartheta)-\frac{1}{2}(1+\varphi)(\frac{1}{2}-\frac{1}{r_{n,\sigma}}))t}.   
\end{equation*}
Since $\varphi-(1+\varphi)/2=(\varphi-1)/2$, the exponent of $m$ is non-positive if and only if $\varphi\leq 1$, which already holds by the setup.

The exponent of $\#\fkD$ is non-positive if and only if 
\begin{equation*}
   \vartheta \geq 1-(1+\varphi)(\frac{1}{2}-\frac{1}{r_{n,\sigma}})
=
\frac{(n+\sigma+3)-2(1+\varphi)}{n+\sigma+3}. 
\end{equation*}

Since $a=2/\vartheta$, this is equivalent to
\begin{equation}\label{eq: condition D iota}
    a\leq \frac{2(n+\sigma+3)}{(n+\sigma+3)-2(1+\varphi)}.
\end{equation}

The exponent of $R$ in \eqref{eq: complicated RHS} can be written as $Y + O(\varepsilon_1)$ where
\begin{align*}
    Y &\coloneq
t\Big[\frac{\beta s'}{2}\Big(\frac{1}{2}-\frac{1}{r_{n,\sigma}}\Big)\varphi-\frac{n-1}{2}\Big(\frac{1}{2}-\frac{1}{a}\Big)\Big]+\frac{1}{2}(1-t) \\
&=
t\Big[\frac{2\beta s'\varphi}{2(n+\sigma+3)}-\frac{n-1}{2}\Big(\frac{1}{2}-\frac{1}{a}\Big)-\frac{1}{2}\Big]+\frac{1}{2}.
\end{align*}
This is non-positive if and only if
\begin{equation}\label{eq: condition R}
t\geq \Big[1+(n-1)\Big(\frac{1}{2}-\frac{1}{a}\Big)-\frac{2\beta s'\varphi}{n+\sigma+3}\Big]^{-1}.    
\end{equation}
Choosing $t=t_{n,\sigma}(s,\beta)$ as in \eqref{eq: target exp arrival} and $a=a_{n,\sigma}(s,\beta)$ and $B=B_{n,\sigma}(s,\beta)$ as in \eqref{eq: a n sigma p beta}, it can be checked that \eqref{eq: condition rho}, \eqref{eq: condition D iota} and \eqref{eq: condition R} all hold with equality. Hence, 
\begin{equation*}
   \Big\|\sum_{T\in\T }\sum_{D\in\mathcal{D}(T;\bar\kappa)} \big(H^{\lambda} f_T\big)\cdot \chi_D\Big\|_{L^{q}(E)}\\
    \lesssim_{\varepsilon}
    R^{O(\varepsilon_1)}\tnorm{f}_{L_{R}^p(\mathbb{B}^{n-1})}. 
\end{equation*}

Recall from \eqref{eq: epsilon 1} that $0 < \varepsilon_1 < \varepsilon_{\circ} \varepsilon$, where $\varepsilon_{\circ} > 0$ is a small dimensional constant which we are free to choose. Recalling \eqref{eq: Q mult 4}, provided $\bdC_{\varepsilon} \geq 1$ is chosen sufficiently large and $\varepsilon_{\circ} > 0$ is chosen sufficiently small, we conclude that
\begin{equation*}
    \|H^{\lambda} f\|_{L^q(B_R)} \leq \bdC_{\varepsilon} R^{\varepsilon} \tnorm{f}_{L_{R}^p(\mathbb{B}^{n-1})},
\end{equation*}
which establishes the desired bound in the two-ends case. This closes the induction and completes the proof. \qed




\section{Proof of Corollary~\ref{cor: beyond univ}}\label{subsec: cor}

The objective of this section is to prove the following:

\begin{proposition}\label{prop: improved uniform geometric bound}
    Under the hypotheses of Corollary~\ref{cor: beyond univ}, for any $1 < s  < \frac{n+1}{n-1}$ there exists $\beta_0(s)<\beta(n,s)$ such that $\bdK_{s \to s}^*(\phi;\beta_0(s))$ holds.
\end{proposition}

Indeed, using this result, we can prove Corollary~\ref{cor: beyond univ}. First, notice the following relations between the exponents.

\begin{lemma}\label{lem: relation exponents 0}
Let $1 < s < \frac{n+1}{n-1}$. If $\beta<\beta(n,s)$, then $q_{n,\sigma}(s,\beta)<q_{n,\sigma}$.        
\end{lemma}

\begin{proof}
    Since $\beta(n, s) > \beta$, it is clear from \eqref{eq: A exp} that $A(s, \beta) > 0$. Since $t_{n,\sigma}\in (0,1)$, then $1-A(s,\beta)<1-t_{n,\sigma}A(s,\beta)$, so from \eqref{eq: u exp} we see that $u_{n, \sigma}(s, \beta)<1$. Hence, it follows from \eqref{eq: target exp arrival} that $u_{n,\sigma}(s,\beta) < 1$ implies that $q_{n,\sigma}(s,\beta)<q_{n,\sigma}$, as required.    \qedhere    
\end{proof}

\begin{lemma}\label{lem: relation exponents 1}
Let $1 \leq s \leq \frac{n+1}{n-1}$ and $\beta\geq 0$. Then $q_{n,\sigma}(s,\beta)\geq p_{n,\sigma}(s,\beta)$.    
\end{lemma}

\begin{proof}
    By definition,
    \[
    q_{n,\sigma}(s,\beta)\geq p_{n,\sigma}(s,\beta)
    \iff 
    r_{n,\sigma}\geq a_{n,\sigma}(s,\beta).
    \]
    Since we can express
    \[
    a_{n,\sigma}(s,\beta)
    =
    \frac{r_{n,\sigma}}{1+\frac{2}{n+\sigma-1}[1-B(s,\beta)]},
    \]
    it is enough to prove that $B(s,\beta)\leq 1$. Using that $s\leq \frac{n+1}{n-1}$, we get 
    \[
    B(s,\beta)
    =
    \frac{s-1}{4}\frac{n-1}{1-A(s,\beta)}
    \leq
    \frac{1}{2}\frac{1}{1-A(s,\beta)}.
    \]
    And 
    \begin{equation}\label{eq: relation exponents 1 1}
    \frac{1}{2}\frac{1}{1-A(s,\beta)}<1
    \iff 
    \frac{s}{2}[\beta(n,s)-\beta]<\frac{1}{2}.
    \end{equation}
    Since $\beta\geq 0$ and $s\leq \frac{n+1}{n-1}$, we obtain
    \[
    \frac{s}{2}\beta(n,s)
    =
    \frac{n-1}{4}(s-1)
    \leq 
    \frac{n-1}{4}\Big(\frac{n+1}{n-1}-1\Big)
    =
    \frac{1}{2}.
    \]
    Hence, \eqref{eq: relation exponents 1 1} holds, so $B(s,\beta)\leq 1$ and the result is proved.
\end{proof}

\begin{proof}[Proof of Corollary~\ref{cor: beyond univ}]
    Using Theorem \ref{thm: main} together with Proposition~\ref{prop: improved uniform geometric bound}, we obtain $\bdH_{p_{\phi} \to q_{\phi}}(\phi)$ for $p_{\phi}\coloneqq p_{n,\sigma}(s,\beta_0(s))$ and $q_{\phi}\coloneqq q_{n,\sigma}(s,\beta_0(s))$. By Lemma~\ref{lem: relation exponents 0} and Lemma~\ref{lem: relation exponents 1}, $p_{\phi}<q_{\phi}<q_{n,\sigma}$. Reasoning as in the proof of Corollary~\ref{cor: univ off diag} by the finite-measure support, we obtain that $\bdH_{q_{\phi} \to q_{\phi}}(\phi)$ holds. Interpolating with the trivial bound $\bdH_{\infty \to \infty}(\phi)$, we conclude the proof.
\end{proof}

Hence, it remains to prove Proposition \ref{prop: improved uniform geometric bound}.




\subsection{Reformulation of the geometric estimates}

For the rest of this section, we assume the hypotheses of Corollary~\ref{cor: beyond univ} and we shall fix $1 < s  < \frac{n+1}{n-1}$ arbitrarily.

For $K \geq 1$, define $\phi^K(x,t;y) \coloneq K\phi(x/K,t/K;y)$ as in \eqref{eq: rescaled}. Given $\delta > 0$, we say a set $T \subseteq \R^n$ is a \textit{$1 \times \delta$ tube associated to $\phi^K$} if it is of the form
\begin{equation*}
T = T_{y, \omega} = \big\{(x,t)\in \B^{n-1}\times \B^1:|\partial_y\phi^K(x,t;y)-\omega|< \delta \big\}
\end{equation*}
for some $y \in \delta \Z^{n-1}\cap\B^{n-1}$ and $\omega\in \delta \Z^{n-1}$. Given $1 \leq s \leq \infty$, $1 \leq R \leq \lambda$ and $\beta \geq 0$, $\varepsilon > 0$, after rescaling, we are interested in inequalities of the form
\begin{equation}
    \label{eq: rescaled Kakeya}
\big\|\sum_{T\in\T}\chi_T\big\|_{L^s(\R^n)}\lesssim_{\phi,\eps} R^{(n+1)/(2s)+\beta/2+\eps}R^{-n/s}[\#\T]^{1/s},
\end{equation}
where $\T$ is a set of $1 \times R^{-1/2}$ tubes associated to $\phi^{\lambda/R}$. In particular, for fixed $1 \leq s \leq \infty$ and $\beta \geq 0$, note: 

$\bdK_{s \to s}^*(\phi; \beta)$ holds provided \eqref{eq: rescaled Kakeya} is valid for all $\varepsilon > 0$ whenever $\T$ is direction-separated in the sense that $|y_1 - y_2| \geq R^{-1/2}$ whenever $T_{y_1, \omega_1}$, $T_{y_2, \omega_2} \in \T$ are distinct.

To prove  \eqref{eq: rescaled Kakeya}, we consider two different parameter regimes. In particular, let $\varepsilon_0 > 0$ be a fixed constant, chosen sufficiently small for the forthcoming purposes of the proof. We shall argue differently depending on whether $\lambda \leq R^{1 + \varepsilon_0}$ or $R^{1 + \varepsilon_0} < \lambda$.  




\subsection{Large \texorpdfstring{$R$}{} regime: \texorpdfstring{$R \leq \lambda \leq R^{1+\varepsilon_0}$}{}} Using a rescaling argument, Theorem~2.9 of \cite{CGGHIW} (after using the duality relation \cite[(6.2)]{CGGHIW}) immediately establishes the desired Kakeya estimates in this case, provided $\eps_0 > 0$ is chosen sufficiently small. 

Rescaling by a factor of $R/\lambda$ and replacing each resulting tube by a superset living at unit scale, to show \eqref{eq: rescaled Kakeya} it suffices to prove
\begin{equation*}
    \big\|\sum_{T\in \T}\chi_{\wt{T}}\big\|_{L^s(\R^n)}\lesssim_{\phi,\eps}R^{(n+1)/(2s)+\beta/2+\eps}\lambda^{-n/s}[\#\T]^{1/s},
\end{equation*}
where
\begin{equation*}
 \wt{T}\coloneq\big\{(x,t)\in \B^{n-1}\times \B^1:|\partial_y\phi(x,t;y)-(R/\lambda) \omega|<R^{1/2}/\lambda\big\} \quad \textrm{for $T = T_{y, \omega} \in \T$.}  
\end{equation*}
 Let $\wt{\T} \coloneq \{\wt{T} : T \in \T\}$, so that $\wt{\T}$ is a family of $1 \times R^{1/2}/\lambda$ tubes. If $\T$ is direction-separated, then $\wt{\T}$ is direction-separated, since $|y_1 - y_2|\geq R^{-1/2}$ implies $|y_1 - y_2|\geq R^{1/2}/\lambda$ for $y_1$, $y_2 \in \R^{n-1}$.

\begin{theorem}[\textit{cf.} {\cite[Theorem 2.9]{CGGHIW}}]\label{thm: improv geometric bounds}
    Let $n\geq 3$ be odd and $\phi\colon \D^n\to \R$ be a real analytic H\"ormander-type phase. For $\delta > 0$ suppose that $\phi$ is translation-invariant and satisfies Hypothesis I) and $\T$ is a direction-separated family of $1\times\delta$ tubes  associated to $\phi$. 
    For every $1<s<\frac{n+1}{n-1}$, there exists $\beta_{\phi}(s)<\beta(n,s)$ such that 
    \begin{equation}\label{eq: delta multiplicity bound}
\big\|\sum_{T\in\T}\chi_{T}\big\|_{L^s(\R^n)}
    \lesssim_{\phi}
   \delta^{(n-1)/s-\beta_{\phi}(s)}[\#\T]^{1/s}.
    \end{equation}
   The parameter $\beta_{\phi}(s)$ and implied constant in \eqref{eq: delta multiplicity bound} are independent of $\delta >0$ and the choice of family $\T$. 
\end{theorem}

In view of our earlier observations, we may apply Theorem~\ref{thm: improv geometric bounds} to the family of tubes $\wt{\T}$ with $\delta=R^{1/2}/\lambda$. Using $\lambda\leq R^{1+\eps_0}$, this gives
\begin{align*}
 \big\|\sum_{T \in \T}\chi_{\tilde T} \big\|_{L^s(\R^n)} &\lesssim_{\phi,\eps}(\lambda/R^{1/2})^{-(n-1)/s+\beta_\phi(s)}[\#\T]^{1/s} \\
 & \leq R^{(n+1)/(2s)+\beta_\phi(s)/2+\eps_0(1/s + \beta_\phi(s))}\lambda^{-n/s}[\#\T]^{1/s}.
\end{align*}
To conclude the proof in this case, choose $\eps_0>0$ sufficiently small so that $\beta_\phi(s)+2\eps_0(1/s + \beta_\phi(s)) <\beta(n,s)$. Thus, we obtain the desired inequality \eqref{eq: rescaled Kakeya}.




\subsection{Small \texorpdfstring{$R$}{} regime: \texorpdfstring{$\lambda\geq R^{1+\varepsilon_0}$}{}}\label{subsec: small R} The scaling argument used above breaks down in the remaining regime, and so here we must argue differently. We first show that any rescaled phase $\phi^K$ satisfies a non-concentration condition in the spirit of \cite{KR2018} and \cite{CGGHIW}, where the strength of the non-concentration depends on $K$. This is a universal property of translation-invariant phases, in the sense that it does not depend on Hypothesis I). Once we have the non-concentration property, we can use the theory from \cite{CGGHIW} to show the desired Kakeya estimates hold. Unfortunately, one complication is that it is necessary to show that the arguments used in \cite{CGGHIW} are independent of the parameter $K$.




\subsubsection*{Properties of the phase and uniform derivative bounds} We begin with some technical observations regarding the uniformity of various derivatives of the phase $\phi^K$ and associated functions with respect to the $K$ parameter. For any H\"ormander-type phase and $K \geq 1$, it is clear that $\partial^2_{xy}\phi^K(x,t;y)=\partial^2_{xy}\phi(x/K,t/K;y)$ and 
\begin{equation*}
    \|\partial_{(x,t)}^{\alpha} \phi^K\|_{L^{\infty}(\D^n)} \leq \|\partial_{(x,t)}^{\alpha} \phi\|_{L^{\infty}(\D^n)} \qquad \textrm{for all $\alpha \in \N_0^{n-1} \times \N_0$ with $\alpha \neq 0$.} 
\end{equation*}
Recall, however, we are also assuming $\phi$ is translation-invariant and therefore
\begin{equation*}
    \phi(x,t;y)=\langle x,y\rangle+\psi(t;y)\qquad \textrm{where} \qquad \psi(0;y)=0.
\end{equation*}
This implies that $\partial_y^{\beta}\phi(x,0;y) = 0$ for all $(x,y) \in \B^{n-1} \times \B^{n-1}$ whenever $\beta \in \N_0^{n-1}$ with $|\beta| \geq 2$. It then follows from the mean value theorem that  
\begin{equation}\label{eq: unif derivs}
\|\partial^{\beta}_y\phi^K\|_{L^{\infty}(\D^n)} \leq \|\partial_t\partial^{\beta}_y\phi\|_{L^{\infty}(\D^n)} \qquad \textrm{for all $\beta \in \N_0^{n-1}$ with $|\beta| \geq 2$.}
\end{equation}

Given $K \geq 1$, we define 
\begin{equation}\label{eq: core curve}
   \Psi^K(\omega,t;y) \coloneq \omega-K\partial_y\psi(t/K;y) \quad \textrm{so that} \quad \partial_y\phi^K\big(\Psi^K(\omega,t;y), t; y \big) = \omega. 
\end{equation}
It is a straightforward matter, arguing as above, to show that the derivatives of $\Psi^K$ are also bounded, independent of $K$. 




\subsubsection*{Non-concentration} We recall the following terminology from \cite{CGGHIW}.

\begin{definition}
    Given a finite family $P_1,\ldots,P_m\in \R[X_1,\ldots,X_n]$ of non-zero polynomials, we let $Z(P_1,\ldots,P_m)$ denote their common zero set
    \begin{equation}\label{eq: dfn zero set}
        Z(P_1,\ldots,P_m)
        \coloneqq 
        \bigcap_{j=1}^m\{x\in\R^n: P_j(x)=0\}.
    \end{equation}
    We say that a set of the form \eqref{eq: dfn zero set} is a\textit{ transverse complete intersection} if the vectors $\nabla P_1(x),\ldots,\nabla P_m(x)$ are linearly independent for all $x\in Z(P_1,\ldots,P_m)$. In this case, $Z(P_1,\ldots,P_m)$ is a smooth submanifold of $\R^n$ of codimension $m$.
\end{definition}

\begin{definition}
    Given $0<\delta\leq \rho<1$, $0\leq m\leq n$ integers and $d\in\N$, we say a set $G\subseteq \R^n$ is a $(\delta,\rho,d)$-grain of codimension $m$ if 
    \[
    G=N_{\delta}Z(P_1,\ldots,P_m)\cap B(x_G,\rho)
    \]
    for non-zero polynomials $P_1,\ldots,P_m\in\R[X_1,\ldots,X_n]$ of degree at most $d$ such that $Z(P_1,\ldots,P_m)$ is a codimension $m$ transverse complete intersection and $x_G\in\R^n$. Here, $N_{\delta}Z(P_1,\ldots,P_m)$ denotes the (Euclidean) $\delta$-neighbourhood of $Z(P_1,\ldots,P_m)$. 
\end{definition}

\begin{definition}
Let $0<\delta<1$, $0\le m\le n$, $K \geq 1$ and $\kappa>0$,
$B_1,B_2\geq 0$.
We say a family $\T$ of $1 \times \delta$ tubes associated to $\phi^K$ is
$(\kappa,B_1,B_2)$-non-concentrated in codimension $m$ grains if the
following holds.

We have $\#\T\le\delta^{-(n-1)}$ and, for all
$d\in\mathbb{N}$ and $\varepsilon>0$, the inequality
\begin{equation}\label{eq:cgg-6.21}
    \#\bigl\{
        T\in\T:|T\cap G|\ge\lambda|T|
    \bigr\}
    \lesssim_{\phi,\varepsilon,d}
    (\lambda/\rho)^{-B_1}\rho^{-B_2}
    \delta^{-(n-1)+\kappa-\varepsilon}
\end{equation}
holds whenever $G\subseteq\R^n$ is a
$(\delta,\rho,d)$-grain of codimension $m$ for
$\delta\le\lambda\le\rho\le1$, where the implied constant in \eqref{eq:cgg-6.21}
depends on $\phi$, $n$, $\varepsilon$ and $d$ only; in particular, it is independent of $K$.
\end{definition}

The key non-concentration result is as follows.

\begin{theorem}[Non-concentration for grains]\label{thm: non-concentration} 
Let $n\geq 3$, $0 < \kappa< 1$ and $K \geq \delta^{-\kappa}$. Suppose that $\T$ is a direction-separated  family of $1 \times \delta$ tubes associated to $\phi^K$. Then $\T$ is $(\tfrac{\kappa}{n-2},n,n)$-non-concentrated in codimension 1 grains. 
\end{theorem}

The proof of Theorem~\ref{thm: non-concentration} involves a variant of the argument used in \cite{KR2018} and \cite[Theorem 4.6]{CGGHIW}. We postpone the details until \S\ref{subsec: nonconc proof} below.




\subsubsection*{Multilinear theory} We can use the theory developed  \cite{CGGHIW} to derive multilinear Kakeya estimates from the non-concentration inequalities from Theorem~\ref{thm: non-concentration}. Given a $1 \times \delta$ tube $T = T_{y, \omega}$ associated to $\phi^K$, recalling \eqref{eq: core curve}, we let 
\begin{equation*}
    \Gamma_T^K(t) \coloneq \big(\Psi^K(\omega, t; y), t \big) \qquad \textrm{for $t \in (-1,1)$}
\end{equation*}
denote the (parametrisation) of the \textit{core curve} of $T_{y, \omega}$.

\begin{definition}
Let $0<\delta<\nu<1$.
\begin{enumerate}
\renewcommand{\labelenumi}{\roman{enumi})}

\item
Given $2\le k\le n$, we say
$(T_1,\ldots,T_k)\in\T^k$ is $\nu$-transversal if\footnote{
We remark that \eqref{eq:cgg-6.19} agrees with \cite[(6.3)]{CGGHIW} for an
appropriate choice of $\nu$ when the tubes are defined by a
phase function.
}
\begin{equation}\label{eq:cgg-6.19}
    \inf_{\substack{
        x_\ell\in\Gamma_{T_\ell}^K,\ |x-x_\ell|<\delta\\
        1\le\ell\le k
    }}
    \left|
        \bigwedge_{\ell=1}^k T_{x_\ell}\Gamma_{T_\ell}^K
    \right|
    \ge\nu>0
    \qquad
    \text{for all }x\in T_1\cap\cdots\cap T_k.
\end{equation}
Here $T_{x_\ell}\Gamma_{T_\ell}^K$ denotes the tangent space to
$\Gamma_{T_\ell}^K$ at $x_\ell$ and the infimum is taken over all
tuples
\begin{equation*}
     (x_1,\ldots,x_k)
    \in\Gamma_{T_1}^K\times\cdots\times\Gamma_{T_k}^K
\qquad \textrm{with $|x-x_\ell|<\delta$ for $1\le\ell\le k$.}   
\end{equation*}

\item
We say $\T_1,\ldots,\T_k\subseteq\T$
are $\nu$-transversal if every $k$-tuple $(T_1,\ldots,T_k)$ with
$T_\ell\in\T_\ell$, $1\le\ell\le k$, is
$\nu$-transversal.
\end{enumerate}
\end{definition}

The key multilinear inequality reads as follows.

\begin{theorem}[{Multilinear Kakeya estimate; \cite[Theorem 6.13]{CGGHIW}}]\label{thm: MK}
For $2\le k\le n$, $1<s<\frac{k}{k-1}$, $\kappa>0$ and  $K \geq \delta^{-\kappa}$  there exists a constant $A=A_{n,k,s}\ge1$, and some $\eta(s) = \eta_{n,k, \kappa}(s)>0$ such that
the following holds. Suppose that $\T$ is a direction-separated  family of $1 \times \delta$ tubes associated to $\phi^K$. Then the inequality
\begin{equation}
    \Big\|
        \prod_{\ell=1}^k
        \big|
            \sum_{T_\ell\in\T_\ell}\chi_{T_\ell}
        \big|^{1/k}
    \Big\|_{L^s(\R^n)}
    \lesssim_{\phi,s}
    \nu^{-A}\delta^{-(n-k)/s'+\eta(s)}
    \big(
        \sum_{T\in\T}|T|
    \big)^{1/s}
    \label{eq:cgg-6.22}
\end{equation}
holds whenever
$\T_1,\ldots,\T_k\subseteq\T$
are $\nu$-transversal families for some $\nu>0$.
\end{theorem}

Importantly, the implied constant in \eqref{eq:cgg-6.22} and the exponents $A$ and $\eta(s)$ are independent of $K$.

\begin{proof}[Proof (of Theorem~\ref{thm: MK})] The result follows directly from \cite[Theorem 6.13]{CGGHIW} and Theorem~\ref{thm: non-concentration}. In particular, Theorem~\ref{thm: non-concentration} ensures the family of tubes $\T$ satisfies the required non-concentration hypothesis of \cite[Theorem 6.13]{CGGHIW}. The uniformity in $K$ arises from the uniform bounds for the derivatives of $\Psi^K$: see the discussion following the statement of \cite[Theorem 6.13]{CGGHIW}.
\end{proof}





\subsubsection*{From multilinear to linear Kakeya estimates} Let $\T$ be a family of $1 \times R^{-1/2}$ tubes as in \eqref{eq: rescaled Kakeya}. Recall we are assuming $R^{1 + \varepsilon_0} \leq \lambda$; we let $K \coloneq \lambda/R$ and $\delta \coloneq R^{-1/2}$ and apply the above results with these parameter values. To conclude the proof of Corollary~\ref{cor: beyond univ}, we show that the resulting multilinear Kakeya estimates from Theorem~\ref{thm: MK} imply linear counterparts. To do this, we use a primordial version of the Bourgain--Guth broad/narrow argument \cite{bourgain--guth}. This is carried out in \cite{CGGHIW}; however, it remains to verify that the estimates involved are independent of the scaling parameter $K$.

Henceforth, fix $K \coloneq \lambda/R$. Let $L \geq 1$ and $\Theta(L^{-1})$ be as defined in \S\ref{sec: notation}. Given a set of $1 \times \delta$ tubes associated to $\phi^K$, we define $\T[\theta]=\{T_{y, \omega} \in\T:y \in \theta \}$ for $\theta \in \Theta(L^{-1})$. We further let $\Q_{L^{-1}}$  be a tiling of $\B^{n-1}\times \B^1$ by $L^{-1}$-cubes. 

\begin{lemma}[Broad/narrow decomposition] Let $K$, $L \geq 1$, $1 \leq s < \infty$ and suppose $\T$ is a family of $1 \times \delta$ tubes associated to $\phi^K$. For each $Q\in\Q_{L^{-1}}$, there exists $\T_1[Q],\dots,\T_k[Q]\subseteq\T$ which are $L^{-(k-1)}$-transversal and such that 
\begin{multline}\label{eq: broad/narrow}
    \Big\|\sum_{T\in\T}\chi_T\Big\|_{L^s(\mathbb R^n)}
    \lesssim_\phi
    L^{(k-2)/s'}
    \Big(
        \sum_{\theta\in\Theta(L^{-1})}
        \Big\|
            \sum_{T\in\T[\theta]}\chi_T
        \Big\|_{L^s(\mathbb R^n)}^s
    \Big)^{1/s}
    \\
    +
    L^{(n-1)(1+k/s)}
    \Big(
        \sum_{Q\in\mathcal Q_{L^{-1}}}
        \Big\|
            \prod_{j=1}^k
            \big|
                \sum_{T\in\T_j[Q]}\chi_T
            \big|^{1/k}
        \Big\|_{L^p(Q)}^s
    \Big)^{1/s}.
\end{multline}
\end{lemma}

\begin{proof} The proof is exactly the same as in \cite[Corollary 6.8]{CGGHIW}, except we must justify why the implicit constant is uniform in $K$. First, the caps $\Theta(L^{-1};V;\bdx)$ in \cite[(6.5)]{CGGHIW} are defined using a constant $C_\phi\geq 1$. This constant appears in the estimate \begin{equation*}
    \#\Theta(L^{-1};V;\bdx)\lesssim C_\phi^{n-k+1}L^{n-2}, 
\end{equation*}
in the proof of their Corollary 6.8. This constant is later chosen sufficiently large, so that $$|\bdz-\bdx|_\infty<2L^{-1}\implies |G_\phi(\bdz;y)-G_\phi(\bdx;y)|\lesssim C_\phi L^{-1},$$ where $G_\phi$ is the Gauss map for their phase $\phi$. For our rescaled phases $\phi^K$, we have 
\begin{equation*}
|G_{\phi^K}(\bdz;y)-G_{\phi^K}(\bdx;y)|\lesssim\|\partial_{\bdx}G_{\phi^K}\||\bdz-\bdx| \lesssim K^{-1}\|\partial_{\bdx} G_{\phi}\| |\bdz-\bdx|,
\end{equation*}
where we have used $G_{\phi^K}(\bdx;y)=G_{\phi}(\bdx/K;y)$. Thus, we may also choose the constant $C_\phi$ large enough independently of $K$, such that $C_\phi\gtrsim \|\partial_{\bdx}G_{\phi}\|$.

In addition, given $|y_1 - y_2| \lesssim L^{-1}$ and $|\Psi^K(\omega_1, t; y_1) - \Psi^K(\omega_2, t; y_2)| \lesssim L^{-1}$, we need also to note that
\begin{equation*}
   |\omega_1 - \omega_2| = \big|(\partial_y \phi^K)(\Psi^K(\omega_1, t; y_1),t;y_1) - (\partial_y \phi^K)(\Psi^K(\omega_2, t; y_2),t;y_2)\big| \lesssim_\phi L^{-1},
\end{equation*}
where the implied constant is independent of $K$. This is a consequence of the uniform bounds on the derivatives of the phase, as discussed above (in particular, we use \eqref{eq: unif derivs} which relies on the translation-invariance of $\phi$). 
\end{proof}

It remains to bound the terms appearing on the right-hand side of \eqref{eq: broad/narrow}. For the first term, we use the following variant of the universal Kakeya inequality \cite[Theorem 21]{WisewellThesis}.

\begin{lemma}
    \label{lemma: universal kakeya}
    Let $n\geq 3$, $K \geq 1$ and $\T$ be a family of distinct $1\times\delta$-tubes associated to $\phi^K$, where $\phi$ is translation-invariant. Then for $1\leq s\leq\tfrac{n+1}{n-1}$ and all $\eps>0$, we have
\begin{equation}
    \label{eq: universal kakeya}
    \big\|\sum_{T\in\T}\chi_T\big\|_{L^s(\R^n)}\lesssim_{\phi,\eps}\delta^{(n-1)/s-\eps}[\#\T]^{(s+1)/(2s)}.
\end{equation}
\end{lemma}

Importantly, the implied constant in \eqref{eq: universal kakeya} is independent of $K$. We discuss the proof of Lemma~\ref{lemma: universal kakeya}, and in particular this uniformity, in Appendix~\ref{sec: tube intersect} below. We remark that the above lemma holds for all families $\T$ of distinct $1 \times \delta$ tubes $T = T_{y, \omega}$ (that is, direction-separation is not required), with  $y \in \delta \Z^{n-1}\cap\B^{n-1}$ and $\omega\in \delta \Z^{n-1}$. We also remark that here the phase $\phi$ is again assumed to be translation-invariant. 

Direction-separation implies $\#\T[\theta]\lesssim (L\delta)^{-(n-1)}$ for all $\theta \in \Theta(L^{-1})$. We apply Lemma~\ref{lemma: universal kakeya} to each family $\T[\theta]$ and take $\ell^s(\Theta(L^{-1}))$-sum to deduce that
\begin{align}
\nonumber
\Big(\sum_{\theta\in\Theta(L^{-1})}\big\|\sum_{T\in\T[\theta]}\chi_T\big\|^s_{L^s(\R^n)}\Big)^{1/s}&\lesssim_{\phi,\eps} \delta^{(n-1)/s-\eps}\Big(\sum_{\theta\in\Theta(L^{-1})}\big[\#\T[\theta]\big]^{(s+1)/2}\Big)^{1/s}\\
\label{eq: B/N 1}
    &\lesssim L^{-(n-1)/2s'}\delta^{(n-1)/s-\beta(n,s)-\eps}[\#\T]^{1/s}.
\end{align}

By the hypothesis of the small $R$ regime, $K \coloneq \lambda/R \geq R^{\varepsilon_0} = \delta^{-2\eps_0}$. We may apply Theorem \ref{thm: MK} to the multilinear term in \eqref{eq: broad/narrow} with $\kappa=2\eps_0$ and $k = (n+1)/2$ (using the hypothesis that $n$ is odd). From this, we deduce that
\begin{equation}\label{eq: B/N 2}
    \Big\|
            \prod_{j=1}^k
            \big|
            \sum_{T\in\T_j[Q]}\chi_T
            \big|^{1/k}\Big\|_{L^s(Q)}\lesssim_\phi L^{A(k-1)}\delta^{(n-1)/s-\beta(n,s)+\eta(s)}[\#\T]^{1/s}.
\end{equation}

Combining \eqref{eq: B/N 1} and \eqref{eq: B/N 2}, we deduce that 
\begin{equation*}
\big\|\sum_{T\in\T}\chi_T\big\|_{L^s(\R^n)}\lesssim_{\phi,\eps} \delta^{(n-1)/s-\beta(n,s)-\eps}\big(L^{-(n+3-2k)/2s'}+L^{(n-1)(A+1+k/s)}\delta^{\eta(s)}\big)[\#\T]^{1/s}.
\end{equation*}
Since $k=(n+1)/2$, it follows that $n+3-2k>0$. The above inequality holds for all $L \geq 1$. We choose $L$ to satisfy 
\begin{equation*}
    L^{(n-1)(A+1+k/s)+(n+3-2k)/2s'}=\delta^{-\eta(s)}.
\end{equation*}
Recalling $\delta = R^{-1/2}$, this gives an estimate of the desired form \eqref{eq: rescaled Kakeya}.




\subsection{Universal non-concentration for scaled phases}\label{subsec: nonconc proof} Here we describe how to adapt the strategy of \cite[Theorem 4.6]{CGGHIW} (which is a variant of the Katz--Rogers argument \cite{KR2018}) to prove Theorem~\ref{thm: non-concentration}.

Fix a $(\delta,\rho,d)$-grain $G$ of codimension $m \coloneq 1$ and $\varepsilon > 0$ and let $\T$ be a collection of $1 \times \delta$ tubes associated to a translation-invariant, Hörmander-type phase $\phi^K$. The uniform smoothness bounds for the $\Psi^K$ ensure that the polynomial
approximation in Step 1 of the proof of
\cite[Theorem 4.6]{CGGHIW} applies to $\Psi^K$ with degree of Taylor approximation depending only on $n$ and $\varepsilon$, and with a
smallness threshold $\delta_0=\delta_0(\phi,n,\varepsilon)>0$
independent of $K$. The resulting polynomial parametrisations
therefore have degree bounded independently of $K$.
Consequently, the semialgebraic decomposition in Step 2 also
has complexity bounds independent of $K$. Following these steps, we obtain
an interval $\bar{I}\subset \R$, a finitely overlapping collection $\mathcal{J}(\lambda)$ of open intervals of length $(c_{d,\varepsilon}/10)\lambda$ that partitions $\bar{I}$, and a set of (portions of) $\delta$-tubes 
\[
\T[J]
\coloneqq 
\{T\in\T: \wt{T}[\delta,J]\subset N_{2\delta}G\}
\]
that approximate\footnote{More explicitly, for $T\in\T$ and $J\in \mathcal{J}(\lambda)$, $\widetilde{T}[\delta,J]$ denotes 
\[
\wt{T}[\delta,J]
\coloneqq 
\{(x,t)\in\B^{n-1}\times J: |x-\wt{\gamma}_T(t)|<\delta\},
\]
with $\widetilde{\gamma}_T$ a high-order Taylor approximation of the core curve of $T$. } the original tubes $T\in \T$, such that $(\kappa_{\circ},n,n)$-non-concentration of $\T$ reduces to the following estimate (see  \cite[(4.9)]{CGGHIW}):
\begin{equation*}
    \sum_{J\in \mathcal{J}(\lambda)}\#\T[J]\lesssim_{\phi,d,\varepsilon}\lambda^{-n}\delta^{-(n-1)+\kappa_{\circ}-\varepsilon/2}.
\end{equation*}
This is implied by
\begin{equation}\label{eq: non-conc 0}
    \lambda^{-n}|N_{3\delta}G \cap (\R^{n-1} \times J)| + \delta^{\kappa_{\circ}} \gtrsim_{\phi, d, \varepsilon} \delta^{n-1 + \varepsilon/2} \#\T[J].
\end{equation}
Indeed, using that the intervals $J$ have finite overlap, that $\#\mathcal{J}(\lambda)\lesssim_{d,\varepsilon}\lambda^{-1}$ and that $|N_{3\delta}G|\lesssim \delta^m=\delta^1$ by Wongkew's theorem \cite{Wongkew1993}, we get that \eqref{eq: non-conc 0} implies
\begin{align*}
    \sum_{J\in\mathcal{J}(\lambda)}\#\T[J]
    &\lesssim_{\phi,d,\varepsilon}
    (|N_{3\delta}G|\lambda^{-n}
    +
    \lambda^{-1}\delta^{\kappa_{\circ}})
    \delta^{-(n-1)-\varepsilon/2}
    \lesssim 
    \lambda^{-n}\delta^{\kappa_{\circ}-(n-1)-\varepsilon/2}.
\end{align*}

The goal is now to prove \eqref{eq: non-conc 0} for an arbitrary, fixed $J\in\mathcal{J}(\lambda)$. The uniform first derivative bounds for the polynomial
approximation of $\Psi^K$ ensure that the
argument in Step 3 of \cite[Theorem 4.6]{CGGHIW} has constants
independent of $K$. Moreover, the semialgebraic sets involved
have complexity bounded independently of $K$. Applying the
semialgebraic selection, Gromov parametrisation and polynomial
approximation in Steps 3--4, we therefore obtain an open set $\wt{B} \subset [0,1]^{n-1}$ and polynomial mappings\footnote{The polynomials $(\wt{P},\wt{D})$ are Taylor approximations of maps $(P,D)$ that parametrise
\[
\mathcal{L}(J)
\coloneqq 
\{
(\omega,y)\in\Omega_{\circ}\times Y_{\circ}:
\wt{T}_{y,\omega}[\delta/2,J]\subseteq N_{2\delta}G
\}.
\]} $\wt{D},\wt{P} \colon \R^{n-1} \to \R^{n-1}$ such that 
\begin{equation}\label{eq: non-conc 1}
    \|\partial_x\wt{D}(x)\| + \|\partial_x \wt{P}(x)\| \lesssim 1 \qquad \textrm{for all $x \in \wt{B}$}.
\end{equation}
For $W:=\wt{D}(\wt{B})$,
note that $\wt{D}|_{\wt{B}}\colon \wt{B}\to W$ is a diffeomorphism.

The above objects satisfy 
    \begin{equation}\label{eq: non-conc 1.5}
        |W| 
        \gtrsim_{\phi, d, \varepsilon}
        \delta^{n-1 +\varepsilon/2} \# \T[J].
    \end{equation}
The degrees of $\widetilde P,\widetilde D$ and the constants
in \eqref{eq: non-conc 1} and \eqref{eq: non-conc 1.5} are bounded independently of $K$, although the
maps themselves may depend on $K$.
The right-hand side is exactly the expression appearing in \eqref{eq: non-conc 0}. Hence, we only need to prove 
    \begin{equation}\label{eq: non-conc 1.5.1}
        |W| 
        \lesssim_{\phi, d, \varepsilon}
        \lambda^{-n}|N_{3\delta}G\cap (\R^{n-1}\times J)|+\delta^{\kappa_{\circ}}.
\end{equation}
We have
\[
|W|
=
|\wt{D}(\wt{B})|
=
\int_{\wt{D}(\wt{B})}1\, \ud y
=
\int_{\wt{B}}|\det \partial_{x}\wt{D}(x)|\, \ud x.
\]
We split the domain of integration into two parts. Let $C_{\circ}\coloneqq C_{\circ}(n,\phi)\geq 1$  be a constant that we will determine later. Define
$$
U(\delta^{\kappa_{\circ}})
\coloneqq 
\{x\in\wt{B}: |\det \partial_x \wt{D}(x)|\geq C_{\circ}\delta^{\kappa_{\circ}}\}.
$$
The integral over $\wt{B}\setminus{U(\delta^{\kappa_{\circ}})}$ is easily bounded:
\[
\int_{\wt{B}\setminus{U(\delta^{\kappa_{\circ}})}}|\det \partial_x \wt{D}(x)|
\leq 
|\wt{B}|C_{\circ}\delta^{\kappa_{\circ}}
\leq 
C_{\circ}\delta^{\kappa_{\circ}},
\]
having used that $\wt{B}\subset [0,1]^{n-1}$. So to obtain \eqref{eq: non-conc 1.5.1} it is enough to prove 
\begin{equation}\label{eq: non-conc 1.5.2}
\int_{U(\delta^{\kappa_{\circ}})}|\det \partial_x \wt{D}(x)|\, \ud x
\lesssim 
\lambda^{-n}|N_{3\delta}G\cap (\R^{n-1}\times J)|.    
\end{equation}
Notice that, if $C_{\circ}\delta^{\kappa_{\circ}}> \sup_{x\in\wt{B}}|\det \partial_x \wt{D}(x)|$, then $U(\delta^{\kappa_{\circ}})=\varnothing$, and therefore \eqref{eq: non-conc 1.5.2} holds trivially. Hence, we can assume that $C_{\circ}\delta^{\kappa_{\circ}}\leq \sup_{x\in\wt{B}}|\det \partial_{x}\wt{D}(x)|$. Notice that, in this case, we can use \eqref{eq: non-conc 1} to bound
\begin{equation}\label{eq: non-conc 1.5.3}
    C_{\circ}\delta^{\kappa_{\circ}}
    \leq 
    \sup_{x\in\wt{B}}|\det \partial_x \wt{D}(x)|
    \leq 
    \sup_{x\in\wt{B}}\|\partial_{x}\wt{D}(x)\|^{n-1}
    \lesssim 
    1.
\end{equation}

To get \eqref{eq: non-conc 1.5.2}, we again follow \cite{CGGHIW} and
define\footnote{Since $(P,D)$ parametrises the set of pairs of centres-directions whose corresponding tubes defined over $J$ are in $N_{2\delta}G$, the map $\fkp$ can be thought of as taking directions-to-centres.}
$\fkp(y) \coloneq \widetilde{P}\circ\widetilde{D}^{-1}(y)$ and $\fkP(y) \coloneq \partial_y \fkp(y)$ and
    \begin{align*}
A(t;y) &\coloneq \fkP(y) - (\partial_{yy}^2\phi^K)\big(\Psi^K(\fkp(y), t;y), t; y\big) \\
&\coloneq \fkP(y) - K(\partial_{yy}^2\phi)\big(\Psi(\fkp(y)/K, t/K;y), t/K; y\big)
\qquad \textrm{ for all $y\in W$.}
    \end{align*}
From \cite[(4.24)]{CGGHIW},
\begin{align}\label{eq: non-conc 2}
    |N_{3\delta}G \cap (\R^{n-1} \times J)| + \delta^{10n}
    &\gtrsim
    \int_J \int_{W} |\det A(t;y)|\,\ud y \ud t 
    \\&\geq 
    \int_{U(\delta_{\kappa_{\circ}})}\left(\int_J|\det A(t;\wt{D}(x))|\, \ud t\right)|\det\partial_x \wt{D}(x)|\, \ud x.\nonumber
\end{align}
Indeed, translation-invariance gives $$\Psi^K(\fkp(y),t;y)=\fkp(y)-K\partial_y\psi(t/K;y),$$ and so the Jacobian determinant of the mapping $(y,t)\mapsto (\Psi^K(\fkp(y),t;y),t)$ is precisely $A(t;y)$.
Notice that \eqref{eq: non-conc 2} is a slightly simpler and more explicit version of \cite[(4.24)]{CGGHIW}. We may work with the simplified \eqref{eq: non-conc 2} here since we are only considering non-concentration in codimension $1$ grains, which corresponds to setting $m=1$ in the argument in \cite{CGGHIW}. This significantly simplifies the analysis in Step 6 of the proof of \cite[Theorem 4.6]{CGGHIW} and avoids the appearance of certain error terms. We further remark that the (harmless) $\delta^{10n}$ term in \eqref{eq: non-conc 2} arises from polynomial approximating, as in \cite[(4.21)]{CGGHIW}.

Hence, to obtain \eqref{eq: non-conc 1.5.2} it is enough to check that 
\[
\int_{J}|\det A(t,\wt{D}(x))|\, dt
\gtrsim 
\lambda^n,
\]
which is implied by 
\begin{equation*}
    \inf_{x\in U(\delta^{\kappa_{\circ}})}|\{t\in J: |\det A(t,\wt{D}(x))|> c\lambda^{n-1}\}\gtrsim \lambda,
\end{equation*}
which in turn is implied by  
\begin{equation}\label{eq: sublevel set estimate}
\sup_{x\in U(\delta^{\kappa_{\circ}})}|\{t\in J:|\det A(t,\wt{D}(x))|\leq \mu\}|
\lesssim 
\mu^{\frac{1}{n-1}} \quad \textrm{for all $0 < \mu < 1$.}
\end{equation}
To prove this sublevel set estimate, we use van der Corput's lemma (see \cite[Lemma 3.3]{Christ1985}): if we denote
\[
u_x(t)\coloneqq \det A(t;\wt{D}(x)),
\]
then \eqref{eq: sublevel set estimate} is implied by 
\begin{equation}\label{eq: van der Corput condition}
\inf_{x\in U(\delta^{\kappa_{\circ}})}\inf_{t\in J}|(u_x)^{(n-1)}(t)|\geq 1.    
\end{equation}

Let $\mathcal{H}$ denote the Hessian 
\begin{equation*}
\mathcal{H}(\omega,s;y)\coloneq\partial^2_{yy}\phi(\Psi(\omega,s;y),s;y).  
\end{equation*}
Recall, by translation-invariance, we have $\phi(0,0; y) = 0$ for all $y \in \B^{n-1}$, and so $\mathcal{H}(0,0;y)=0$ for all $y \in \B^{n-1}$. For fixed $y \in \B^{n-1}$, we Taylor expand the function $(\omega,t)\mapsto \mathcal{H}(\omega,t;y)$ around $(0,0)$, obtaining 
\begin{align*}
    \mathcal{H}(\omega,t,y)
    = 
    L(y)\cdot \omega 
    +
    J(y)\cdot t
    +
    R(\omega,t,y),
\end{align*}
with $L(y)\coloneqq\partial_\omega \mathcal{H}(0,0,y)$ and $J(y)\coloneqq \partial_t \mathcal{H}(0,0,y)$. Notice that the remainder satisfies
\begin{equation}
    \label{eq: remainder bound}
    \|R(\omega,t;y)\|\lesssim |\omega|^2+t^2.
\end{equation}
If $\mathcal{H}^K(\omega,t,y)\coloneqq K\cdot \mathcal{H}(\frac{\omega}{K},\frac{t}{K},y)$ and $R^K(\omega,t,y)\coloneqq K\cdot R(\frac{\omega}{K},\frac{t}{K},y)$, then
\[
    \mathcal{H}^K(\omega,t,y)
    =
    L(y)\omega+J(y)t+R^K(\omega,t,y).
\]

    We have that
    \begin{equation}\label{eq: Bound Derivatives Remainder}
    \|\partial_t^{(j)}R^K(\omega,t;y)\|\lesssim K^{-1}\qquad \textrm{for all $0\leq j\leq n-1$.}
    \end{equation}
    Indeed, for $j=0$ this follows from \eqref{eq: remainder bound}. For $j=1$, it follows from the definition of $R$ and the mean value theorem. And for $j\geq 2$, we have from the definition of $R$ that
    \[
    \partial_t^{(j)}R^K(\omega,t,y)
    =
    K^{1-j}\partial_s^{(j)}\mathcal{H}\Big(\frac{\omega}{K},\frac{t}{K},y\Big),
    \]
    so
    \[
    \|\partial_t^{(j)}R^K(\omega,t,y)\|
    \lesssim 
    K^{1-j}\|\partial_s^{(j)}\mathcal{H}\|
    \lesssim_{\phi}
    K^{1-j}.
    \]
    
We claim that $J(y)$ is invertible and satisfies 
\begin{equation}\label{eq: properties J}
|\det J(y)|\gtrsim 1,\qquad \|J(y)^{-1}\|\lesssim 1,    
\end{equation}
uniformly in $y \in \B^{n-1}$. To see this, note, by the chain rule,
    \[
    J(y)
    =
    \partial_s \mathcal{H}(0,0,y)
    =
    \partial_{\bdx}\partial_{yy}^2\phi(0,0;y)\cdot (\partial_s\Psi(0,0;y),1).
    \]
    Since $\partial_y\phi(\Psi(0,s;y),s,y)=0$, differentiating with respect to $s$ and evaluating at $s=0$ we have
    \begin{align*}
        0
        &=
        [\partial_{\omega}\partial_y \phi(\Psi(0,s,y),s,y)\cdot \partial_s\Psi(0,s,y)+\partial_s\partial_y\phi(\Psi(0,s,y),s,y)]_{s=0}
        \\&=
        \partial_y\partial_{\bdx}\phi(0,0,y)\cdot (\partial_s\Psi(0,0,y),1).
    \end{align*}
    That is, the vector $(\partial_s\Psi(0,0,y),1)$  is orthogonal to every vector $\partial_{y_j}\partial_{\bdx}\phi(0,0,y)$, $j=1,\ldots,n-1$. Hence, it must be proportional to the unit normal vector at that point, which is $G(0,0,y)$. That is,
    \[
    (\partial_s\Psi(0,0,y);1)
    =
    a(y)G(0,0,y),
    \]
    for some scalar $a(y)$ that satisfies $|a(y)|=|(\partial_s\Psi(0,0,y),1)|\simeq 1$. Therefore,
    \begin{align*}
    J(y_0)
    &=
    \partial_{\bdx}\partial_{yy}^2\phi(0,0,y_0)\cdot a(y_0)G(0,0,y_0)
    \\&=
    a(y_0)\left.\partial_{yy}^2\langle \partial_{\bdx}\phi(0,0,y),G(0,0,y_0)\rangle\right|_{y=y_0}.
    \end{align*}
    Since $\phi$ is of Hörmander-type, by H2 we have $|\det(J(y_0))|\gtrsim 1$ uniformly in $y_0$. This gives us the first claim in \eqref{eq: properties J}.

    For the second claim in \eqref{eq: properties J}, we just apply Cramer's rule $A^{-1}=\frac{\operatorname{adj}(A)}{\det(A)}$ to obtain that
    \[
    \|J(y)^{-1}\|
    =
    \frac{\|\operatorname{adj} J(y)\|}{|\det J(y)|}
    \lesssim 
    1,
    \]
    having used the bound $|\det(J)|\gtrsim 1$ proved in the first item and the fact that $\|\operatorname{adj}J(y)\|\lesssim 1$. Indeed, this is true because the entries of the adjugate $\operatorname{adj}(J(y))$ are sums of products of $(n-2)$-tuples of elements (up to sign) of 
    \[
    J(y)=a(y_0)\left.\partial_{yy}^2\langle \partial_{\bdx}\phi(0,0,y),G(0,0,y_0)\rangle\right|_{y=y_0},
    \]
    and these elements are continuous functions on a compact set: hence we can bound them by an absolute constant that only depends on $\phi$.

The matrix $J(y)$ is defined using the original phase $\phi$
and is independent of $K$. Its bounds in \eqref{eq: properties J} hold uniformly
on the fixed parameter neighbourhood. They therefore also
hold for $J_x=J(\widetilde D(x))$, even though
$\widetilde D$ may depend on $K$.

With the two claims, we now can prove \eqref{eq: van der Corput condition}. Recalling the definition of $A(t;y)$, we have
\begin{align*}
    A(t;y) 
    &=
    \fkP(y)-K\cdot \mathcal{H}\Big(\frac{\fkp(y)}{K},\frac{t}{K},y\Big)
    =
    \fkP(y)-\mathcal{H}^K(\fkp(y),t,y),
\end{align*}
where recall that $\fkp(y)=\wt{P}\circ\wt{D}^{-1}(y)$. Hence, if we define 
\[
J_x\coloneqq J(\wt{D}(x))=\partial_s \mathcal{H}(0,0,\wt{D}(x)),
\]
we get
\begin{align*}
    A(t,\wt{D}(x))
    &=
    \fkP(\wt{D}(x))-L(\wt{D}(x))\wt{P}(x)-t\cdot J_x-R^K(\wt{P}(x),t,\wt{D}(x))
    \\&=
    J_x[J_x^{-1}\fkP(\wt{D}(x))-J_x^{-1}L(\wt{D}(x))\wt{P}(x)-J_x^{-1}R^K(\wt{P}(x),t,\wt{D}(x))-tI_{n-1}
    ]
    \\&\eqqcolon 
    J_x[M^K(x,t)-tI_{n-1}].
\end{align*}
We decompose the matrix $M^K(x,t)$ into a part that only depends on $x$ and a part that depends on both $x$ and $t$, that is,
\[
M^K(x,t)
=
M_0(x)-M_1^K(x,t)
\]
with
\begin{align*}
    &M_0(x)\coloneqq J_x^{-1}[\fkP(\wt{D}(x))-L(\wt{D}(x))\wt{P}(x)],
    \\&
    M_1^K(x,t)\coloneqq J_x^{-1}R^K(\wt{P}(x),t,\wt{D}(x)).
\end{align*}
To lower bound the determinant of $A$, we first prove bounds for the norms of $M_0$ and $M_1^K$. We start with $M_0$. Notice that $|\det \partial_x \wt{D}(x)|^{-1}\leq (C_{\circ}\delta^{\kappa_{\circ}})^{-1}$ on $U(\delta^{\kappa_{\circ}})$. Using this and \eqref{eq: non-conc 1}, we can bound
\[
\|\fkP(\wt{D}(x))\|
\lesssim (C_{\circ}\delta^{\kappa_{\circ}})^{-1}.
\]
Also notice that, by \eqref{eq: non-conc 1} and \eqref{eq: non-conc 1.5.3},
\[
\|L(\wt{D}(x))\wt{P}(x)\|
\lesssim_{\phi}
1
\lesssim 
C_{\circ}^{-1}\delta^{-\kappa_{\circ}}.
\]  
Using these bounds and the fact that $\|J_x^{-1}\|\lesssim 1$ uniformly in $x$ by \eqref{eq: properties J}, we obtain
\begin{equation}\label{eq: Bound M0}
    \sup_{x\in U(\delta^{\kappa_{\circ}}),t\in J}\|M_0(x)\|
    \lesssim_{\phi}
    C_{\circ}^{-1}\delta^{-\kappa_{\circ}},    
\end{equation}
Now we prove an analogous bound for $M_1^K(x,t)$. We have that
\[
\partial_t^{(j)}M_1^K(x,t)
=
J_x^{-1}\cdot \partial_t^{(j)}R^K(\wt{P}(x),t,\wt{D}(x)).
\]
So, by \eqref{eq: Bound Derivatives Remainder} and \eqref{eq: properties J},
\begin{equation}\label{eq: Bound M1}
    \|\partial_t^{(j)}M_1^K(x,t)\|
    \lesssim
    K^{-1}\qquad \textrm{for all $0 \leq j \leq n-1$.}
\end{equation}

Combining these bounds with the following lemma, we can finish the proof.

\begin{lemma}\label{lemma: Jonathan's argument}
    Let $M^K(x,t)\coloneqq M_0(x)- M_1^K(x,t)$ be the difference of $(n-1)\times (n-1)$-matrices satisfying the bounds \eqref{eq: Bound M0} and \eqref{eq: Bound M1}. Then,
    \begin{equation}\label{eq: lower bound vx derivatives}
     \inf_{x\in U(\delta^{\kappa_{\circ}})}\inf_{t\in J}|(v_x)^{(n-1)}(t)|\geq 1    \qquad \textrm{where} \qquad
    v_x(t)\coloneqq \det(M^K(x,t)-tI_{n-1}).   
    \end{equation}
\end{lemma}

Recall that
\[
u_x(t)
=
\det(A(t,\wt{D}(x))
=
\det(J_x)v_x(t).
\]
Since $|\det(J_x)|\gtrsim 1$ by \eqref{eq: properties J}, we have that
\[
|u_x^{(n-1)}(t)|
=
|\det(J_x)|\cdot |v_x^{(n-1)}(t)|
\gtrsim 
1.
\]
This concludes the argument, except that it remains to prove Lemma \ref{lemma: Jonathan's argument}.

\begin{proof}[Proof (of Lemma \ref{lemma: Jonathan's argument})]
    Observe that $v_x(t)$ is the characteristic polynomial of $M^K(x,t)$. Hence, 
    \begin{equation}\label{eq: Principal Minor Expansion}
    v_x(t)
    =
    (-1)^{n-1}t^{n-1}+\sum_{m=0}^{n-2}\alpha_m^K(x,t)t^m    
    \end{equation}
    being the coefficients $\alpha_m^K(x,t)$ (up to sign) the sum of the principal minors of size $(n-m-1)\times (n-m-1)$ of $M^K(x,t)$.

    Recall that $\kappa_\circ=\kappa/(n-2)$ and $K\geq\delta^{-\kappa}$. By the multilinearity of the determinant and the uniform bounds \eqref{eq: Bound M0} and \eqref{eq: Bound M1} for the matrix entries, one can check that
    \begin{equation*}
        \sup_{x \in U(\delta^{\kappa_{\circ}})} \sup_{t \in J}|\partial_t^{(j)} \alpha_m^K(x, t)| \lesssim K^{-1}C_\circ^{-(n-2)}\delta^{-(n-2)\kappa_\circ}\lesssim C_\circ^{-1} \qquad \textrm{for $1 \leq j \leq n-1$.}
    \end{equation*}
    Differentiating \eqref{eq: Principal Minor Expansion} we obtain 
    \begin{equation*}
        \pm (v_x)^{(n-1)}(t) = (n-1)! + O_{\phi, \varepsilon}(C_{\circ}^{-1}). 
    \end{equation*}
    Thus, provided $C_{\circ}$ is chosen sufficiently large, depending only on $\phi$ and $\varepsilon$, we deduce the desired bound \eqref{eq: lower bound vx derivatives}. This completes the proof. 
\end{proof}



\section{Sharpness of the universal off-diagonal estimates}
In this section, we prove the necessity of the two conditions in \eqref{eq: universal range off diagonal}. 

Given $\bdB \colon \R \to \operatorname{Sym}(n-1, \R)$ a function taking values in the space of real, symmetric $(n-1) \times (n-1)$ matrices and $\bdC \in \operatorname{Sym}(n-1, \R)$ a fixed real symmetric $(n-1) \times (n-1)$ matrix, consider
\begin{equation}\label{eq: structure of phase}
\phi(x,t;y)
\coloneq 
\langle x,y\rangle
+
\frac{1}{2}\langle \bdB(t)y,y\rangle
+
\frac{1}{2}\langle \bdC y,y\rangle.
\end{equation}
If $\det \bdB'(t) \neq 0$ for $|t| <1$, then $\phi$ is a H\"ormander-type phase and the signature of $\phi$ corresponds to the signature of $\bdB'(0)$. The sharp example is the tensor product of an elliptic term and a hyperbolic term.

The hyperbolic term is due to Minxing Shen, and was used in \cite[Theorem 3.3]{CGGHIW}. In dimension $n-\sigma$, define a H\"ormander-type phase $\phi_{\hyp}$ as in \eqref{eq: structure of phase} by taking $\bdB=\bdB_{\hyp}^{n-\sigma}$ and $\bdC=\bdC_{\hyp}^{n-\sigma}$ where
\begin{equation*}
\bdB_{\hyp}^{n-\sigma}(t) 
\coloneq 
\underbrace{\begin{pmatrix}
   -t^2  & t \\
    t & 0
\end{pmatrix}
\oplus 
\cdots 
\oplus
\begin{pmatrix}
   -t^2  & t \\
    t & 0
\end{pmatrix}}_{\frac{n-\sigma-1}{2}-\mathrm{fold}},   
\quad
\bdC_{\hyp}^{n-\sigma}
\coloneq 
\underbrace{
\begin{pmatrix}
   1  & 0 \\
   0 & 0
\end{pmatrix}
\oplus 
\cdots 
\oplus
\begin{pmatrix}
   1  & 0 \\
   0 & 0
\end{pmatrix}}_{\frac{n-\sigma-1}{2}-\mathrm{fold}}.    
\end{equation*}
It is then a straightforward matter to check $\phi_{\hyp}$ has signature $0$. In addition, let $a_{\hyp} \in C^{\infty}_c(\D^{n-\sigma})$ be nonnegative with $a_{\hyp}(0,0;0) = 1$ and $\supp a_{\hyp} \subseteq \bbD_{1/4}^{n-\sigma}$. 

For the elliptic term, we use constructions from \cite[\S2]{GHI}. In dimension $\sigma+1$, define a H\"ormander-type phase $\phi_{\elliptic}$ as in \eqref{eq: structure of phase} by taking $\bdB=\bdB_{\elliptic}^{\sigma+1}$ and $\bdC=\bdC_{\elliptic}^{\sigma+1}$ where
\begin{equation*}
\bdB_{\elliptic}^{\sigma+1}(t)
\coloneq
\underbrace{
\begin{pmatrix}
    t+t^3 & -t^2 \\
    -t^2 & t
\end{pmatrix}
\oplus 
\cdots 
\oplus
\begin{pmatrix}
    t+t^3 & -t^2 \\
    -t^2 & t
\end{pmatrix}
}_{\lfloor\frac{\sigma}{2}\rfloor-\mathrm{fold}}
\oplus 
\begin{pmatrix}
    t
\end{pmatrix}^{\ast}, \quad
\bdC_{\elliptic}^{\sigma+1}
\coloneq
\mathrm{Id}\oplus 
\begin{pmatrix} 
0
\end{pmatrix}^{\ast},   
\end{equation*}
here $\mathrm{Id}$ denotes the $2\lfloor\frac{\sigma}{2}\rfloor \times 2\lfloor\frac{\sigma}{2}\rfloor$ identity matrix and the asterisk $^{\ast}$ indicates that the final block only appears if $\sigma+1$ is even. In particular, when $\sigma+1=2$, we have $$\bdB_{\elliptic}^2(t)=(t),\quad \bdC_{\elliptic}^2=(0).$$ It is a straightforward matter to check $\phi_{\elliptic}$ has signature $\sigma$. In addition, let $a_{\elliptic} \in C^{\infty}_c(\D^{\sigma+1})$ be nonnegative with $a_{\elliptic}(0,0;0) = 1$ and $\supp a_{\elliptic} \subseteq \bbD_{1/4}^{\sigma+1}$.

We proceed as in \cite{HI2022} by taking the tensor product of the examples in the minimal and maximal signature cases. Define a H\"ormander-type phase $\phi_{\mix}$ as in \eqref{eq: structure of phase} by taking $\bdB = \bdB_{\mix}^{n,\sigma}$ and $\bdC = \bdC_{\mix}^{n,\sigma}$ where
\begin{equation*}
\bdB_{\mix}^{n,\sigma}
\coloneq
\bdB_{\hyp}^{n-\sigma}
\oplus 
\bdB_{\elliptic}^{\sigma+1} \qquad \textrm{and} \qquad
\bdC_{\mix}^{n,\sigma}
\coloneq
\bdC_{\hyp}^{n-\sigma}
\oplus 
\bdC_{\elliptic}^{\sigma+1}.
\end{equation*}
Then, for $x=(x^1,x^2)$, $y=(y^1,y^2)\in\R^{n-\sigma-1}\times\R^\sigma$,
\begin{equation*}
    \phi_{\mix}(x,t;y)
    =\phi_{\hyp}(x^1,t,y^1)+\phi_{\elliptic}(x^2,t,y^2)
\end{equation*}
and $\phi_{\mix}$ has signature $\sigma$. In addition, let $a_{\mix}(x,t;y) \coloneq a_{\hyp}(x^1,t,y^1) a_{\elliptic}(x^2,t,y^2)$.
We let $H_{\mix}^{\lambda}$ denote the H\"ormander-type operator associated to the pair $(\phi_{\mix},a_{\mix})$. The operators $H_{\hyp}^{\lambda}$ and $H_{\elliptic}^{\lambda}$ are associated to $(\phi_{\hyp}, a_{\hyp})$ and $(\phi_{\elliptic}, a_{\elliptic})$ in dimensions $n-\sigma$ and $\sigma+1$, respectively.

\begin{lemma}\label{lem: off diag sharp}
Let $n\geq 3$ and $\sigma\in\{0,1,\ldots,n-1\}$ be such that $n-\sigma$ is odd.
If $\bdH_{p\to q}(\phi_{\mathrm{mix}})$ holds, then $p$ and $q$ satisfy \eqref{eq: universal range off diagonal}.
\end{lemma}

\begin{proof}
Let $h_{\circ}\in\cS(\R^{n-\sigma-1})$ be defined via the Fourier transform by
\[
\widehat{h}_{\circ}(y^1)
\coloneqq 
e^{-i\lambda|y^1|^2/2}\psi(y^1),\quad y^1\in\R^{n-\sigma-1},
\]
for a nonnegative smooth cutoff $\psi\in C_c^{\infty}(\B^{n-\sigma-1})$ satisfying $\psi(0) = 1$.
Define
\begin{equation*}
f_\lambda(y^1,y^2)\coloneq\widehat{h}_\circ(y^1)\psi_\lambda(y^2) \qquad \textrm{for $y = (y^1, y^2) \in \R^{n-\sigma-1} \times \R^{\sigma}$,}
\end{equation*}
where $\psi_\lambda(y^2)=\psi(\lambda^{1/2}y^2)$ for $\psi\in C^{\infty}_c(\R^{\sigma})$ nonnegative, supported on $\B^{\sigma}$ with $\psi(0)=1$. It is easy to see that 
\begin{equation}\label{eq: off diag sharp 1}
\|f_\lambda\|_{L^p(\R^{n-1})}\sim\lambda^{-\frac{\sigma}{2p}}.
\end{equation}
By Fubini's theorem,
\begin{equation*}
  H_{\mix}^\lambda f_\lambda(x,t)=H_{\hyp}^\lambda\widehat{h}_\circ(x^1,t)\cdot H_{\elliptic}^\lambda\psi_\lambda(x^2,t),
\end{equation*}
and for $q<\infty$,
\begin{equation*}
    \|H_{\mix}^\lambda f_\lambda\|_{L^q(\R^n)}^q
    =\int_{\R} \|H_{\hyp}^\lambda \widehat h_\circ\|_{L^q(\R^{n-\sigma-1}\times\{t\})}^q
    \cdot\|H_{\elliptic}^\lambda\psi_\lambda\|_{L^q(\R^\sigma\times\{t\})}^q\,\ud t.
\end{equation*}
Consider the variety
\begin{equation*}
\cZ_{n-\sigma}
\coloneq
\big\{(x,t) \in \R^{n-\sigma-1} \times \R : x_{2j-1}=-t\cdot x_{2j},\quad j=1,\ldots,\big\lfloor\tfrac{n-\sigma-1}{2}\big\rfloor\big\}    
\end{equation*}
and its $\lambda$-scaling
\begin{equation*}
    \cZ_{n-\sigma}^\lambda\coloneqq
    \{(\tilde{x},\tilde{t})\in\R^{n-\sigma-1}\times\R:
      (\tilde{x}/\lambda,\tilde{t}/\lambda)\in\cZ_{n-\sigma}\}.
\end{equation*}
Then there exists a constant $c>0$ such that
\begin{equation}\label{eq: hyp example 1}
|H^\lambda_{\hyp}\widehat h_{\circ}(\bdx)|
\gtrsim 
\lambda^{-\frac{n-\sigma-1}{4}}
\quad 
\textrm{for all $\bdx\in N_c \cZ_{n-\sigma}^{\lambda}\cap B(0,c\lambda)$,}
\end{equation}
where $N_c \cZ_{n-\sigma}^{\lambda}$ denotes the $c$-neighbourhood of $\cZ_{n-\sigma}^{\lambda}$. For the full details, see \cite[pp. 1446--1448]{CGGHIW}.
For $|t|\leq c\lambda/2$ consider the slices
\begin{equation*}
  \cZ_{n-\sigma}^{\lambda}[t]
\coloneq 
\cZ_{n-\sigma}^{\lambda}\cap (\R^{n-\sigma-1}\times\{t\}).  
\end{equation*}
Each of these sets is a $\frac{n-\sigma-1}{2}$-dimensional affine space contained in $\cZ_{n-\sigma}^{\lambda}$. Using the bound \eqref{eq: hyp example 1}, we have
\begin{equation*}
  \|H_{\hyp}^{\lambda}\widehat h_{\circ}\|_{L^q(\R^{n-\sigma-1}\times\{t\})}
\gtrsim 
\lambda^{-\frac{n-\sigma-1}{4}}|N_c \cZ_{n-\sigma}^{\lambda}[t]\cap B(0,c\lambda)|^{1/q}
\sim
\lambda^{-\frac{n-\sigma-1}{2}(\frac{1}{2} - \frac{1}{q})}.  
\end{equation*}
For the elliptic term, we have
\begin{equation*}
 H_{\elliptic}^\lambda\psi_\lambda(x^2,t)=\lambda^{-\sigma/2}
 \int_{\B^\sigma} e^{i(\lambda^{-1/2}\langle x^2,\eta\rangle+
 \frac12\langle(\bdB^{\sigma+1}_{\elliptic}(t/\lambda)+\bdC^{\sigma+1}_{\elliptic})\eta,\eta\rangle)}
 a_{\elliptic}^{\lambda}(x^2,t;\lambda^{-1/2}\eta)\psi(\eta)\,\ud\eta.
\end{equation*}
There exists some absolute constant $c>0$ such that
\begin{equation*}
  |\bdB^{\sigma+1}_{\mathrm{ell}}(t/\lambda)\eta|\leq\|\bdB^{\sigma+1}_{\mathrm{ell}}(t/\lambda)\||\eta|\leq 10^{-2}  \quad \textrm{and} \quad \lambda^{-1/2}|\langle x^2,\eta\rangle|\leq\lambda^{-1/2}|x^2||\eta|\leq 10^{-2}
\end{equation*}
and therefore
\begin{equation*}
    \big|\lambda^{-1/2}\langle x^2,\eta\rangle+\frac12\langle(\bdB^{\sigma+1}_{\mathrm{ell}}(t/\lambda)+\bdC^{\sigma+1}_{\mathrm{ell}})\eta,\eta\rangle\big| \leq \frac{1}{2} + \frac{2}{100} < \frac{\pi}{2}
\end{equation*}
 for all $|t| \leq c \lambda$, $|x^2| \leq c \lambda^{1/2}$ and all $\eta \in \B^{\sigma}$. It follows that
\begin{equation*}
    \|H_{\elliptic}^\lambda\psi_\lambda\|_{L^q(\R^\sigma\times\{t\})}\gtrsim \lambda^{-\frac{\sigma}{2}(1-1/q)} \qquad \textrm{for all $|t|\leq c\lambda$.}
\end{equation*}
Combining the lower bounds for the two terms, and recalling \eqref{eq: off diag sharp 1}, we have
\begin{equation*}
    \|H_{\mix}^\lambda f_\lambda\|_{L^q(\R^n)}\gtrsim \lambda^{-\frac{n+\sigma-1}{4}+\frac{n+1}{2q}} \quad\text{so}\quad \frac{\|H_{\mix}^\lambda f_\lambda\|_{L^q(\R^n)}}{\|f_\lambda\|_{L^p(\R^{n-1})}}\gtrsim \lambda^{-\frac{n+\sigma-1}{4}+\frac{n+1}{2q}+\frac{\sigma}{2p}}.
\end{equation*}
Thus, $\bdH_{p\to q}(\phi_{\mathrm{mix}})$ can hold only if the first condition in \eqref{eq: universal range off diagonal} is satisfied.

On the other hand, the necessity of the second condition comes from considering a standard Knapp example. We test against $g_\lambda(y)\coloneq\psi(\lambda^{1/2}y)$ for $\psi\in C^{\infty}_c(\R^{n-1})$ nonnegative, supported on $\B^{n-1}$ with $\psi(0)=1$. Repeating the calculation involving $\psi_\lambda$ above, we get that $$
\frac{\|H^\lambda_{\mathrm{mix}} g_\lambda\|_q}{\|g_\lambda\|_p}
\gtrsim \lambda^{\frac{n+1}{2q}-\frac{n-1}{2p'}},$$ which forces
$\tfrac{n+1}{q} \leq \tfrac{n-1}{p'}$.

\end{proof}






\appendix




\section{Tube intersection and universal Kakeya estimates}\label{sec: tube intersect}

Let $\phi \colon \D^n \to \R$ be a H\"ormander-type phase and $0 < \delta < 1$. We consider $\delta$-tubes associated to $\phi$, which are sets of the form
\begin{equation*}
    T = \big\{(x,t) \in \B^{n-1} \times \B^1 : |\partial_y \phi(x, t; y) - \omega| < \delta\big\}.
\end{equation*}
For $T$ as above, we write $\omega(T) \coloneq \omega$ and $y(T) \coloneq y$. Denote $$J_\phi(\bdx;y_0)\coloneq\partial^2_{yy}\langle\partial_\bdx\phi(\bdx;y),G(\bdx,y_0)\rangle\vert_{y=y_0}.$$
The following tube intersection lemma is a generalization of \cite[Lemma 20]{WisewellThesis}.
\begin{lemma}\label{lem: tube intersect} Let $A$, $B \geq 1$ and $\phi \colon \D^n \to \R$ be a H\"ormander-type phase satisfying 
    \begin{enumerate}[P1)]
    \item $\|\partial_\bdx\phi\|_{C^3(\D^n)},\|\partial^2_{yy}\phi\|_{L^{\infty}(\D^n)}\leq A$, 
    \item $\|(\partial^2_{x y}\phi)^{-1}\|_{L^{\infty}(\D^n)},\|J_\phi\|_{L^{\infty}(\D^n)},\|J_\phi^{-1}\|_\infty\leq B.$
\end{enumerate}
There exists $\rho=\rho_{A,B}\in(0,1)$ for which the following holds. Let $T_1$, $T_2$ be two $\delta$-tubes associated to $\phi$. Then
\begin{equation}\label{eq: tube intersection}
    |y(T_2) - y(T_1)| + |\omega(T_2) - \omega(T_1)| \leq C_{A,B} \frac{\delta}{\delta+\diam{T_1\cap T_2}},
\end{equation}
provided $|y(T_2)-y(T_1)|,\,\diam{T_1\cap T_2}\leq\rho$. 
\end{lemma}
\begin{proof}
    Let $y_i \coloneq y(T_i)$, $\omega_i \coloneq \omega(T_i)$, and $\Delta \coloneq \diam{T_1\cap T_2}$. We focus only on the non-trivial case $\Delta\geq cB\delta$, for an appropriately large constant $c$, possibly depending on $n$. Indeed, for $\Delta<cB\delta$, the claim follows easily. 
    
    There exists $\bdx_1,\bdx_2\in\B^{n-1} \times \B^1$ with $|\bdx_2-\bdx_1|\sim\Delta$, satisfying 
    \begin{equation}
        \label{eq: tube intersect 1}
        |\partial_y\phi(\bdx_j;y_i)-\omega_i|<\delta,\quad i,j\in\{1,2\}.
    \end{equation}
    Differencing \eqref{eq: tube intersect 1} in $i$, we get
    \begin{equation}
        \label{eq: tube intersect 2}
        |(\partial_y\phi(\bdx_j;y_2)-\partial_y\phi(\bdx_j;y_1))-(\omega_2-\omega_1)|<2\delta,
    \end{equation}
    Differencing \eqref{eq: tube intersect 2} in $j$, we get
   \begin{equation*}
    |(\partial_y\phi(\bdx_2;y_2)-\partial_y\phi(\bdx_2;y_1))-(\partial_y\phi(\bdx_1;y_2)-\partial_y\phi(\bdx_1;y_1))|<4\delta.
    \end{equation*}
    Applying the multivariate Taylor theorem and P1),
    \begin{equation*}
        \big| \big( (\partial_{y \bdx}^2 \phi)(\bdx_1; y_2) - (\partial_{y \bdx}^2 \phi)(\bdx_1; y_1)\big)(\bdx_2 - \bdx_1) \big| \lesssim \delta + A|y_1 - y_2||\bdx_1 - \bdx_2|^2,
    \end{equation*}
    where we have also applied the mean value theorem in the $y$ variable to estimate the remainder. It is convenient to normalise the vector $\bdx_2 - \bdx_1$, giving
    \begin{equation}\label{eq: tube intersect 3}
        \Big| \big( (\partial_{y \bdx}^2 \phi)(\bdx_1; y_2) - (\partial_{y \bdx}^2 \phi)(\bdx_1; y_1)\big)\frac{\bdx_2 - \bdx_1}{|\bdx_2 - \bdx_1|}\Big|  \lesssim \frac{\delta}{\Delta} + A\Delta |y_1 - y_2|.
    \end{equation}
    We now apply the mean value theorem in the $y$ variable to the difference on the left-hand side of \eqref{eq: tube intersect 3}. After changing order of derivatives and exploiting linearity, this yields
    \begin{equation}\label{eq: tube intersect 4}
       \Big| \partial_{yy}^2 \inn{(\partial_{\bdx} \phi)(\bdx_1; y)}{\frac{\bdx_2 - \bdx_1}{|\bdx_2 - \bdx_1|}}\Big|_{y = y_1}(y_2 - y_1)\Big|  \lesssim \frac{\delta}{\Delta} + A|y_1 - y_2|\big(|y_1 - y_2| + \Delta\big).
    \end{equation}

\medskip\noindent\textbf{Claim:} Let $\mathrm{proj}_{(\bdx_2 - \bdx_1)^{\perp}}$ denote the orthogonal projection onto the hyperplane orthogonal to $\bdx_2 - \bdx_1$. Then
\begin{equation}\label{eq: tube intersect 5}
   \Big| \mathrm{proj}_{(\bdx_2 - \bdx_1)^{\perp}} G(\bdx_1; y_1)\Big| \lesssim_{A,B} \frac{\delta}{\Delta} + \Delta.
\end{equation}

We temporarily assume the claim. Form the orthogonal decomposition of $G(\bdx_1; y_1)$ with respect to the subspace spanned by $\bdx_2 - \bdx_1$ and its complement. Combining \eqref{eq: tube intersect 4} and \eqref{eq: tube intersect 5} we find 
\begin{equation*}
\big|\partial_{yy}^2\langle\partial_{\bdx}\phi(\bdx_1;y),G(\bdx_1;y_1)\rangle\vert_{y=y_1}(y_2-y_1)\big| \lesssim_{A,B} \frac{\delta}{\Delta} + |y_1 - y_2|\big(|y_1 - y_2| + \Delta\big).
\end{equation*}
Using the hypothesis $|y_1 - y_2| + \Delta \leq 2\rho$, this further implies
\begin{equation}\label{eq: tube intersect 6}
\big|\partial_{yy}^2\langle\partial_{\bdx}\phi(\bdx_1;y),G(\bdx_1;y_1)\rangle\vert_{y=y_1}(y_2-y_1)\big| \lesssim_{A,B} \frac{\delta}{\Delta}+ \rho|y_2-y_1|.
\end{equation}

By hypothesis P2), we have
\begin{equation}\label{eq: tube intersect 7}
\big|\partial_{yy}^2\langle\partial_{\bdx}\phi(\bdx_1;y),G(\bdx_1;y_1)\rangle\vert_{y=y_1}(y_2-y_1)\big|\gtrsim B^{-1}|y_2-y_1|.
\end{equation}

Combining \eqref{eq: tube intersect 6} and \eqref{eq: tube intersect 7}, there exists a constant $C_{A,B}\geq 1$ such that
\begin{equation*}
    |y_2-y_1|\leq C_{A,B}\frac{\delta}{\Delta}+C_{A,B}\rho|y_2-y_1|.
\end{equation*}
If we choose $0 <\rho\leq (2C_{A,B})^{-1}$, then $|y_2-y_1|\lesssim\frac\delta\Delta$. On the other hand, applying the triangle and mean value theorem to \eqref{eq: tube intersect 2} gives
\begin{equation*}
    |\omega_2-\omega_1| \lesssim\delta+|\partial_y\phi(\bdx_j;y_2)-\partial_y\phi(\bdx_j;y_1)| \lesssim\delta+|y_2-y_1|\|\partial^2_{yy}\phi\|_\infty.
\end{equation*}
Applying P1) to control $\|\partial^2_{yy}\phi\|_\infty$, and using the bound on $|y_2-y_1|$ concludes the argument.

Now we prove the claim. By the implicit function theorem, provided we are sufficiently localised, we can ensure that the tube $T_1$ admits a core curve. That is, $T_1$ is contained in an $O(\delta)$ neighbourhood of some smooth curve $\Gamma_1$. This curve $\Gamma_1$ can be parametrised as $\Gamma(t) = (\gamma(t),t)$ for a smooth function $\gamma$ satisfying 
\begin{equation*}
(\partial_{y_j}\partial_{\bdx}\phi)(\Gamma(t);y_1) \cdot
\Gamma'(t)
= 0
 \qquad \textrm{for $1 \leq j \leq n-1$.}
\end{equation*}
Write $\bdx_j=(x_j,t_j)$ and pick a point $\bdz_j=(z_j,t_j)$ from the core curve with the same height. It follows that 
$$|\partial_y\phi(x_j,t_j;y_1)-\partial_y\phi(z_j,t_j;y_1)|<\delta.$$
Applying the fundamental theorem of calculus and using P2), we find that 
$$|\bdz_j-\bdx_j|\lesssim B\delta,\quad j=1,2.$$ 
Since $B\delta\lesssim\Delta$, we have $$\big||\bdz_2-\bdz_1|-|\bdx_2-\bdx_1|\big|\lesssim_n B\delta\implies |\bdz_2-\bdz_1|\sim \Delta.$$ It follows that $$\Big|\mathrm{proj}_{(\bdx_2 - \bdx_1)^{\perp}}u-\mathrm{proj}_{(\bdz_2 - \bdz_1)^{\perp}}u\Big|\lesssim\frac{B\delta}{\Delta}|u|.$$
On the other hand, P1) and P2) give
$$|G(\bdx_1;y_1)-G(\bdz_1;y_1)|\lesssim_{A,B}|\bdx_j-\bdz_j|\lesssim B\delta.$$
Thus, it suffices to prove 
\begin{equation*}
   \big| \mathrm{proj}_{(\bdz_2 - \bdz_1)^{\perp}} G(\bdz_1; y_1)\big| \lesssim_{A,B} \Delta.
\end{equation*}
Note that since $\bdz_j\in \Gamma_1$, there exist $t_j \in \R$ such that $\bdz_j=\Gamma(t_j)$  for $j = 1$, $2$. 

Applying Taylor's theorem component-wise, 
\begin{equation}\label{eq: tube intersect 9}
    \bdz_2 - \bdz_1 = \Gamma(t_2) - \Gamma(t_1) = (t_2-t_1)\Gamma'(t_1)+O(A^3B^3|t_2-t_1|^2),
\end{equation}
Now \eqref{eq: tube intersect 9} can be rearranged to give
\begin{equation*}
   \Gamma'(t_1) = \frac{|\bdz_2 - \bdz_1|}{t_2 - t_1} \cdot \frac{\bdz_2 - \bdz_1}{|\bdz_2 - \bdz_1|} + O(A^3B^3\Delta).
\end{equation*}
Since $|\Gamma'(t_1)|\geq 1$, projecting $\Gamma'(t_1)$ onto $(\bdz_2-\bdz_1)^\perp$ gives $$ \Big| \mathrm{proj}_{(\bdz_2 - \bdz_1)^{\perp}} G(\bdz_1; y_1)\Big| \lesssim A^3B^3\Delta,$$
noting that $\Gamma'(t_1)$ is parallel to $G(\bdz_1;y_1)$.
\end{proof}

Using the above lemma, we can prove the quantitative universal Kakeya estimate from Lemma~\ref{lemma: universal kakeya}.

\begin{proof}[Proof (of Lemma~\ref{lemma: universal kakeya})] Let $A$, $B \geq 1$ and suppose $\phi \colon \D^n \to \R$ satisfies the hypotheses of Lemma~\ref{lem: tube intersect}. It suffices to show that the conclusion of Lemma~\ref{lemma: universal kakeya} holds whenever the tubes $\T$ are associated to $\phi$. Indeed, we claim that for any translation-invariant phase $\phi$ and $K \geq 1$ as in Lemma~\ref{lemma: universal kakeya}, the phase $\phi^K$ satisfies P1) and P2) for some constants $A$, $B \geq 1$ independent of $K$. This largely follows from the discussion at the beginning of \S\ref{subsec: small R}. We only note here that, for P2), $G_{\phi^K}(\bdx;y)=G_{\phi}(\bdx/K;y)$ and therefore
\begin{align*}
    J_{\phi^K}(\bdx;y_0)=\partial^2_{yy}\langle \partial_{\bdx}\phi(\bdx/K;y),G_{\phi}(\bdx/K;y_0)\rangle\vert_{y=y_0}=J_{\phi}(\bdx/K;y_0).
\end{align*}

    Let $\rho=\rho_{A,B}$ be as in the conclusion of Lemma~\ref{lem: tube intersect}. We may restrict both the $L^s$ norm and the directions $y$ to balls of radius $\rho_{A,B}$, so that \eqref{eq: tube intersection} holds. Tile the spatial domain by a grid of $\delta$-cubes $\cP$. For $1\leq r < \#\T$ dyadic and $q\in\cP$, denote $$\T(q)\coloneq\{T\in\T:q\cap T\neq\emptyset\},\quad\cP_r\coloneq\{q\in\cP:r \leq \#\T(q) \leq 2r\}.$$
    Then 
    \begin{equation}
        \label{eq: rich cubes}
    \big\|\sum_{T\in\T}\chi_T\big\|_{L^s(B_{\rho})}^s\sim \delta^n \sum_{\substack{1 \leq r < \#\T \\ \mathrm{dyadic}}}r^s\#\cP_r,
    \end{equation} where the sum is over all dyadic numbers $1\leq r\leq\#\T$. 

    \noindent\textbf{Claim:} $\displaystyle\#\cP_r\lesssim_{A,B}\delta^{-1}\frac{\#\T^{n/(n-1)}}{r^{(n+1)/(n-1)}}$.
    
    Assuming the claim we can conclude the proof as follows. Taking geometric mean between the claimed incidence bound and the trivial bound $\#\cP_r\lesssim\delta^{-1}\frac{\#\T}{r}$, we obtain
    \begin{equation*}
\#\cP_r\lesssim\delta^{-1}\frac{\#\T^{(s+1)/2}}{r^s},\qquad 1\leq s\leq\frac{n+1}{n-1}. 
    \end{equation*}
    We use the resulting bound in \eqref{eq: rich cubes}. The dyadic summation introduces a logarithmic loss. Taking $s$ root establishes \eqref{eq: universal kakeya}. 

    Now we prove the claim. Consider the following subfamily of large density $$\T_1\coloneq\{T\in\T:\#(\cP_r\cap T)\geq \lambda/\delta\},\quad\lambda\coloneq\frac{\delta r\#\cP_r}{2\#\T}.$$ The incidences between $\cP_r$ and $\T_2\coloneq\T\setminus\T_1$ satisfy 
    \begin{equation*}
       \cI(\cP_r,\T_2) \coloneq\sum_{T\in\T_2}\#(\cP_r\cap T)\leq r\#\cP_r/2, 
    \end{equation*}
    and so 
    \begin{equation}
        \label{eq: T_1 incidences}
        r\#\cP_r\leq 2\cI(\cP_r,\T_1) =2\sum_{q\in\cP_r}\#\T_1(q),
    \end{equation}
     where $$\T_1(q)\coloneq\{T\in\T_1:q\cap T\neq\emptyset\}.$$ By \eqref{eq: T_1 incidences}, there exists $q^*\in\cP_r$ satisfying $\#\T_1(q^*)\geq r/2$. In view of Lemma \ref{lem: tube intersect}, if two tubes are directionally $d$-separated, then their intersection is contained in a ball of radius $C_{A,B}\delta/d$. Take $d\coloneq2C_{A,B}\delta/\lambda$, and split $\T_1(q^*)$ into at most $(d/\delta)^{n-1}$ subfamilies of directionally $d$-separated tubes. One of those subfamilies satisfies $$\#\T_1'(q^*)\geq (\delta/d)^{n-1}r/2.$$ Using this, the disjointedness of these tubes outside $B(q^*,\lambda/2)$, and the density of the tubes in $\T_1$, we get
     \begin{equation*}
         \Big|\bigcup_{T\in\T_1'(q^*)}\cP_r\cap T\setminus B(q^*,\lambda/2)\Big|=\sum_{{T\in\T_1'(q^*)}}|\cP_r\cap T\setminus B(q^*,\lambda/2)|\gtrsim_{A,B}r\lambda^n\delta^{n-1},
     \end{equation*}
     so that, $\#\cP_r\gtrsim_{A,B} \delta^{-1}r\lambda^n.$ Plugging the value of $\lambda$ and rearranging the resulting bound proves the claim.
\end{proof}

\bibliographystyle{alpha}
\bibliography{references}

\end{document}